\documentclass[11pt]{article}
\usepackage{amsmath,amssymb,amsfonts,amsthm,amscd,graphicx,psfrag,epsfig}
\usepackage{hyperref}
\hypersetup{colorlinks, linkcolor=blue}
\usepackage[all]{xy}           
\usepackage{marginnote}  
\usepackage{float}
\usepackage{appendix}
\usepackage{rotating}
\usepackage{multirow}
\usepackage{subcaption}
\newtheorem{theorem}{Theorem}[section]

\newtheorem{theorem-definition}[theorem]{Theorem-Definition}
\newtheorem{theorem-construction}[theorem]{Theorem-Construction}
\newtheorem{lemma-definition}[theorem]{Lemma--Definition}
\newtheorem{lemma-construction}[theorem]{Lemma--Construction}
\newtheorem{lemma}[theorem]{Lemma}
\newtheorem{proposition}[theorem]{Proposition}
\newtheorem{corollary}[theorem]{Corollary}
\newtheorem{conjecture}[theorem]{Conjecture}
\newtheorem{definition}[theorem]{Definition}

\newcommand{\old}[1]{}

\newcommand{\Z}{{\mathbb Z}}
\renewcommand{\S}{{\mathbb S}}
\newcommand{\R}{{\mathbb R}}
\newcommand{\Q}{{\mathbb Q}}
\newcommand{\C}{{\mathbb C}}

\newcommand{\lms}{\longmapsto}
\newcommand{\lra}{\longrightarrow}
\newcommand{\hra}{\hookrightarrow}

\newcommand{\be}{\begin{equation}}
\newcommand{\ee}{\end{equation}}
\newcommand{\bt}{\begin{theorem}}

\newcommand{\et}{\end{theorem}}
\newcommand{\bd}{\begin{definition}}
\newcommand{\ed}{\end{definition}}
\newcommand{\bp}{\begin{proposition}}
\newcommand{\ep}{\end{proposition}}

\newcommand{\bl}{\begin{lemma}}
\newcommand{\el}{\end{lemma}}
\newcommand{\bc}{\begin{corollary}}
\newcommand{\ec}{\end{corollary}}
\newcommand{\bcon}{\begin{conjecture}}
\newcommand{\econ}{\end{conjecture}}
\newcommand{\la}{\label}

\begin{document}
\date{}

\title{Ideal webs, moduli spaces of local systems, and 3d Calabi-Yau categories}
\author{A.B. Goncharov}

\maketitle

\tableofcontents 

\begin{abstract}

A {\it decorated surface} $S$ is an oriented surface with {punctures}, and  
a finite  set of {marked points} on the boundary, considered  modulo isotopy. We assume that 
each boundary component has a marked point. 
We introduce {\it ideal bipartite graphs}   on $S$.  
Each of them is related to a group $G$ of type  ${\rm A}_{m}$ or $GL_m$, and 
 gives rise to cluster coordinate systems  
on certain  moduli spaces of $G$-local systems on $S$. 
These coordinate systems  
generalize the ones assigned in  \cite{FG1} to ideal triangulations of $S$. 

A bipartite graph ${W}$  on $S$ gives rise to 
 a quiver with a canonical potential. The latter determines a triangulated 3d  
Calabi-Yau $A_\infty$-category ${\cal C}_{W}$  
with a {\it cluster collection ${\cal S}_W$} -- a generating collection of spherical 
objects of special kind \cite{KS1}. 

Let ${W}$ be an ideal bipartite graph on $S$ of type $G$. 
We define an extension  $\Gamma_{G, S}$  of the mapping class group of $S$, and  prove that it acts by 
symmetries of the category ${\cal C}_{W}$. 

There is a family of open CY threefolds over the {\it universal Hitchin base} ${\cal B}_{G, S}$, whose 
intermidiate Jacobians 
describe Hitchin's integrable system \cite{DDDHP}, \cite{DDP}, \cite{G}, 
\cite{KS3}.  
We conjecture that the 3d CY category with  cluster collection $({\cal C}_{W}, {\cal S}_{W})$ 
is equivalent to a full subcategory of the Fukaya category of a generic threefold of the family, 
equipped with a cluster collection of special Lagrangian spheres. 
For $G=SL_2$ a substantial part of the story
 is already known thanks to Bridgeland, Keller, Labardini-Fragoso, Nagao, Smith, and others, see \cite{BrS}, \cite{S}.  

We hope that ideal bipartite graphs
 provide special examples of the Gaiotto-Moore-Neitzke spectral networks \cite{GMN4}.
\end{abstract}

\section{Introduction and main constructions and results}

\subsection{Ideal webs and cluster coordinates on moduli spaces of local systems}  
Recall that a {\it decorated surface} $S$ is an oriented topological surface with {\it punctures} and  {\it marked points} on the boundary, considered  modulo isotopy.  
The set of {\it special points} is the union of the set of punctures and marked points. 
We assume that each boundary component has a special point,  the  number $n$ of special points is $>0$, and  
 that  $S$ is hyperbolic: if the genus $g(S)=0$, then $n\geq 3$. 
Denote by $\Gamma_S$ the mapping class group of $S$. 
\begin{figure}[ht]
\centerline{\epsfbox{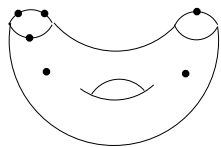}}
\caption{A decorated surface: punctures inside, special points on the boundary.}
\label{gra0}
\end{figure}

There is a pair of moduli spaces assigned to $S$, closely related to the moduli space of 
$m$-dimensional vector bundles with flat connections on $S$, defined  in \cite{FG1}:
\be \la{8.15.11.1}
\mbox{a $K_2$-moduli space ${\cal A}_{SL_{m}, S}$} ~~ \mbox{and} ~~ \mbox{a Poisson 
moduli space ${\cal X}_{PGL_{m}, S}$}.
\ee 

As the names suggest, the  space ${\cal A}_{SL_m, S}$ is equipped with a $\Gamma_S$-invariant class in 
$K_2$, and the space ${\cal X}_{PGL_m, S}$ is equipped with a $\Gamma_S$-invariant Poisson structure.\footnote{There are similar pairs of moduli spaces  for any split semisimple group $G$.}

An {\it ideal triangulation} of  $S$ is a triangulation  with vertices at the special points. 
In \cite{FG1} each of  the spaces (\ref{8.15.11.1}) was equipped 
with a collection of rational coordinate systems,  called {\it cluster coordinate systems},   
assigned to {ideal triangulations} of $S$. 
In particular this implies that these 
spaces  are rational.\footnote{For general $S$ the 
moduli space of $m$-dimensional vector bundles 
with flat connections on $S$ is not rational, and thus can not have a rational coordinate system.}  
The ring of regular functions on 
the space ${\cal A}_{SL_m, S}$ is closely related to 
(and often coincides with) a cluster algebra of Fomin-Zelevinsky \cite{FZI}.  

One of the crucial features
 of the cluster structure is that a single cluster coordinate system on a given space determines the cluster structure, and  
gives rise to a huge, usually infinite, 
collection of cluster coordinate systems, obtained by arbitrary sequences of {\it mutations}. 

When $G$ is of type ${\rm A}_1$, the coordinate systems related to ideal triangulations 
provide most of the cluster coordinate systems. Any two ideal triangulations are related by a sequence of flips. 
A flip of an ideal triangulation gives rise to a cluster mutation 
of the corresponding cluster coordinate system. Therefore any two of them are related by a sequence of 
mutations.

When $G$ is of type ${\rm A}_m$, $m>1$, the situation is dramatically different.  
First, the coordinate systems related to ideal triangulations form only a tiny 
part of all cluster coordinate systems. Second, 
although the coordinate systems 
assigned to any two ideal triangulations are still related by a sequence of mutations, 
the intermediate  cluster coordinate systems are not related to any 
ideal triangulation.

In this paper we introduce a wider and much more flexible class of 
cluster coordinate systems on the moduli spaces (\ref{8.15.11.1}). 
They are  
assigned to  
geometric data on $S$, which we call 
{\it ${\rm A}_{m}$-webs on $S$}.  We also introduce closely related and equally important 
{\it ${\rm A}_{m}^*$-webs on $S$}.

\begin{figure}[ht]
\centerline{\epsfbox{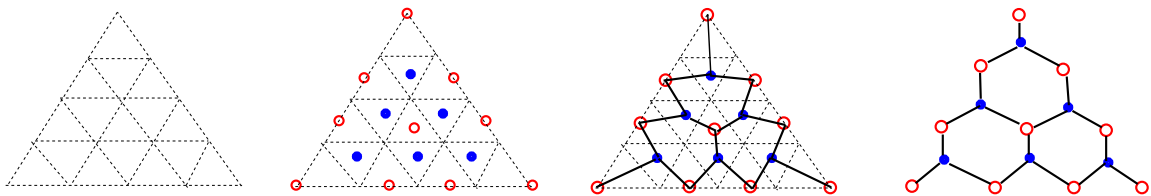}}
\caption{A bipartite graph associated with a 4-triangulation of a triangle.}
\label{gra10i}
\end{figure}

\paragraph{An example of an ${\rm A}_{m}$-web.} Let $T$ be an ideal triangulation of $S$. 
Given a triangle of $T$, we subdivide it into $(m+1)^2$ small triangles 
by drawing three families of $m$ equidistant lines, parallel to 
 the sides of the triangle, 
as shown on the left  of Figure. \ref{gra10i}.   

For every triangle  of $T$ there are two kinds of small triangles: 
the ``looking up'' and ``looking down'' triangles. We put a $\bullet$-vertex
 into the center of each of the 
``down'' triangles. 

Let us color  in red all ``up'' triangles.  
We consider the obtained red domains -- some of them are unions of red triangles, and 
put a $\circ$-vertex into the center of each of them. 
Alternatively, the set of $\circ$-vertices is described as follows. 
Put a $\circ$-vertex into every special point, into the centers of little segments on the sides of $t$ 
which do not contain the special points, and into the 
centers of interior ``up'' triangles, as shown on the middle picture on Figure \ref{gra10i}. 

A $\circ$-vertex and a $\bullet$-vertex 
 are {\it neighbors} if the corresponding domains share an edge. 

Recall that a {\it bipartite graph} is a graph with vertices of two kinds,  so that  
each edge connects vertices of different kinds.  
So connecting the neighbors,  
 we get a bipartite surface graph $\Gamma_{{\rm A}_{m}}(T)$, called the {\it bipartite ${\rm A}_{m}$-graph} assigned to $T$, see Figure \ref{gra10i}. 

\begin{figure}[ht]
\centerline{\epsfbox{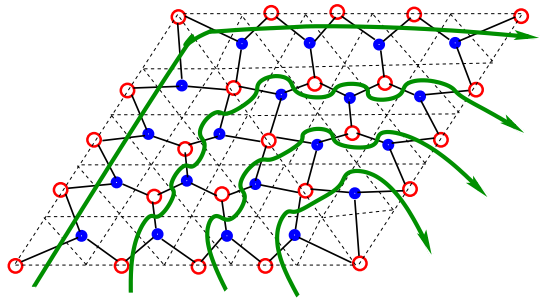}}
\caption{The strands of the  web assigned to
 a 5-triangulation of a pair of adjacent triangles, which go around one of the four special points - 
the bottom right vertex of the quadrilateral.}
\label{gra11i}
\end{figure}

We associate to any bipartite surface graph a collection of its {\it zig-zag strands}. Namely, take 
the paths on the graph which go left at the $\bullet$-vertices, and right at the $\circ$ ones, and 
push them a bit off the vertices, getting a (green) web shown on Figure \ref{gra11i}. 
The strands are oriented so that we go around 
$\bullet$-vertices clockwise. 
A zig-zag strand inside of a triangle $t$ 
goes along one of its sides.

The strands associated to the  bipartite graph $\Gamma_{{\rm A}_{m}}(T)$
enjoy the following properties: 

\begin{enumerate}

\item The directions of the strands intersecting any given strand alternate. 

\item The connected components of the complement to the union of strands are contractible. 

\item The part of surface on the right of each strand is a disc, containing  
 a \underline{single} special point. 

\end{enumerate}

\paragraph{Webs and ideal webs.} The alternation property 
is a well known feature of the collection of zig-zag strands 
 assigned to any bipartite graph on a surface. Adding  contractibility of the domains, and  
axiomatizing this, one defines a {\it web}
 on a surface as a connected collection of oriented strands in generic position, considered modulo isotopy,
 satisfying 
conditions (1) and (2). 
So the collection of zig-zag strands assigned to a bipartite graph form a web, and 
it is well known that this way we get a bijection
 between webs and bipartite graphs on surfaces whose faces are discs -- see \cite{GK} or Section \ref{ssec2.1}. 
$$
\mbox{So webs are just an alternative way to think about bipartite graphs on surfaces.} 
$$
The third condition is new. 
We define a {\it strict ideal web} on a decorated surface $S$ as a web 
satisfying  condition (3), plus a more technical minimality condition. 
An {\it ideal web} is a web which becomes strict ideal on the universal  cover of $S$. 
See  Section \ref{ssec2.3} for the precise definition. 

In this paper we refer to webs rather then to bipartite graphs since the crucial third condition 
is formulated using the zig-zag strands. 

\begin{figure}[ht]
\centerline{\epsfbox{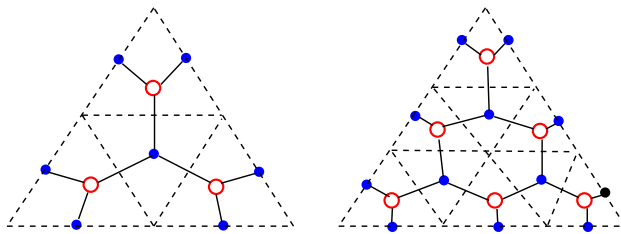}}
\caption{An ${\rm A}_{1}^*$-graph, and an ${\rm A}_{2}^*$-graph in an internal ideal triangle.}
\label{gc1}
\end{figure}

\begin{figure}[ht]
\centerline{\epsfbox{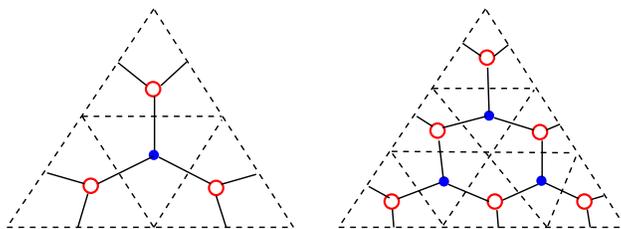}}
\caption{An ${\rm A}_{1}^*$-graph, and an ${\rm A}_{2}^*$-graph on a decorated surface given by a triangle.}
\label{gc1m}
\end{figure}

\paragraph{Another example: the ${\rm A}_{m}^*$-web assigned to ideal triangulations.}
Given an {$m+1$-triangulation} of $T$ of $S$,  
let us define the {\it bipartite ${\rm A}_{m}^*$-graph}  $\Gamma_{{\rm A}_{m}^*}(T)$  
 assigned to $T$, see Fig \ref{gc1}. 
Take a triangle $t$ of the triangulation $T$. 
Assign the $\circ$-vertices to the centers of all ``looking up'' triangles of the $(m+1)$-triangulation of $T$, 
and the $\bullet$-vertices to the centers of all ``looking down'' triangles, as well as to the 
 centers of $m+1$ little segments on 
each of the side of $t$. So every edge of the triangulation $T$ carries $m+1$ $\bullet$-vertices. 
Connect them into a bipartite graph as shown on Fig \ref{gc1}. 
Notice that the special points of $S$  are no longer among the vertices of the graph $\Gamma_{{\rm A}_{m}^*}(T)$. 
This procedure creates $\bullet$-vertices on the boundary of the surface. 
We remove the $\bullet$-coloring from all external vertices. For example, 
if the surface $S$ is just a triangle, we get a graph shown on Fig \ref{gc1m}

Shrinking the 2-valent $\bullet$-vertices of the bipartite graph $\Gamma_{{\rm A}_{m}^*}(T)$ we get the 
graph $\Gamma_{{\rm A}_{m}}(T)$. 

\paragraph{Elementary moves.} There are three 
 elementary local transformations of bipartite surface 
graphs / webs. The crucial one is the two by two move.

\begin{figure}[ht]
\centerline{\epsfbox{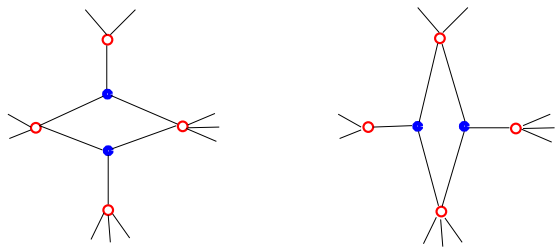}}
\caption{A two by two move. 
Flipping the colors of vertices delivers another two by two move. }
\label{gc2a}
\end{figure}

{\it Two by two move.} It is shown on Fig \ref{gc2a}. It affects  four
 vertices of one color, and two  
of the other. 
It amounts to a transformation of webs, 
also known as a two by two move, 
see Fig \ref{gra12}. 

\begin{figure}[ht]
\centerline{\epsfbox{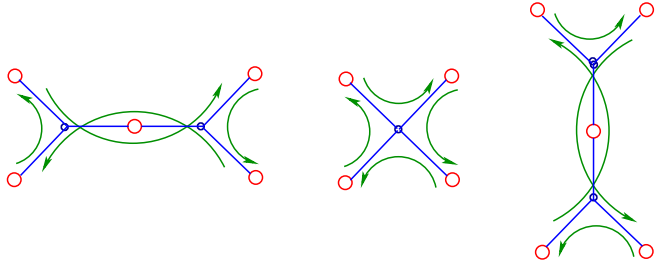}}
\caption{Shrinking a 2-valent $\circ$-vertex, and then expanding a 4-valent $\bullet$-vertex}
\label{bc04}
\end{figure}

{\it Shrinking a $2$-valent vertex.} It is the first move shown  on Fig \ref{bc04}. 
It amounts to resolving a bigon move of webs. 

{\it Exapanding a $k$-valent vertex, $k>3$.} It is the second move  on Fig \ref{bc04}. 
It amounts to creating a bigon move of webs. 
 
\paragraph{A guide to classes of webs 
studied in the paper.}

We define
several subclasses of webs:
\be \la{2.19.12.1}
\{\mbox{${\rm A}_{m}$ and ${\rm A}_{m}^*$-webs}\} \stackrel{}{=} \{\mbox{\it rank $m+1$ ideal webs}\}  \subset \{\mbox{\it spectral webs}\}.
\ee

\begin{itemize}

\item
{\it Spectral webs} are defined on any oriented surface $S$. 
A spectral web ${\cal W}$ on $S$ 
gives rise to a new surface $\Sigma_{\cal W}$, called the {\it spectral 
surface}, equipped with a ramified cover map  
$\Sigma_{\cal W} \lra S$, called the {\it spectral cover} -- see Section \ref{ssec2.2}. 
Its degree  is the {\it rank} of ${\cal W}$. 

\item
When $S$ is a decorated surface, we introduce a subclass 
of spectral webs, called {\it ideal webs}, which play the central role in our story-- 
see Section \ref{ssec2.3}.

\item 
There is an involution on the webs, whose orbits are {\it triple crossing diagrams}. 
 It interchanges the 
${\rm A}_{m}$-web with ${\rm A}_{m}^*$-webs.  -- see Section \ref{sec2.4.1}.
\end{itemize}

\paragraph{Cluster coordinate systems assigned to ${\rm A}_{m-1}$-webs and ${\rm A}^*_{m-1}$-webs.} 
An ${\rm A}_{m-1}$-web ${\cal W}$ on $S$ 
 gives rise to a cluster coordinate system 
on each of the spaces (\ref{8.15.11.1}). They
generalize the ones assigned in  \cite{FG1} to ideal triangulations of $S$. 
Two by two moves amounts to cluster mutations. 

Precisely, a face of a web is called {\it internal} if it does not intersect the boundary of $S$, and does not contain special points. Otherwise a face is called 
{\it external}. 
We define regular functions $\{A_{F}^{{\cal W}}\}$ on  the space 
${\cal A}_{SL_{m}, S}$, parametrised by the faces $F$ of the web ${\cal W}$. They form a cluster coordinate system. 
The functions $\{A_{F}^{{\cal W}}\}$ corresponding to the external faces of ${\cal W}$ 
are the frozen cluster coordinates: they do not mutate. 

We define rational functions $\{X_{F}^{{\cal W}}\}$ on the space ${\cal X}_{PGL_{m}, S}$  
parametrised by the internal faces $F$ of the web ${\cal W}$; they form a cluster Poisson 
coordinate system.

\vskip 3mm
An {\it decorated flag} in an $m$-dimensional vector space $V_m$ is a flag 
$F_0 \subset F_1 \subset \ldots \subset F_m$ where ${\rm dim}F_k=k$ equipped with 
non-zero vectors $v_i \in F_i/F_{i-1}-\{0\}$, $i=1, ..., m$.
Denote by ${\cal A}_{GL_m}$ the moduli space of all decorated flags 
in $V_m$. So ${\cal A}_{GL_m} = GL_m/U$ where $U$ is the upper triangular unipotent subgroup in $GL_m$.\footnote{There are two versions of decorated flags used in this paper: one of them includes the factor $v_m$ in the data, the other does not. The first is 
decorated flag, that is a point of ${\rm G}/{\rm U}$,  for ${\rm G}=GL_m$, the second for ${\rm G}=SL_m$.} 

Pick a volume form $\Omega_m$ in $V_m$. Consider the moduli space ${\rm Conf}_n({\cal A}^*_{SL_m})$ parametrising 
$n$-tuples of decorated flags in $V_m$ modulo the diagonal action of the group $SL_m = {\rm Aut}(V_m, \Omega_m)$:
\be \la{A*space}
{\rm Conf}_n({\cal A}^*_{SL_m}) = ({\cal A}_{GL_m} \times \ldots \times {\cal A}_{GL_m} )/{\rm Aut}(V_m, \Omega_m).
\ee

In the case when $S$  is a disc with $n$ special points on the boundary,   
an ${\rm A}_{m-1}^*$-web  ${\cal W}$ on $S$ describes a cluster coordinate system 
on the moduli space (\ref{A*space}). 
The external faces of ${\cal W}$ parametrise the frozen coordinates. 

For any decorated surface $S$ there is an analog of the  space (\ref{A*space}), 
denoted by ${\cal A}^*_{SL_m, S}$. It 
parametrises twisted $SL_m$-local systems on $S$ with an additional data: 
a choice of a decorated flag near each of the special points, invariant under the monodromy around the point. 
An ${\rm A}_{m-1}^*$-web  ${\cal W}$ on $S$ describes  a cluster coordinate system 
on the space ${\cal A}^*_{SL_m, S}$. Its external faces parametrise the frozen coordinates. 

\vskip 3mm
On the other hand, given a ${\rm A}_{m-1}^*$-web  ${\cal W}$ on $S$, we describe 
generic points of the spaces ${\cal A}_{GL_m, S}$ and ${\cal X}_{GL_m, S}$ 
via flat line bundles on the spectral surface $\Sigma_{S, \cal W}$. 
This description admits a generalisation to non-commutative local systems on $S$ \cite{GKo}.

\subsection{Presenting an ${\rm A}_{m}$-flip as a composition of two by two moves} \la{2by2moves}

Starting from now,  we consider only triangulations without self-folded triangles.  
 
We have seen that each ideal triangulation of $S$ gives rise to a canonical  ${\rm A}_{m}$-web. 
Our next goal is to connect the   ${\rm A}_{m}$-webs assigned to two 
different triangulations by a sequence of elementary moves. 
Any two ideal triangulations without self-folded triangles can be connected by a sequence of flips that only involves ideal triangulations without self-folded 
triangles. 
Therefore it suffices to show that every flip in the sequence is a cluster transformation. 
So it is sufficient to 
do this for a single flip. We say that the  ${\rm A}_{m}$-webs 
related to an original ideal triangulation and a flipped one are related by an 
{\it ${\rm A}_{m}$-flip}. 

\begin{figure}[ht]
\centerline{\epsfbox{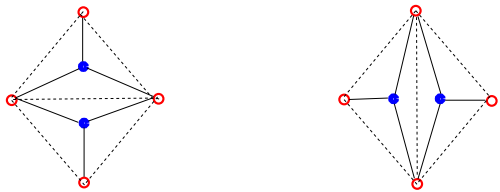}}
\caption{Two by two move corresponding to a flip of triangulations.}
\label{gra13a}
\end{figure}

Figure \ref{gra13a} tells that an ${\rm A}_{1}$-flip is the same thing as two by two move. 

Proposition \ref{slm} tells how to decompose an ${\rm A}_{m}$-flip for any $m$ 
as a composition of elementary moves. The existence of such a decomposition follows from  
Theorem \ref{8.13.11.1}. Here we  present an ${\rm A}_{m}$-flip as a \underline{specific}  
composition of $m+2 \choose 3$ two by two moves, and shrinking / expanding  moves.

\bp \la{slm}
There is a decomposition of an ${\rm A}_{m}$-flip into a composition of elementary moves, including 
$m+2 \choose 3$ two by two moves, 
depicted on Fig \ref{gra13a} for $m=1$, Fig \ref{bc1} for $m=2$, and on 
Fig \ref{bc0}-\ref{bc02} for $m=3$. 
\ep 

\begin{proof} The pattern is clear from the $m=3$ and $m=4$ examples. 

\begin{figure}[ht]
\centerline{\epsfbox{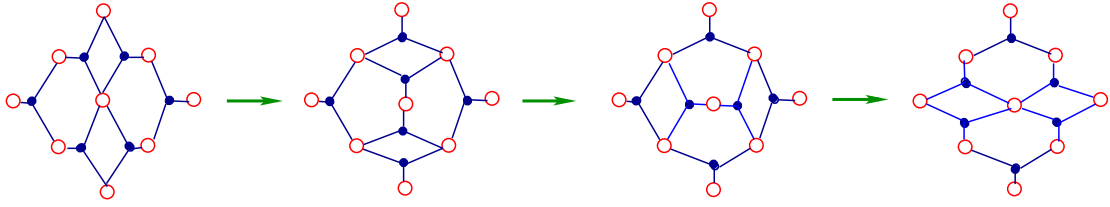}}
\caption{Presenting  an ${\rm A}_{2}$-flip as a sequence of elementary moves. We 
transform the ${\rm A_2}$-web related to one triangulation  
into the 
${\rm A_2}$-web related to  the 
other. We apply $2\times 1 = 2$ two by two moves, shrink and expand the central 
vertex, and 
apply $1\times 2 =2$ two by two moves.}
\label{bc1}
\end{figure}

\newpage 
Figures \ref{bc0}-\ref{bc02} show how to  
decompose a transformation of ${\rm A}_{3}$-webs arising from 
a flip of the vertical diagonal  into a composition of 
$3\times 1 + 2 \times 2 + 1\times 3 = 3+4+3=10$ two by two moves.

\begin{figure}[ht]
\centerline{\epsfbox{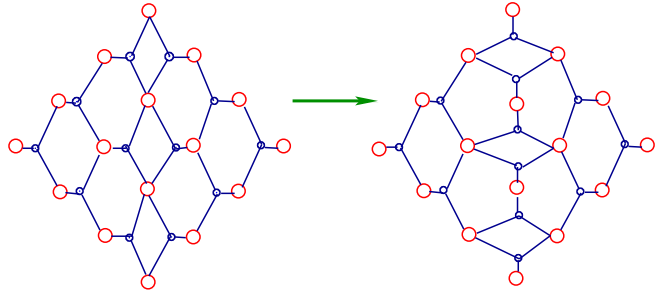}}
\caption{Step 1: Applying $3\times 1 =3$ two by two moves.}
\label{bc0}
\end{figure}

\noindent To go from Step 1 to Step 2 we shrink vertically and expand horisontally 
two  of the ``middle'' $\circ$-vertices. 

\begin{figure}[ht]
\centerline{\epsfbox{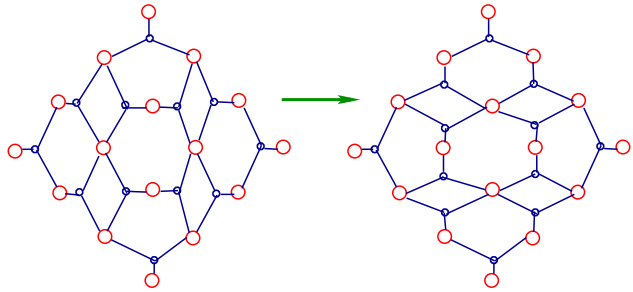}}
\caption{Step 2: Applying $2\times 2 =4$ two by two moves.}
\label{bc01}
\end{figure} 
\noindent To go from Step 2 to Step 3 we shrink vertically and expand horisontally 
the other two   ``middle'' $\circ$-vertices. 
\begin{figure}[ht]
\centerline{\epsfbox{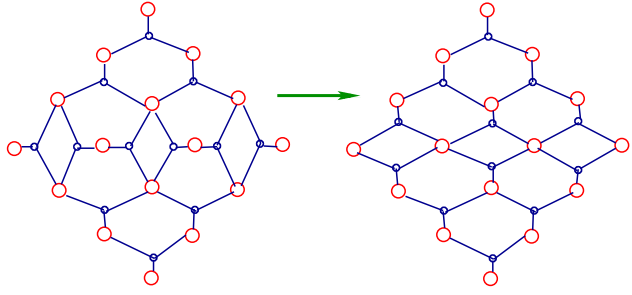}}
\caption{Step 3: Applying $1\times 3 =3$ two by two moves.}
\label{bc02}
\end{figure}

Similarly, a transformation of ${\rm A}_{m}$-webs arising from 
a flip of a triangulation is a composition of 
$m\times 1 + (m-1) \times 2 + \ldots + 1\times m = {m+2 \choose 3}$ two by two moves. 
Notice that the 
$\circ$-vertices do not change their position in the process. \end{proof}

{\bf Remark}. Given a positive integer $m$,  we assigned in \cite{FG1} to an ideal triangulation a quiver,  
 called the ${\rm A}_{m}$-quiver. A flip gives rise to a transformation of ${\rm A}_{m}$-quivers, 
presented in \cite{FG1}, Section 10, as an explicit composition 
 of cluster mutations.  
Proposition \ref{slm} tells that all intermediate quivers 
are described by  ${\rm A}_{m}$-webs, and the sequence 
of cluster mutations match the sequence of two by two moves. 
Each two by two move corresponds to gluing of an octahedron. 
This collection of two by two moves  is parametrised 
by the octaherda of the {\it ${\rm A}_{m}$-decomposition of a tetrahedra}, see  {\it loc. cit.}.


\subsection{${\rm CY}_3$ categories  categorifying cluster varieties} \la{anover}
It was discovered more then a decade ago that many interesting rational varieties which 
appear in representation theory, geometry  and physics 
come with a natural cluster structure of certain kind.  Examples are given by 
 Grassmannians \cite{GSV}, double Bruhat cells \cite{BFZ}, different versions of 
the classical Teichmuller space  \cite{GSV1},  \cite{FG1}, and more generally different 
versions of moduli spaces of $G$-local systems on decorated surfaces, 
e.g. (\ref{8.15.11.1}), \cite{FG1}.

So far cluster structures were used mostly to describe certain features 
of the relevant spaces, like Poisson / $K_2$-structures, 
canonical bases in the space of regular functions, 
non-commutative $q$-deformation, 
quantisation, etc. 

We suggest that {\it each cluster variety which appears in representation theory, geometry or physics 
admits a canonical categorification, given by  a 
3d Calabi-Yau category}.

To formulate this precisely, 
let us first review some  basic definitions. 

\paragraph{3d CY categories.} A 3d Calabi-Yau category over a characteristic zero field $k$ is a 
$k$-linear triangulated $A_\infty$-category ${\cal C}$ with a non-degenerate symmetric 
degree $-3$ pairing on ${\rm Hom}$'s: 
\be\la{4.18.12.1uh}
(\bullet, \bullet): {\rm Hom}_{\cal C}^\bullet(E, F) \otimes {\rm Hom}_{\cal C}^\bullet(F, E) \lra k[-3].
\ee

The  $A_\infty$-structure on a 3d CY category is encoded by the 
cyclic higher composition maps:
\be\la{4.18.12.2}
m_{E_1, ..., E_{n}}: \bigotimes_{1 \leq i \leq n}{\rm Hom}_{\cal C}^\bullet(E_i, E_{i+1}[1]) \lra k.
\ee
These maps are 
 defined for any  objects $E_1, \ldots , E_{n}$ of the category ${\cal C}$. They are graded 
 cyclically invariant, and satisfy certain quadratic relations which we omit.

An object $S$ of a category ${\cal C}$ is {\it spherical} if  
$
{\rm Ext}_{\cal C}^\bullet(S, S)= H^\bullet(S^3, k).
$ 
\bd \cite{KS1} A collection ${\cal S} = \{S_i\}_{i\in I}$ of spherical generators of a 3d CY category ${\cal C}$ 
is a {\rm cluster collection} if   
for any two different $i,j \in I$ one has 
\be \la{ccc}
 ~~~~\mbox{\rm ${\rm Ext}_{\cal C}^a(S_i, S_j)$ can be non 
zero only for a single $a\in \{1,2\}$}.   
\ee
\ed 

A 3d CY category ${\cal C}$ with a cluster collection ${\cal S}$ gives rise to 
a quiver $Q_{{\cal C}, {\cal S}}$. 
Its vertices are the isomorphism classes $[S_i]$ of objects of the cluster collection, 
and the number of arrows $[S_i] \to [S_j]$ equals  ${\rm dim}{\rm Ext}_{\cal C}^1(S_i, S_j)$.  
Thus the quiver $Q_{{\cal C}, {\cal S}}$ has no loops and length two cycles. 

Let us assume that $m_1=0$. Then the cyclic composition maps, restricted to the 
cyclic product of 
${\rm Ext}^1(E_i, E_{i+1})[1]$, determine a {\it potential} $P_{{\cal C}, {\cal S}}$ of the quiver 
$Q_{{\cal C}, {\cal S}}$, that is a linear functional on the vector spaces spanned by the cycles in the quiver.
 The original 3d CY category ${\cal C}$ with a cluster collection ${\cal S}$, considered 
modulo autoequivalences  preserving the cluster collection ${\cal S}$, 
can be reconstructed from the quiver with potential $(Q_{{\cal C}, {\cal S}}, P_{{\cal C}, {\cal S}})$ \cite{KS1}. 

A spherical object $S$ of a CY category 
 gives rise to the 
Seidel-Thomas  \cite{ST} reflection  functor ${\rm R}_{S}$, 
acting by an autoequivelence 
of the category:
$$
R_S(X):= {\rm Cone}\Bigl({\rm Ext}^\bullet(S, X) \otimes S \lra X\Bigr).
$$  

Given elements $0, i\in I$,  we  write $i<0$  if ${\rm Ext}^1(S_i, S_0)$ is non-zero,
 and $i>0$ otherwise.  
\bd
A mutation of a cluster collection ${\cal S} = \{S_i\}_{i\in I}$ in a 3d CY  category ${\cal C}$ at an object $S_0$ 
is a new spherical collection ${\cal S}' =\{S'_i\}$ in the same category ${\cal C}$ given by 
$$
S'_i=S_i, ~~i<0, ~~~~S'_0=S_0[-1], ~~~~S'_i=R_{S_0}(S_i), ~~i>0,
$$
\ed
Since the shift and reflection functors transform spherical objects to spherical ones, 
$\{S'_i\}_{i\in I}$ is a collection of spherical objects. Since the objects $\{S_i\}$ 
generate the triangulated category ${\cal C}$, 
the objects $\{S'_i\}$  are generators. 
However $\{S'_i\}_{i\in I}$ is not necessarily a cluster collection. 

 Mutations of cluster collections, being projected to 
$K_0({\cal C}) $, recover mutations of bases (\ref{12.12.04.2a}).

\vskip 3mm

Mutations of cluster collections  
translate into {\it mutations of quivers with generic potentials}  introduced earlier by
 Derksen, Weyman and Zelevinsky \cite{DWZ}. 

Given a 3d CY category with a cluster collection $({\cal C}, {\cal S})$, denote by 
${\rm Auteq}({\cal C}, {\cal S})$ the quotient of the group of autoequivelances of 
${\cal C}$ by the subgroup 
 preserving a cluster collection ${\cal S}$ in ${\cal C}$. 

\paragraph{Categorification of cluster varieties.} 
Returning to cluster varieties, we suggest the following metamathematical principle: 

\begin{itemize}

\item {\it Each cluster variety   ${\cal V}$
 which appears
in representation theory, geometry  and physics,  
gives rise to  a  combinatorially defined $3$-dimensional 
Calabi-Yau triangulated ${\rm A}_{\infty}$-category 
${\cal C}_{\cal V}$. 

The category ${\cal C}_{\cal V}$ comes 
with an array of cluster collections $\{{\cal S}_\alpha\}$, related by cluster mutations. 

There is a  group ${\Gamma}_{\cal V}$ of 
symmetries of ${\cal C}_{\cal V}$, understood as a subgroup}  
$
\Gamma_{\cal V} \subset  {\rm Auteq}({\cal C}_{\cal V}, {\cal S}).
$
\end{itemize}

One  should have the following additional features of this picture: 

\begin{itemize}
\item  A family of  complex open Calabi-Yau threefolds ${\cal Y}_b$  
over a base ${\cal B}$, which is smooth over ${\cal B}_0 \subset {\cal B}$, 
and a fully faithful functor
$
\varphi_{\cal V}: {\cal C}_{\cal V} \lra {\cal F}({\cal Y}_b)
$ 
to the 
Fukaya category ${\cal F}({\cal Y}_b)$ of  ${\cal Y}_b$, where $b \in {\cal B}_0$. 
It sends objects of cluster collections to the objects of the Fukaya category represented 
by special Lagrangian spheres. 

\item 
The group 
$\pi_1({\cal B}_0, b)$ acts by symmetries of the category ${\cal F}({\cal Y}_b)$. 
The equivalence ${\varphi}_{\cal V}$ intertwines 
the actions of the groups $\Gamma_{\cal V}$ and $\pi_1({\cal B}_0, b)$:
\be \la{equiva}
{\varphi}_{\cal V}: ({\cal C}_{\cal V}, \Gamma_{\cal V}) \stackrel{}{\lra} ({\cal F}({\cal Y}_b), \pi_1({\cal B}_0, b)).
\ee
\end{itemize}

A cluster collection ${\cal S}_\alpha$ in the category ${\cal C}_{\cal V}$ determines a canonical 
$t$-structure on ${\cal C}_{\cal V}$, whose heart is an abelian category denoted by $A({\cal S}_\alpha)$. 

An array of cluster collections ${\cal S} = \{{\cal S}_\alpha\}$ in the category ${\cal C}_{\cal V}$  
determines 
an open domain in the space 
of  stability conditions ${\rm Stab}({\cal C}_{\cal V})$ on the category 
${\cal C}_{\cal V}$:
$$
{\rm Stab}_{{\cal S}}({\cal C}_{\cal V}) \subset {\rm Stab}({\cal C}_{\cal V}). 
$$
Namely, ${\rm Stab}_{{\cal S}}({\cal C}_{\cal V})$ consists of stability conditions whose heart is 
one of the hearts $A({\cal S}_\alpha)$.  
If any two of the cluster collections ${\cal S}_\alpha$ are related 
 by a sequence of mutations, the domain ${\rm Stab}_{{\cal S}}({\cal C}_{\cal V})$ is connected. 

\begin{itemize}
\item  One should have  an embedding of    
${\cal B}$ into the stack of stability conditions ${\rm Stab}({\cal C}_{\cal V})$:
$$
i: {\cal B}\hra {\rm Stab}({\cal C}_{\cal V}).
$$

\item Then the domain ${\rm Stab}_{\cal S}({\cal C}_{\cal V})\subset 
{\rm Stab}({\cal C}_{\cal V})$ determines a subspace ${\cal B}_{\cal S}\subset 
{\cal B}_0$:
$$
\begin{array}{ccc}
{\cal B}_{\cal S}& \hra &{\rm Stab}_{\cal S}({\cal C}_{\cal V})\\
\downarrow &&\downarrow \\
{\cal B}_0 &\hra &{\rm Stab}({\cal C}_{\cal V})
\end{array}
$$
One should have an isomorphism
$$
\Gamma_{\cal V} = \pi_1({\cal B}_{\cal S}).
$$
Combined with the  map induced by the embedding ${\cal B}_{\cal S} \hra {\cal B}_0$, it gives rise 
to a map  
$ 
\Gamma_{\cal V} \lra \pi_1({\cal B}_0, b), 
$ 
intertwined by the functor (\ref{equiva}).   

\end{itemize}

The pair $({\cal C}_{\cal V}, \Gamma_{\cal V})$ is easily described 
combinatorially. 
Contrary to this, the Fukaya category ${\cal F}({\cal Y}_b)$, whose definition we do not 
discuss, is a complicated object. 
The pair $({\cal C}_{\cal V}, \Gamma_{\cal V})$ serves as a combinatorial model 
of a full subcategory of  the Fukaya category ${\cal F}({\cal Y}_b)$. 

The homological mirror symmetry conjecture of M. Kontsevich predicts 
equivalence of a Fukaya category of a Calabi-Yau threefold ${\cal Y}$ with an ${\rm A}_\infty$-category 
of complexes of 
coherent sheaves on the mirror ${\cal Y}^\vee$. 

The above picture adds a new ingredient: a combinatorial description of the ${\rm A}_\infty$-category.  
 It can be viewed as a manifestation of a general principle:
All interesting categories are equivalent to Fukaya categories.

The main message of our paper can be formulated as follows:

\begin{itemize}

\item 
{\it Ideal webs on decorated surfaces $S$ play the same role for moduli spaces 
of $G$-local systems on $S$ where $G$ is of type ${\rm A}_m$ as ideal triangulations for the groups of type ${\rm A}_1$. 
In particular ideal webs provide an explicit description of a ``tame part'' of the 
corresponding cluster atlas, \underline{canonical potentials} for the related quivers, 
and therefore an explicitly defined 3d CY category with an array of cluster collections,  
together with its symmetry group.}
\end{itemize}

${\rm A}_{m}$-webs on  $S$ 
provide a minimal $\Gamma_S$-equivariant 
collection of cluster coordinate systems 
on each of the  moduli spaces  (\ref{8.15.11.1}), connected in the following sense: 
any two of them are related by a sequence of elementary transformations  
involving only the coordinate systems from the  collection. 
Similarly, ${\rm A}_{m}$-webs are  a minimal 
 topological data on $S$ which allow to define 
combinatorially a 3d CY category with a large symmetry group. 

\vskip 3mm
The fundamental 
geometric structures on surfaces generalising ideal triangulations 
are the Gaiotto-Moore-Neitzke spectral networks, 
introduced in \cite{GMN4}, and studied further in \cite{GMN5}. 
We discuss spectral networks in comparison with ideal webs in Section \ref{sec1.8}. 

\paragraph{Our main example: $G$-local systems on decorated surfaces.} 
In this example, which is the one we elaborate in this 
paper, the ingrediants of the above picture look as follows:

\begin{itemize}

\item Varieties ${\cal V}$ are  
variants of the moduli spaces of $G$-local systems on a decorated surface $S$, 
where $G$ is $SL_m$ or $GL_m$,  equipped with 
cluster variety structures defined in \cite{FG1}. 

\item The quivers assigned to ideal webs on $S$ 
describe $\Gamma_S$- equivariant  cluster structures of these moduli spaces. These quivers 
are equipped with \underline{canonical potentials}. 

Therefore they give rise to 
a combinatorially defined 3d CY category with a $\Gamma_{G, S}$-equivariant 
 array of cluster collections. The  group  $\Gamma_{G, S} $ is an extension of 
the mapping class group $\Gamma_S$ by a generalised braid group, see Section \ref{SSec2.3}. 

\item The base ${\cal B}$ is the base 
 of the universal Hitchin integrable system for the pair $(G, S)$. 


\item The family of open CY threefolds over ${\cal B}$ is the one related to the Hitchin integrable system 
of the pair $(G, S)$ by \cite{DDDHP}, \cite{DDP}, \cite{G}, 
\cite{KS3}.  

\end{itemize}

\paragraph{The case $G=SL_2$.}  In this case, when $S$ is 
sufficiently general, the existence of an equivalence $\varphi$  
is proved by Smith \cite{S}, building on 
the works of 
Derksen, Weyman and Zelevinsky \cite{DWZ}, Gaiotto-Moore-Neitzke \cite{GMN1}-\cite{GMN3}, Kontsevich and 
Soibelman \cite{KS1}, 
\cite{KS3},
Labardini-Fragoso \cite{LF08},  \cite{LF12}, Keller and Yang \cite{KY}, 
Bridgeland and Smith \cite{BrS}, and others. 
In particular: 

Labardini-Fragoso \cite{LF08}, \cite{LF12} 
introduced and studied canonical potentials related to surfaces with ideal triangulations 
 as well as the tagged ideal triangulations  \cite{FST}. 

Bridgeland and Smith \cite{BrS} proved that the spaces of stability conditions 
for the combinatorial categories provided by Labardini-Fragoso's work
are naturally isomorphic to 
the spaces of quadratic differentials of certain kind on the corresponding Riemann surfaces. 

K. Nagao in his lecture \cite{N12}   
 outlined a 
combinatorial action  of a mapping class group.

\paragraph{Another example: dimer model integrable systems.} The dimer integrable system is assigned to 
a Newton polygon ${\Bbb N}$, and was studied in \cite{GK}. It gives rise to a cluster Poisson variety 
${\cal X}_{\Bbb N}$ and a map given by the dimer model Hamiltonians 
$$
{\cal X}_{\Bbb N} \lra {\cal B}_{\Bbb N}.  
$$
The story goes as follows. Instead of ideal bipartite graphs on a decorated surface 
$S$ we consider 
minimal bipartite graphs $\Gamma$ on a torus ${\rm T}^2$. Each of them 
is related to a convex polygon ${\Bbb N} \subset \Z^2$. The set of  
minimal bipartite graphs  on a torus  modulo isotopy with a given Newton polygon  
is finite, and any two of them are related by 
the elementary moves. We assign to such a graph $\Gamma$ the moduli space ${\cal L}_\Gamma$ 
of $GL_1$-local systems 
on $\Gamma$, and glue them via birational cluster Poisson transformations ${\cal L}_\Gamma \to 
{\cal L}_{\Gamma'}$ corresponding to two by two 
moves into a cluster variety ${\cal X}_{\Bbb N}$. 
The base  ${\cal B}_{\Bbb N}$ is essentially the space of polynomials $P(x,y)$ with the Newton polygon ${\Bbb N}$. 

Just as above, the bipartite graphs $\Gamma$ give rise to quivers with canonical potentials 
related by two by two moves, and 
hence to 3d CY category ${\cal C}_{\Bbb N}$, determined just by 
the polygon ${\Bbb N}$,   with an array of cluster collections. 
The corresponding CY open threefold is the toric threefold related 
to the cone over the Newton polygon ${\Bbb N}$. The geometry of these threefolds was extensively studied by 
the Physics and Math communities in the past. 

In particular any convex polygon ${\Bbb N}\subset \Z^2$, considered up to translations, gives rise to a group 
$\Gamma_{\Bbb N}$, realised as the categorical symmetry group.  This group is an extension 
of the cluster modular group of the cluster Poisson variety ${\cal X}_{\Bbb N}$ by a generalised braid group. 


\subsection{Concluding remarks} \la{sec1.8}

\paragraph{Bipartite graphs on decorated surfaces with abelian fundamental groups.} 
Our coordinates on the moduli spaces (\ref{8.15.11.1}) 
have many common features with Postnikov's cluster 
coordinates on the Grassmannians \cite{P}, see also \cite{S}
 -- the latter 
are given by a combinatorial data (directed networks on the disc) closely related to 
bipartite graphs on the disc. 
In particular, in both cases we have 
a subatlas of the cluster atlas closed under two by two moves, and the corresponding 
cluster coordinates are given simply by minors. 
Postnikov's story was partially generalized to the case of an annulus in \cite{GSV}. 
Minimal bipartite graphs on decorated surfaces, 
with special attention to the 
minimal bipartite graphs on the torus, were studied 
in   \cite{GK}. 

Minimal bipartite graphs on decorated surfaces, considered  modulo the two by two moves,  
admit a nice classification when the  surface has an abelian fundamental group, i.e. either a disc,  an annulus, or a torus.  Namely,  
the set of equivalance classes is finite, and has an explicit 
description. For the disc this follows from \cite{T} and \cite{P}, the annulus was considered in \cite{GSV}, 
and the torus in \cite{GK}. 

For an arbitrary surface 
the classification is a wild problem. 
Ideal bipartite  graphs provide an interesting subclass of  
bipartite surface graphs. 

\paragraph{Gaiotto-Moore-Neitzke spectral networks.} It would be very interesting to 
compare ideal webs with the Gaiotto-Moore-Neitzke spectral networks \cite{GMN4}. 

$SL_m$-spectral networks describe a geometric object  
on a Riemann surface $C$ related to a point of the Hitchin base, i.e. 
a collection $(t_2, ..., t_m)$ where $t_k \in \Omega^{\otimes k}(C)$ are holomorphic $k$-differentials on $C$. 
It is a remarkable geometric object, rather complicated,  whose propetries are not fully 
established yet.

Ideal webs are defined on topological surfaces. 
They parametrise cluster coordinate systems, and give rise to quivers with potentials, and therefore to  
3d CY categories with cluster collections. It is not clear whether 
one can associate similar features to an arbitrary spectral network. 

Ideal webs and spectral networks were developed independently, 
and relations between them are far from been clear. 
In \cite{GMN5} the ideal webs assigned to ideal triangulations are related to spectral networks. 
I believe  that all ideal webs are nice special examples of spectral networks. 

\paragraph{Ideal webs and scattering amplitudes.} Bipartite graphs on the disc,
 as well as slightly more general {\it on shell scattering diagrams} are in the center of the 
approach to the scattering amplitudes in the N=4 
super Yang-Mills theory developed in \cite{ABCGPT}. We hope that similar constructions 
related to ideal webs on arbitrary surfaces lead to a physically meaningful story. 
The parallels between the two stories are highlighted below: 
\begin{center}
\begin{tabular}{ c | c | c }
\bf{Scattering Amplitudes in the $N=4$ SUYM}&{\bf Moduli spaces of local systems}\\
\hline\hline
Disc with $n$ special points on the boundary &An arbitrary decorated surface \\
\hline
Grassmannian $Gr(m, n)~ \stackrel{\sim}{=}$ Configurations &Moduli space of \\
of $n$ vectors in an $m$-dimensional space&$GL_m$-local systems on surfaces \\
\hline
On shell scattering diagrams &  Ideal webs on surfaces\\
\hline
Scattering Amplitudes in the $N=4$ SUYM &  Similar integrals given by amalgamation\\
\end{tabular}\end{center}

\paragraph{The structure of the paper. }
In Section \ref{sec1.4} we present, for the convenience of the reader, the background material 
on quivers with potentials and 3d Calabi-Yau categories with cluster collections. 
In particular we recall the Kontsevich-Soibelman correspondence 
between these two objects (\cite{KS1}, Section 8). 
We only very briefly recall an alternative approach 
of Keller-Yang \cite{KY} via Ginzburg algebras \cite{G}. 
In Section \ref{SSec2} we define symmetry groups of 3d CY categories with cluster collections. 
The story is very similar to the definition of the cluster modular group \cite{FG2}. 
In Section \ref{SSec2.3} we define an extended mapping class group. 
In Section \ref{SSec2.4} we formulate  the conjecture 
relating the combinatorial 3d CY categories coming from ideal webs 
and Fukaya categories of open CY threefolds related to reductive groups of type A. 

Sections \ref{SSec3} - \ref{SSec5} are entirely geometric: we study 
ideal webs and introduce the related cluster coordinates. 
They are completely independent from the other parts of the paper. 


\paragraph{Acknowledgements.} 
The geometric part of the paper, Sections 2-4,  describing ideal webs and cluster coordinates was written at the IHES during 
the Summer 2011. I am grateful to 
Maxim Kontsevich for useful discussions, and to 
the IHES for the  hospitality and support. 
I was supported by the NSF grants DMS-1059129 and DMS-1301776. 

The whole story was discussed at talks at the Mirror Symmetry conference 
at Miami (January 2012) and Cluster Algebras conference at Northeastern 
Univeristy (April 2012). I am grateful to participants for useful discussions.

\section{Ideal webs on decorated surfaces} \la{SSec3}

\subsection{Webs and bipartite graphs on decorated surfaces \cite{GK}}\la{ssec2.1}

In Section \ref{ssec2.1} we assume  that $S$ is an oriented surface 
with or without boundary. 
We do not assume that it is a decorated surface.

A {\it strand} on $S$ is either an oriented 
loop, or 
an oriented path connecting two boundary points.

\bd \la{Def1}
A {\bf web ${\cal W}$} on an oriented surface $S$  
is a union of 
strands $\{\gamma\}$ in generic position on $S$, considered modulo isotopy,  such that: 

\begin{enumerate}
\item As we move along a strand, 
the orientations of the crossing strands alternate.  

\item The subset ${\cal W}$ is connected. The components of $S-{{\cal W}}$, called {\rm domains}, are discs.

\end{enumerate}
\ed

Our next goal is to show that webs on $S$ can be encoded by bipartite graphs on $S$. 

Any vertex of a graph is either {\it internal}, that is of valency $\geq 2$, or {\it external}. 
\bd

A {\em bipartite graph on a surface $S$}, possibly with a boundary, is a graph with external vertices at the boundary, 
whose internal vertices and internal edges form an abstract bipartite graph. The external vertices are not colored. 
\ed

Lemma \ref{1.21.12.1a} is borrowed from  \cite{GK}. We provide a proof for the convenience of the reader. 

\bl \la{1.21.12.1a} 
Let $S$ be an oriented surface, possibly with boundary. 
There is a bijection between the isotopy classes of 
webs and connected bipartite graphs, whose faces 
are discs, on $S$. 
\el

\begin{proof} A path on the boundary of a domain 
on an oriented surface
is {\it positively oriented} if its orientation coincides with the 
 orientation of the boundary,  provided by the surface orientation. Otherwise it is {\it negatively oriented}. 
We picture positive orientation as clockwise. 

A web ${\cal W}$ on $S$ gives rise to a bipartite graph $\Gamma$ on 
$S$ as follows. 

Condition $1)$ implies 
that there are three types of the domains: 
\begin{enumerate}

\item {\it $\bullet$-domains}: the sides are oriented positively.

\item {\it $\circ$-domains}: the sides are oriented negatively.

\item {\it Faces}: the directions of the sides alternate. 
\end{enumerate}

Since each domain is a disc (Condition 2), 
assigning a $\circ$ (respectively $\bullet$) central 
point to each $\circ$ (respectively $\bullet$) domain, we get a 
bipartite graph $\Gamma$ on $S$. 
Its edges correspond to the crossing points of the strands, 
see Figure \ref{bgc2}. Indeed, since the strands  are 
in generic position, every crossing point of ${\cal W}$  
is an intersection of two arrows. 
So it determines the $\circ$  and the $\bullet$-domains. 

\begin{figure}[ht]
\centerline{\epsfbox{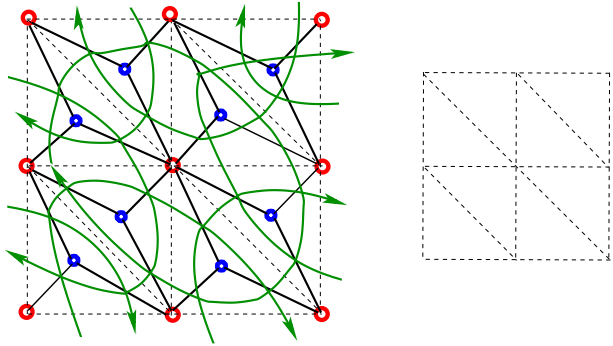}}
\caption{Zig-zag strands (green)   for the 
bipartite graph assigned to a triangulation (shown by dotted lines) 
of a torus with four punctures. }
\label{bgc2}
\end{figure}

\vskip 3mm
Let us now go in the opposite direction, and assign to a bipartite graph $\Gamma$ on $S$ a web 
${\cal W}$. 

A {\it zig-zag path} on a bipartite ribbon graph  
$\Gamma$ is an oriented path on   
$\Gamma$ which turns maximally left at $\bullet$-vertices,  
and maximally right at $\circ$-vertices. We isotope 
zig-zag paths  
slightly to strands going clockwise 
(respectively counterclockwise)
around the $\bullet$ (respectively the $\circ$) vertices. 
Namely, at each vertex $v$ on a zig-zag path $\gamma$ we push the path $\gamma$  a bit 
inside of the angle formed by the edges incident to $v$.
These strands are called {\it zig-zag strands}. 

A bipartite graph $\Gamma$ on $S$ gives rise to a web 
${\cal W}$, given by the zig-zag strands of $\Gamma$. They satisfy Conditions 1 and 2. 
\end{proof}

\begin{figure}[ht]
\centerline{\epsfbox{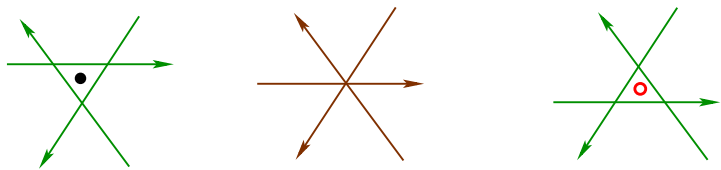}}
\caption{A triple point diagram (middle) and the two associated webs (left and right).}
\label{gc4}
\end{figure}

\paragraph{Webs and triple point 
diagrams.} 
A web on a surface is closely related to a slight generalisation of Dylan
 Thurston's {\it triple point 
diagram on the disc}:

\bd [\cite{T}]
A {\bf triple point 
diagram on the disc} is a collection of oriented strands in the disc 
whose intersection points are as in the middle of Figure \ref{gc4}, 
and the endpoints of the strands are distinct points on the boundary of the disc. 
\ed

We consider triple point 
diagrams on an oriented 
 surface with boundary rather then on a disc, and 
resolve all triple point crossings.

An  important difference is that there are \underline{two} non-equivalent ways 
to resolve triple point crossings, shown on the left and on the right on  Figure \ref{gc4}. 
So there is a pair of webs associated with a triple point diagram, and vice verse. 
See Lemma-Construction \ref{inv*} for more details. This is the 
reason we work with webs rather then with triple point diagrams. 

\paragraph{The conjugate surface \cite{GK}.} A {\it ribbon graph} 
is a graph with  a cyclic order 
of the edges at every vertex. 
A {\it face path} on a ribbon graph $\Gamma$ is an oriented  path on the graph 
which turns right at every vertex. 
A closed face path is called a {\it face loop}. 

A graph embedded into an oriented
 surface  has an induced ribbon structure. 
We  assume that its external vertices are on the boundary of $S$. 

A bipartite ribbon graph $\Gamma$  
gives rise to a new ribbon graph $\Gamma^*$, 
the {\it conjugate graph}, obtained by reversing the cyclic orders at 
all $\bullet$-vertices\footnote{Reversing the cyclic order at $\circ$-vertices we get another ribbon graph 
which differs from 
$\Gamma^*$ by the orientation.}. 
Notice that $\Gamma$ and $\Gamma^*$ are the same graphs  -- only their ribbon structures are different. 
There is a bijection:
\be\la{ZZBHa}
\{\mbox{zig-zag strands $\{\gamma_i\}$ on $\Gamma$}\} \longleftrightarrow \{\mbox{face paths 
 on
 $\Gamma^*$}\}.
\ee
The orientations of zig-zag strands match  
the orientations of the face paths on $\Gamma^*$.

A ribbon graph $\Gamma$ gives rise to a surface in a standard way --
we replace the edges of $\Gamma$  by thin ribbons. 
Let $\Sigma$ be the topological surface with boundary 
corresponding to the ribbon graph $\Gamma^*$. 
It can be visualised by taking the surface glued from the ribbons assigned to the edges of $\Gamma$,  
cutting each of the ribbons in the middle, twisting its ends incident to black 
vertices  by $180^\circ$, and gluing it back. 
Then bijection (\ref{ZZBHa}) gives rise to a bijection 
\be\la{ZZBH}
\{\mbox{zig-zag strands $\{\gamma_i\}$ on $\Gamma$}\} \longleftrightarrow \{\mbox{boundary components
 of 
 $\Sigma$}\}.
\ee
Orientations of zig-zag strands match  
orientations of boundary components on $\Sigma$.

 \subsection{Spectral webs and spectral covers}\la{ssec2.2}

Let ${\cal W}$ be a web on a decorated surface $S$. Let $p: \widetilde S \to S$  be the universal cover of  $S$. 
The surface $\widetilde S$ inherits a set of special points: the preimage of the ones on $S$. 
So it is a decorated surface. 
Then $\widetilde {\cal W}:= p^{-1}({\cal W})$ is a web on $\widetilde S$. 
Defining different types of webs ${\cal W}$ on $S$, like 
minimal webs, spectral webs,  ideal webs, 
we work with the web $\widetilde {\cal W}$. We use the same strategy 
introducing cluster coordinates  related to an ideal web on $S$. Since the web $\widetilde {\cal W}$ is 
$\pi_1(S)$-equivariant, all constructions lead to notions and objects 
related to  the surface $S$ itself.   

\begin{figure}[ht]
\centerline{\epsfbox{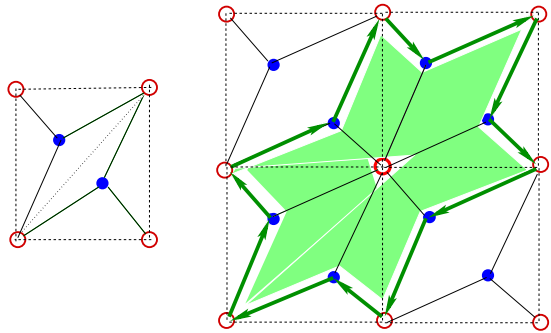}}
\caption{On the left: the bipartite graph for a punctured torus, triangulated by a diagonal.  
On the right: 
the  bipartite graph on its universal cover
with a (green) zig-zag loop $\gamma$. 
The domain  $S_\gamma$ is the interier of the  filled 12-gon. It contains a single $\circ$-vertex. }
\label{bc11}
\end{figure}

\paragraph{Key     example.} Figure \ref{bc11} illustrates how misleading could be 
the geometry of a web considered on $S$ rather then its universal cover. 
   Consider a torus  with a single special point. 
Pick an ideal triangulation. 
The corresponding bipartite graph has a unique $\circ$-vertex, and two $\bullet$-vertices, 
as illustrated on the left of Figure \ref{bc11}. The web assigned to this bipartite graph 
consists of a single selfintersecting zig-zag path. 
However the universal cover, illustrated on  the right of Figure \ref{bc11}, consists of nonselfintersecting 
zig-zag loops. Each of them surrounds a single special point, which in this particular case is 
also a single $\circ$-vertex 
inside of the loop. 
We show in Figure \ref{bc11} one of its zig-zag loops, surrounding a $\circ$-vertex. 
The corresponding domain is a 12-gon, 
filled green. Therefore this 
 web on the torus is an ideal ${\rm A}_1$-web  for the definition given below.

\begin{figure}[ht]
\centerline{\epsfbox{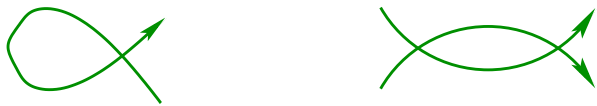}}
\caption{Minimality: no selfcrossings and no parallel bigons on the universal cover. }
\label{bgc4}
\end{figure}

\bd \la{Def2}
A web ${\cal W}$ on a 
surface $S$ 
is {\bf minimal} if the web $\widetilde {\cal W}$ 
  has the following property:

\begin{itemize}
\item Strands have no selfcrossings, and there are no parallel bigons, see Figure \ref{bgc4}.
\end{itemize}
A bipartite graph on a 
surface $S$ is {minimal} if the corresponding web is minimal.
\ed

\bd A web ${\cal W}$ on a 
surface $S$ 
is {\bf spectral} if 
the web $\widetilde {\cal W}$  has the following property: 
\begin{itemize}

\item Each strand $\gamma$ of $\widetilde {\cal W}$ cuts the surface $\widetilde S$ into two 
components:
\be \la{6.20.11.1}
\widetilde S - \gamma = \widetilde S_{\gamma} \cup \widetilde S'_{\gamma}.
\ee

\end{itemize}
\ed

\noindent
The connected components in (\ref{6.20.11.1}) are orientated by the  orientation of $S$.  
The domain $\widetilde S_{\gamma}$ is the one for which its orientation 
induces the original (clockwise on the pictures) orientation of 
 $\gamma$. 

\paragraph{Remark.} Condition (\ref{6.20.11.1}) is equivalent to the condition that 
$\gamma$ does not have selfintersections. Minimality implies condition (\ref{6.20.11.1}), but not 
the other way: we can have parallel bigons on a spectral web.

\vskip 3mm
Given  a spectral web ${\cal W}$ on $S$, let us define 
a {\it spectral surface}  ${\Sigma}_{\cal W}$, together with 
a ramified {\it spectral cover} map
$
\pi: {\Sigma}_{\cal W} \lra S. 
$

Consider the bipartite graph $\widetilde \Gamma \subset \widetilde S$ assigned to the web $\widetilde {\cal W}$. 
Each zig-zag path $\gamma$ on $\widetilde \Gamma$ bounds 
the unique component $\widetilde S_{\gamma}$ on $\widetilde S$. 
Now consider $\widetilde \Gamma$ as an abstract ribbon bipartite graph. 
We attach the  component $\widetilde S_{\gamma}$ to the graph $\widetilde \Gamma$ along the strand 
$\gamma$. 
Performing this procedure for all zig-zag paths, we get 
a  surface ${\Sigma}_{\widetilde {\cal W}}$. 
By construction, each component $\widetilde S_{\gamma}$ becomes a part of 
${\Sigma}_{\widetilde {\cal W}}$, denoted by $F_{\gamma}$.  So 
$$
{\Sigma}_{\widetilde {\cal W}} = \widetilde \Gamma \cup\cup_{\gamma_i}F_{\gamma_i}. 
$$
Using the tautological isomorphisms 
$F_{\gamma}\stackrel{\sim}{\to} \widetilde S_{\gamma}$, we get a map, called {\it spectral cover map}
\be \la{9.22.11.1gc}
\widetilde 
\pi: {\Sigma}_{\widetilde {\cal W}} = \widetilde \Gamma 
\cup \cup_{\gamma_i}F_{\gamma_i} \lra \widetilde S.
\ee

The fundamental group $\pi_1(S)$  acts on 
both $\widetilde S$ and ${\Sigma}_{\widetilde {\cal W}}$, and by the construction its action commutes
with  the 
projection $\pi$. So the group $\pi_1(S)$ acts by automorphisms of the cover. Set 
$$
{\Sigma}_{{\cal W}}:= {\Sigma}_{\widetilde {\cal W}}/\pi_1(S). 
$$
 Then we get a commutative diagram, where the horisontal maps are given by the factorisation 
by the free action of the group $\pi_1(S)$:

$$
\begin{array}{ccc}
{\Sigma}_{\widetilde {\cal W}}&\lra &{\Sigma}_{{\cal W}}:= {\Sigma}_{\widetilde {\cal W}}/\pi_1(S)\\
\widetilde \pi\downarrow &&\downarrow \\
\widetilde S&\lra &S= \widetilde S/\pi_1(S)
\end{array}
$$
Therefore we arrive at a   map, which we call the spectral cover map: 
\be \la{covermap}
\pi: {\Sigma}_{{\cal W}}  \lra S.
\ee
The following theorem summarises its properties.

\bt \la{1.21.12.1} Let ${\cal W}$ be a spectral web on an oriented
 surface $S$. Then the spectral cover 
map (\ref{covermap}) can be ramified only at the $\bullet$-vertices of the associated bipartite graph $\Gamma$. 
The ramification index of a $\bullet$-vertex is its valency minus two. 
\et

\begin{proof} This is a local question, so it is sufficient to argue for the graph $\widetilde \Gamma$ 
on the universal cover. 
The map $\widetilde \pi$ can be ramified only at the vertices of the graph $\widetilde \Gamma$. 
Take a vertex $v$. Let $n$ be the valency of $v$. Then a little neighborhood $U$ of $v$ is a union of $n$ sectors 
$C_1, ..., C_n$ 
formed by pairs of consecutive edges emanating from $v$. A zig-zag path $\gamma$ 
passing through $v$ goes along the boundary of one of the sectors $C_i$. 
If $v$ is a $\bullet$-vertex, the disc $S_{\gamma}\cap U$ contains all these sectors but the one $C_i$. 
Thus the union of $S_{\gamma}\cap U$ over all zig-zag paths $\widetilde 
\gamma$ passing through an $n$-valent 
 $\bullet$-vertex $v$  
covers  $n-1$ times $U-v$.  
It is $1:1$ at the vertex. So the ramification index of the map $\pi$ at a $\bullet$-vertex $v$ is $n-1$. 
If $v$ is a $\circ$-vertex, then $S_{\gamma} \cap U$ 
contains only one of the sectors, $C_i$. So the map $\pi$ is unramified near a $\circ$-vertex. 
\end{proof} 

\bd
The {\bf rank} ${\rm rk}({\cal W})$ of a spectral web ${\cal W}$ is the degree of the map (\ref{covermap}).
\ed

\bl \la{ranklemma}  For an ideal web ${\cal W}$  on a 
 decorated surface $S$ one has 
$$
{\rm deg}(\widetilde \pi) =  {\rm deg}(\pi).
$$
\el

\begin{proof}
The degree ${\rm deg}(\widetilde \pi)$ is the number of points in the fiber $\widetilde \pi^{-1}(x)$ 
for a generic $x \in \widetilde S$. The group $\pi_1(S)$ acts discretely, so 
the factorisation by the action of $\pi_1(S)$ does not change the cardinality of the fiber. 
\end{proof}

\paragraph{Spectral surface $=$ completed conjugate surface.} \la{ssec2.3.5} 
Let us add to the conjugate ribbon graph $\Gamma^*$ 
its faces $F_{\varphi_i}$ bounding the face paths $\varphi_i$. 
Namely,  if $\varphi_i$ is a face loop, then $F_{\varphi_i}$ is a disc bounded by this loop. 
If $\varphi_i$ is an unclosed face path, then $F_{\varphi_i}$ is a disc; the half of its boundary is 
identified with 
the path $\varphi_i$, and the rest becomes boundary of $\Sigma$. We arrive at a  
{\it completed  conjugate surface ${\bf \Sigma}$}:
\be \la{2.13.12.1}
{\bf \Sigma} = \Gamma^* \cup\cup_{\varphi_i}F_{\varphi_i}. 
\ee
The surface ${\bf \Sigma}$ still may have a boundary. 
Topologically it is obtained by filling the holes on $\Sigma$. 
\bl
Let ${\cal W}$ be a spectral web on   $S$. 
Then the spectral surface ${\Sigma}_{{\cal W}, S}$  
is canonically identified with the 
{completed conjugate surface} ${\bf \Sigma}$ of the bipartite 
ribbon graph $\Gamma$ assigned to ${\cal W}$:
\be \la{1.18.12.1}
{\bf \Sigma}  = {\Sigma}_{{\cal W}}.
\ee 
\el

\begin{proof} We construct the spectral surface ${\Sigma}_{{\cal W}}$ 
by attaching the  domain $S_\gamma$ to every zig-zag path $\gamma$ on the bipartite graph $\Gamma$. 
In our case domains $S_\gamma$ are discs. Furthermore, zig-zag paths $\gamma$ on $\Gamma$ are in bijection with the face paths on $\Gamma^*$ -- see (\ref{ZZBHa}). \end{proof}

From now on we make no distinction between the spectral surface ${\Sigma}_{{\cal W}}$ 
of a spectral web and the completed conjugate surface ${\bf \Sigma}$,
and keep the latter notation  only.

The preimages on ${\bf \Sigma}$ of the special points on $S$ are declared to be the special points on ${\bf \Sigma}$. 
\vskip 3mm
Let $\Gamma$ be the bipartite graph assigned to the web ${\cal W}$. 
Let $\varphi$ be a face path on the conjugate graph $\Gamma^*$. Let ${\bf \Sigma}_\varphi$ be the unique domain obtained by cutting 
${\bf \Sigma}$ along $\varphi$, so that the orientation of the boundary of ${\bf \Sigma}_\varphi$ coincides with the 
orientation of the face path $\varphi$. 
Each face path on $\Gamma^*$ surrounds just one special point on ${\bf \Sigma}$, providing a bijection
\be \la{1.21.12.10}
\{\mbox{face paths on $\Gamma^*$}\} \leftrightarrow \{\mbox{special points of ${\bf \Sigma}$}\}. 
\ee
$$
\{\mbox{face path $\varphi$}\} \longmapsto \{\mbox{the unique special point $s_\varphi$ inside of the domain ${\bf \Sigma}_\varphi$}\}. 
$$

\subsection{Ideal webs on decorated surfaces} \la{ssec2.3}

Let us introduce now the notion of an {\it ideal web}, which is the central 
new object in the paper. Ideal webs form a particular class of spectral webs.

From now on, $S$ is a decorated surface 
with $n>0$ special points $\{s_1, ..., s_n\}$. Set 
\be \la{18.1.12.2}
S^\times = S - \{s_1, ..., s_n\}.
\ee 
A bipartite graph on $S$ can have vertices at special points. 
A web on $S$ means a web on $S^\times$.

\paragraph{The double of a decorated surface with boundary.} 
Given an oriented decorated surface $S$ with boundary, its {\it double} $S_{\cal D}$ is  
an oriented decorated surface obtained by gluing $S$ with its mirror $S^\circ$, 
given by the surface $S$ with the opposite orientation,  
into a surface without boundary, by identifying the corresponding components. 
The set of special points on $S_{\cal D}$ 
is inherited from $S$ and $S^\circ$ as follows.  Each special point on the boundary of $S$ gives rise 
to a single special point on the double $S_{\cal D}$. Each special point $s$ inside of $S$ 
gives rise to two special points on 
$S_{\cal D}$, given by $s$ and its mirror. This way we get all special points on the double. 

A web ${\cal W}$ on $S$   gives rise to a mirror web ${\cal W}^\circ$ on $S^\circ$, 
given by the strands of the web ${\cal W}$ 
with the reversed orientations. Gluing the webs ${\cal W}$ and ${\cal W}^\circ$  we get 
a web ${\cal W}_{\cal D}$ on $S_{\cal D}$. 

We can develop the whole 
story for decorated surfaces without boundary first, and then 
just require that the story on decorated  surfaces with boundary is determined by the condition that 
everything going fine on its double. We, nevertheless, will 
spell some of the definitions, e.g. the minimality condition,  
for decorated surfaces with boundary.

\paragraph{Ideal webs.}

\bd\la{Def0}
A strand $\gamma$ on $S^\times$ is an {\bf ideal strand} if $S - \gamma$ has two  components, and

\begin{itemize}

\item
The component   $S_{\gamma}$, see (\ref{6.20.11.1}), is a disc which
 contains a \underline{unique} special point $s$.  
\end{itemize}
\ed

We say that an ideal strand 
$\gamma$ as above is {\it associated to $s$}. 

\bd \la{2.4.12.1a}
A web ${\cal W}$ on a decorated surface $S$ without boundary  is {\bf ideal} if

\begin{itemize}
\item
The web $\widetilde {\cal W}$ is a minimal web of ideal strands on 
$ \widetilde S$.
\end{itemize}
 \ed

Observe that the strands of $\widetilde {\cal W}$ can not have selfintersections due 
to the condition from Definition \ref{Def0} that $\widetilde S_{\gamma}$ is a disc. The minimality 
implies that parallel bigons are also prohibited. 

\bd A web ${\cal W}$ on a 
decorated surface $S$ with boundary is {\bf ideal}  if the web 
${\cal W}_{\cal D}$  on $S_{\cal D}$ is ideal.
\ed

We say that a bipartite graph $\Gamma$ on $S$ is ideal if the corresponding web is ideal.

A web on $S$ is a {\it strict ideal web} if it is already a 
minimal web of ideal strands.

\begin{figure}[ht]
\centerline{\epsfbox{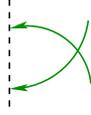}}
\caption{No parallel half-bigons: the boundary of $S$ is punctured. }
\label{gc4a}
\end{figure}

A web on a decorated 
surface $S$ with boundary is ideal if it consists of ideal strands, minimal, and satisfies the following condition: 

\begin{itemize}
\item There are no parallel half-bigons, see Figure \ref{gc4a}. 
\end{itemize}

\bl \la{2.15.12.1}
A strict ideal web on a decorated surface $S$ enjoys the following condition: 
\begin{itemize}
\item The strands associated with a given special point do not intersect. 
\end{itemize}
\el

\begin{proof} Due to, say,  the minimality condition, 
the strands have no selfintersections. So they are simple loops going clockwise around 
the corresponding special points. 
Two such loops going around the same special point do not intersect. Indeed, otherwise, 
since they  have the same orientation, this will 
create  a parallel bigon or half-bigon.  \end{proof}

\bc \la{2.4.12.11}
Let ${\cal W}$ be a  rank $m$ ideal web on $S$. Then the  
strands of the web $\widetilde {\cal W}$ associated with a special point $s \in \widetilde S$ 
form $m$ concentric (half)circles on $\widetilde S$ going clockwise 
around $s$. 
\ec

\paragraph{Key  example continued.} 
   The web in Figure \ref{bc11} assigned to an ideal triangulation of a torus with a single special point 
is an ideal web, but it is not a 
strict ideal web.

\vskip 3mm
Given an ideal web ${\cal W}$ on a decorated surface $S$, 
the {\it codistance} $\langle F, s\rangle$ between a face $F$ of the web $\widetilde {\cal W}$ 
and a special point $s$ on $\widetilde S$ is 
the number of strands associated with $s$ containing $F$. 

\bc \la{1.21.12.1d}  Let ${\cal W}$ be an ideal web  on $S$. 
Then for any face $F$ of the web $\widetilde {\cal W}$ one has 
$$
\sum_s\langle F, s\rangle =  {\rm rk}({\cal W}).
$$
In particular, the number of strands associated with a given special point 
 equals ${\rm rk}({\cal W})$. 
\ec

\begin{proof} Follows from Theorem \ref{1.21.12.1}. \end{proof}

\begin{figure}[ht]
\centerline{\epsfbox{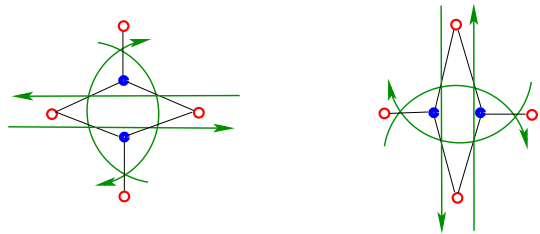}}
\caption{A two by two move of bipartite graphs.}
\label{gc2}
\end{figure}

\begin{figure}[ht]
\centerline{\epsfbox{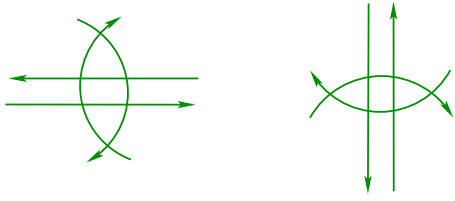}}
\caption{A two by two move of webs.}
\label{gra12}
\end{figure}

\paragraph{Two by two moves.} A two by two move of a  web on $S$ is a local 
transformation involving four zig-zag strands, shown in Figure \ref{gra12}, see \cite{T}, \cite{GK}. 
It clearly takes an ideal web to another ideal web. It preserves the rank. 
Indeed, it decreases by one the codistances 
to two special points and 
increases by one the codistances 
to other two ones. 
A two by two move amounts to a transformation of bipartite graphs, 
also known as a two by two move, 
see Figure \ref{gc2}.

\subsection{Ideal webs arising from ideal triangulations of decorated surfaces}\la{sec2.4.1}

Let $S$ be a decorated surface. 
Given an integer $m\geq 1$, we have assigned to an ideal triangulation of $T$ of $S$ two 
bipartite graphs on $S$,  the graphs $\Gamma_{{\rm A}_{m}}(T)$ and $\Gamma_{{\rm A}^*_{m}}(T)$. 
So there are  webs associated to each of them.  
Below we define  ${\rm A}_{m}$-webs and  ${\rm A}_{m}^*$-webs. Then we  show that 
the associated to these graphs webs are the ${\rm A}_{m}$-webs and ${\rm A}_{m}^*$-webs. 

\paragraph{${\rm A}_{m}$-webs and ${\rm A}_{m}^*$-webs.} 

\bd
An ideal web ${\cal W}$  on $S$ is an {\bf ${\rm A}_{m}$-web} 
(respectively {\bf ${\rm A}_{m}^*$-web}) 
if the web $\widetilde {\cal W}$ on the universal cover $\widetilde S$ has $m$  
(respectively $m+1$) strands around each special point,  and its rank  is $m+1$ .
\ed
As we show below, 
${\rm A}_{m}$-webs are related to moduli spaces of ${\rm A}_m$-local systems on $S$, 
while  ${\rm A}^*_{m}$-webs are related to moduli spaces of $GL_{m+1}$-local systems on $S$. 
The number of strands around each  special point of an ${\rm A}_{m}$ / ${\rm A}_{m}^*$-web 
is the rank of the corresponding reductive group. 

Recall that there are two  ways to resolve 
a triple point diagram, shown on Figure \ref{gc4}. 
So there is an involution acting on the webs. 
Lemma \ref{inv*} tells that it interchanges  ${\rm A}_{m}$-webs with 
${\rm A}_{m}^*$-webs. 

\bl \la{inv*}
There is an involution ${\cal W} \to *{\cal W}$ on the set of isotopy classes of webs, which preserves ideal webs on a decorated surface, and 
interchanges the ${\rm A}_{m+1}$-webs with the  ${\rm A}_{m}^*$-webs. The orbits of this involution are in bijection with 
the triple-point diagrams on $S$. 
\el

\begin{figure}[ht]
\centerline{\epsfbox{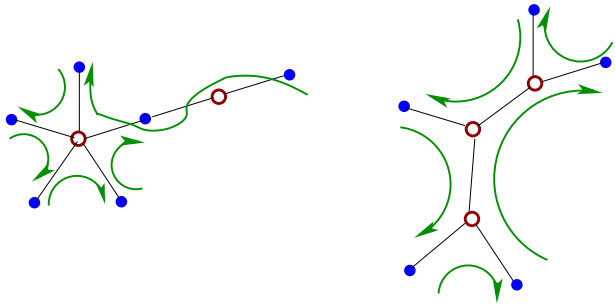}}
\caption{Breaking a $5$-valent $\circ$-vertex into three $3$-valent ones.}
\label{gc6}
\end{figure} 

\begin{proof} 
We replace an arbitrary bipartite graph 
by a bipartite graph  with $3$-valent $\circ$-vertices, so that the corresponding web of 
zig-zag strands does not change, see Figure \ref{gc6}.  Namely, collapse 
the segments containing $2$-valent 
$\circ$-vertices, 
collapsing the two $\bullet$-vertices of such segments into $\bullet$-vertices of higher valency, and 
removing the $2$-valent $\circ$-vertices. Then replace each $k$-valent $\circ$-vertex, $k>3$, 
by a pair of $\circ$-vertices 
of valencies $2$ and $k-2$, by inserting a $2$-valent $\bullet$-vertex between them, 
and keep doing so till no more $k>3$ valent $\circ$-vertices remains.  
Although the obtained bipartite graphs are different, their webs are isotopic.

Assuming all $\circ$-vertices are $3$-valent, we alter the  web
 at each $\circ$-vertex by moving one of the three strands 
surrounding this $\circ$-vertex across the intersection point of the other two. This way we get an involution ${\cal W} \to *{\cal W}$  
on the set of webs, see Figure \ref{gc4}. 
The involution  does not preserve the web rank. It interchanges the ${\rm A}_m$-webs and ${\rm A}_{m}^*$-webs. 

Assuming that all $\circ$-vertices are $3$-valent, and 
shrinking all $\circ$-domains of the corresponding web into points, we get a triple point diagram. 
Conversely, any triple point diagram can be resolved into a web in two different ways, 
so that each triple crossing is replaced by either $\circ$- or $\bullet$-vertex. 
\end{proof}

\begin{figure}[ht]
\centerline{\epsfbox{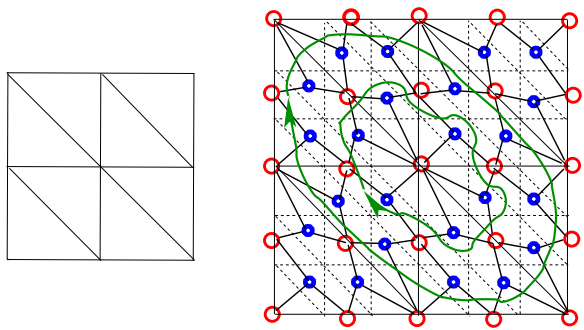}}
\caption{${\rm A}_2$-picture. On the 
left:  an ideal  triangulation of a torus with four special points.  
On the right: its $3$-triangulation $T_3$ -- shown by dotted lines; 
the bipartite graph  assigned to $T_3$ -- with  $\bullet$ and $\circ$- vertices;  
and the two zig-zag loops around the central special point.}
\label{bgc1}
\end{figure}

\paragraph{${\rm A}_{m}$-webs and ${\rm A}_{m-1}^*$-webs assigned to ideal triangulations.}

Recall that an ideal triangulation $T$ of $S$ is regular if every triangle of the triangulation 
has three distinct sides. 

\bl\la{1.20.13.1}
Let $T$ be an ideal regular triangulation  of $S$. Then:

i) The web associated to  the bipartite graph $\Gamma_{{\rm A}_{m}}(T)$ is an ${\rm A}_{m}$-web on $S$. 

ii)  
The web associated to the bipartite graph $\Gamma_{{\rm A}_{m}^*}(T)$ is a ${\rm A}_{m}^*$-web on $S$. 
\el

\begin{proof}

i) It is clear now from Figure \ref{gra11i} that zig-zag strands on the graph are ideal strands. Indeed, 
given a special point $s$ of $S$, 
take the triangles $t_1, ..., t_k$ sharing $s$, ordered clockwise. 
In each triangle $t_i$ consider the zig-zag strands parallel to the side opposite to $s$. 
They form $m$ arcs going clockwise around $s$. Considering all 
triangles sharing the puncture, we get $m$ strands associated with $s$. 
The rank of the obtained ideal web is $m+1$.  

ii) The proof if very similar to i), and thus is omitted. \end{proof}

\bl \la{sl2} There is a bijection between 
${\rm A}_1$-webs and ideal triangulations of $S$, such that 
two by two moves of ideal webs correspond to flips of triangulations,  
see Figure \ref{gra13a}. 
\el

\begin{proof} We can assume that all $\bullet$-vertices are $3$-valent. 
Given an ${\rm A}_1$-web on a decorated surface $S$, 
its zig-zag strands are in the bijection with the set of special vertices of $S$. 
Since the rank of the web is two, each face of the web determines two special points, and 
defines a homotopy class of a path connecting them. 
Let us choose their representatives, and declare them the edges. Then each $\bullet$-vertex 
provides us a triangle, and we arrive at an ideal triangulation of $S$. 
By the construction, the ${\rm A}_1$-web assigned to it is the original web. 
Then $\circ$-vertices are the special 
points. 
The second claim is clear from Figure \ref{gra13a}. 
\end{proof}

\subsection{Quivers associated to bipartite graphs on decorated surfaces} \la{sec2.5} 

We follow   \cite[Section 4.1.3]{GK}, with 
modifications for surfaces with boundaries. 

Let $S$ be a  decorated surface. 
A face of a bipartite graph $\Gamma$ on $S$ is  {\it external} if it 
intersects the boundary   
of  $S$, or its interior  contains a special point. It is called {\it internal} otherwise. 

 Below we use a slightly relaxed definition of a quiver, where we do not require 
the vectors $\{e_i\}$ to form a basis of $\Lambda$, replacing this by the  
condition that they generate the lattice. 
We call the vectors $\{e_i\}$ the {\it basic vectors}. Each basic vector is 
classified as either frozen or non-frozen. 

\bd\la{Quasiquivers}
A quasi-quiver  is a data 
$\Bigl(\Lambda, \{e_i\}, (\ast, \ast)\Bigr)$,  where 
$\Lambda$ is a lattice, $\{e_i\}$ is a 
finite set of vectors of two types -- \underline{non-frozen} and \underline{frozen} -- 
generating the lattice $\Lambda$,  and 
$(\ast, \ast)$ a skew-symmetric $\frac{1}{2}\Z$-valued 
bilinear form on $\Lambda$ such that $(e_i, e_j) \in \Z$ unless both $e_i$ and $e_j$
 are frozen.
\ed

All quivers in this paper are either traditional quivers, 
or quasi-quivers in the sense of Definition \ref{Quasiquivers}  
where the basic vectors $e_i$ satisfy just a single relation: 
$\sum_ie_i=0$. Abusing terminology, we  refer to them as quivers. 
In all examples in this paper $(e_i, e_j) \in \Z$ for any $e_i, e_j$. 

\paragraph{A quiver for a decorated surface $S$  without boundary.} 
So $S$ is a compact surface with a non-empty subset of special points. Let $\Gamma$ be a bipartite  graph on $S$. 

Let us assume that special points can not be   vertices of $\Gamma$. Then
$$
\Gamma \subset S^\times:= S - \{\mbox{\rm special points}\}. 
$$ 
Set
\be \la{10.19.14.1}
\Lambda_\Gamma:= {\rm Ker}\Bigl(H_1(\Gamma;\Z) \stackrel{}{\lra}  H_1(S;\Z) \Bigr). 
\ee
\be \la{10.19.14.1a}
\Lambda^\circ_\Gamma:= {\rm Ker}\Bigl(H_1(\Gamma;\Z) \stackrel{}{\lra}  H_1(S^\times;\Z) \Bigr). 
\ee
The inclusion $S^\times \subset S$ induces a map $H_1(S^\times;\Z)\lra H_1(S;\Z)$. So there is a canonical embedding
$$
\Lambda^\circ_\Gamma \subset \Lambda_\Gamma.
$$

Recall the spectral surface $\Sigma_\Gamma$  of $\Gamma$. 
It is the surface associated with the conjugate bipartite ribbon graph $\Gamma^*$. 
Therefore the surface $\Sigma_{\Gamma}$ is homotopy equivalent to  
$\Gamma^*$. 
Since the abstract graphs $\Gamma^*$ and $\Gamma$ coincide, 
$\Sigma_{\Gamma}$ is homotopy equivalent to $\Gamma$. 
Filling the faces on $\Sigma_\Gamma$ by discs, we get the completed spectral surface ${\bf \Sigma}_\Gamma$. 
The inclusion $\Sigma_{\Gamma}  \hra {\bf \Sigma}_{\Gamma} $ 
provides an exact sequence
\be \la{1.21.12.101}
0 \lra Z_\Gamma \lra H_1(\Sigma_{\Gamma} ; \Z) 
\lra H_1({\bf \Sigma}_{\Gamma} ; \Z)\lra 0.
\ee
The subgroup $Z_\Gamma$ is generated by the zig-zag loops  of $\Gamma$, which are the same as the face 
loops of  $\Sigma_{\Gamma} $. They satisfy a single relation: the sum of all generators equals to zero. 
The intersection form on $H_1({\bf \Sigma}_{\Gamma} ; \Z)$ provides a skewsymmetric 
bilinear form  with the kernel $Z_\Gamma$:
\be \la{1.21.12.100sa}
H_1(\Sigma_{\Gamma} ; \Z)\wedge H_1(\Sigma_{\Gamma} ; \Z) \lra \Z.
\ee
Since $\Gamma$ is homotopy equivalent to $\Sigma_\Gamma$, we get a form on $H_1(\Gamma; \Z)$. 
Restricting it to the $\Lambda_\Gamma$ we 
get a skewsymmetric 
bilinear form 
\be \la{1.21.12.100s123}
(\ast, \ast)_\Gamma: \Lambda_\Gamma \wedge \Lambda_\Gamma \lra \Z.
\ee

Let $\Z[X]$ be the free abelian group generated 
by a set $X$. It has a canonical basis $\{x\}$, $x\in X$. 
 
Since we assumed that 
the faces of $\Gamma$ are 
discs on $S$, there are exact sequences
\be \la{1.21.12.101sds}
0 \lra H_2(S, \Z) \stackrel{r}{\lra}  \Z[\{\mbox{\rm Faces of $\Gamma$}\}] \lra H_1(\Gamma;\Z) 
\stackrel{}{\lra} H_1({S}; \Z)\lra 0.
\ee
\be \la{1.21.12.101sds}
0 \lra   \Z[\{\mbox{\rm Internal faces of $\Gamma$}\}] \lra H_1(\Gamma;\Z) 
\stackrel{}{\lra} H_1({S^\times}; \Z)\lra 0.
\ee
So there are isomorphisms
$$
\Lambda^\circ_\Gamma = \Z[\{\mbox{\rm Internal faces  of $\Gamma$}\}], ~~~~
\Lambda_\Gamma = 
\frac{\Z[\{\mbox{\rm Faces  of $\Gamma$}\}] }{r(H_2(S, \Z))}. 
$$
 There is a basis $\{e_G\}$ of the lattice $\Lambda^\circ_\Gamma$ given by the internal faces $G$ of $\Gamma$,  
and a  collection of vectors $\{e_F\}$ given by all faces $F$ of $\Gamma$, 
satisfying a single relation
$$
\sum_{F: ~\mbox{\rm faces of $\Gamma$}}e_F=0.
$$
The quiver ${\bf q}_\Gamma$  assigned to a bipartite graph $\Gamma$ on $S$ is given by the triple 
$
\Bigl(\Lambda_\Gamma, \{e_F\}, \langle\ast, \ast\rangle_\Gamma\Bigr).
$ The vectors $e_G$ assigned to the internal faces of $\Gamma$ are the non-frozen basic vectors.

\paragraph{A quiver for a decorated surface with boundary.} 
Suppose now that a decorated 
surface $S$ has boundary components. Let $S_{\cal D}$ be the topological double 
of $S$.
 It is an oriented surface without boundary
 obtained by gluing the surface $S$ and its mirror $S^\circ$ along the matching boundary components. 
It inherits the special points from $S$. 

Let $\Gamma$ be any bipartite ribbon graph. Changing the color of each vertex, and 
the cyclic order of the halfedges at every vertex, we arrive 
at its mirror - a bipartite ribbon graph $\Gamma^\circ$. 

Let $\Gamma$ be a bipartite ribbon graph on $S$. So its external vertices are on the boundary, and they are 
 not colored. The mirror graph $\Gamma^\circ$ lives naturally on $S^\circ$. 
We glue the graphs $\Gamma$ and  $\Gamma^\circ$ along  
the matching external vertices into 
a bipartite graph $\Gamma_{\cal D}\subset S_{\cal D}$.  
Let $S^\circ_{\cal D}$ be the surface $S_{\cal D}$ with the opposite orientation. There is a 
canonical involutive isomorphism $\sigma: S_{\cal D} \lra S^\circ_{\cal D}$ flipping  $S$ and $S^\circ$. 
It induces an isomorphism of bipartite ribbon graphs 
$$
\sigma: \Gamma_{\cal D} \stackrel{\sim}{\lra} \Gamma^\circ_{\cal D}. 
$$
Notice that as abstract graphs, $\Gamma_{\cal D}$ coincides with $\Gamma^\circ_{\cal D}$. 
Moreover there is a canonical isomorphism 
of lattices $\Lambda_{\Gamma^\circ_{\cal D}}= \Lambda_{\Gamma_{\cal D}}$, and  the isomorphism $\sigma$ induces an involution 
 $$
\sigma_*: \Lambda_{\Gamma_{\cal D}} \stackrel{\sim}{\lra} \Lambda_{\Gamma^\circ_{\cal D}}= \Lambda_{\Gamma_{\cal D}}. 
$$
\bd
 Let $\Gamma$ be a bipartite graph $\Gamma$ on a decorated surface $S$ with boundary.  

Lattices $(\Lambda_{\Gamma}, \Lambda_{\Gamma}^\circ)$ are
 the coinvariants of the induced involution $\sigma_*$  on the lattices 
$(\Lambda_{\Gamma_{\cal D}}, \Lambda_{\Gamma_{\cal D}}^\circ)$:
$$
\Lambda_{\Gamma}:= \Bigl(\Lambda_{\Gamma_{\cal D}}\Bigr)_{\sigma_*}:= \Lambda_{\Gamma_{\cal D}}/(\sigma_*(l) -l), ~~~~
\Lambda^\circ_{\Gamma}:= \Bigl(\Lambda^\circ_{\Gamma_{\cal D}}\Bigr)_{\sigma_*}.
$$
\ed

\bl Let $\Gamma$ be a bipartite graph on a decorated surface $S$ with a boundary.  
 Then 
$$
\Lambda_{\Gamma} = \frac{\Z[\{\mbox{\rm Faces of $\Gamma$}\}]}{r(H_2(S_{\cal D};\Z))}, ~~~~
\Lambda^\circ_{\Gamma} = \Z[\{\mbox{\rm Internal faces of $\Gamma$}\}].
$$
\el

\begin{proof} The lattice $\Lambda_{\Gamma}$ has canonical vectors associated with the faces of 
$\Gamma$, internal or external. 
Namely, if $F$ is an internal face of $\Gamma$, it is also a face $F_{\cal D}$ of 
$\Gamma_{\cal D}$. If $F$ is an external face of $\Gamma$, gluing it to its mirror $F^\circ$ we get 
a single face $F_{\cal D}$  of $\Gamma_{\cal D}$. Projecting the vectors $e_{F_{\cal D}} \in 
\Lambda_{\Gamma_{\cal D}}$ to the coinvariants 
we get canonical vectors $e_F \in \Lambda_{\Gamma}$.  
\end{proof}

\bd \la{QUIVERG}
The quiver  assigned to a bipartite surface graph $\Gamma$ is given by 
$$
{\bf q}_\Gamma:= \Bigl(\Lambda_\Gamma, \{e_F\}, \langle\ast, \ast\rangle_\Gamma\Bigr).
$$
The vectors $\{e_G\}$ assigned to the internal faces $G$ of $\Gamma$ are the non-frozen basic vectors. 

\noindent
The non-frozen subquiver is given by 
the sublattice $\Lambda^\circ_\Gamma $ with a basis $\{e_G\}$ and the induced form:
$$
{\bf q}^\circ_\Gamma:= \Bigl(\Lambda^\circ_\Gamma, \{e_G\}, {\rm Res}_{\Lambda^\circ_\Gamma}\langle\ast, \ast\rangle_\Gamma\Bigr).
$$
\ed

\bl  {\rm \cite[Lemma 4.5]{GK}}. A two by two move $\Gamma \lra \Gamma'$ centered at an internal face $G$  
leads to a mutation of quivers $
{\bf q}_\Gamma \lra {\bf q}_{\Gamma'}
$ 
in the direction of the vector $e_G$.
\el
\paragraph{A combinatorial description of the form $\langle\ast, \ast\rangle_\Gamma$.} 
Let $F$ be a face of a bipartite  graph $\Gamma$ on a decorated surface $S$, which may have 
a boundary. It gives rise to a face $F_{\cal D}$ of the bipartite graph $\Gamma_{\cal D}$. 
Its boundary, oriented by the orientation of $S_{\cal D}$, 
provides a zig-zag loop $\partial F_{\cal D}$ 
with its canonical orientation on the conjugate surface $\Sigma_{\Gamma_{\cal D}}$.
Given a pair of faces $F, G$ of  $\Gamma$, 
the intersection index  
$\langle F, G\rangle_\Gamma$ of the 
loops $\partial F_{\cal D}$ and $\partial G_{\cal D}$ on the conjugate surface $\Sigma_{\Gamma_{\cal D}}$  
can be calculated as follows.  

Let $E$ be an edge of the graph $\Gamma$. 
We set $\langle F, G\rangle^{(E)} =0$ if $E$ is not a common edge of $F$ and $G$, 
$\langle F, G\rangle^{(E)}  =+1$ if 
the (clockwise on our pictures) 
orientation of $F$ induces $\bullet \to \circ$ orientation of $E$, and 
$\langle F, G\rangle^{(E)} =-1$ otherwise.   

\bl \la{Combl} {\rm \cite[Proposition 8.3]{GK}} One has
$$
\langle F, G\rangle_\Gamma= \sum_{\mbox{\rm $E$: edges of $\Gamma$}} \langle F, G\rangle^{(E)}.
$$
\el 
Let $\Gamma$ be a bipartite  graph on  a decorated surface $S$ with boundary. 
Consider a lattice
$$
\widetilde \Lambda_\Gamma:= \Z[\{\mbox{\rm Faces of $\Gamma$}\}].
$$
 It has a basis $\{e_F\}$ provided by the faces $F$ of $\Gamma$. 
The skew-symmetric bilinear form 
$\langle\ast, \ast\rangle_\Gamma$ is well defined on the lattice $\widetilde \Lambda_\Gamma$. 
We define a traditional quiver $\widetilde {\bf q}_\Gamma$ as a triple 
$$
\widetilde {\bf q}_\Gamma = \Bigl(\widetilde \Lambda_\Gamma, \{e_F\}, \langle\ast, \ast\rangle_\Gamma\Bigr).
$$
The quiver ${\bf q}_\Gamma$ from Definition \ref{QUIVERG} is its quotient by the 
rank one sublattice spanned by the vector 
 $\sum_F e_F$ in the kernel 
of the form $\langle\ast, \ast\rangle_\Gamma$. 
 External faces provide the frozen basis vectors.

\paragraph{Quivers for ${\rm A}_{m}$ / ${\rm A}^*_{m}$ -bipartite graphs.} Let $S$ be a decorated 
surface. 

Let $\Gamma$ be a graph on $S$. We assume that its vertices 
are disjoint from the special points of $S$. 
We say that a face $F$ of  $\Gamma$ is {\it non-special} if it does not contain 
a special point of $S$. 

Let $\Gamma_{{\rm A}_{m}}$ be an ${\rm A}_{m}$-bipartite graph on $S$, and  
$\Gamma_{{\rm A}^*_{m}}$ the associated ${\rm A}^*_{m}$-bipartite graph on $S$. 
Although the special points of $S$ form a subset of the set of $\circ$-vertices of 
$\Gamma_{{\rm A}_{m}}$, 
the  $\Gamma_{{\rm A}^*_{m}}$ is  a graph on the punctured surface $S^\times$. 
Moreover, one has 
\be \la{11.7.14.1}
\{\mbox{\rm Faces of the graph $\Gamma_{{\rm A}_{m}}$}\} = \{\mbox{\rm Non-special faces of the graph $\Gamma_{{\rm A}^*_{m}}$}\}. 
\ee

It is clear from the way we glue the spectral surface from the discs associated with the zig-zag loops 
that there is a canonical isomorphism of surfaces 
$$
{\bf \Sigma}_{\Gamma_{{\rm A}_{m}}} = {\bf \Sigma}_{\Gamma_{{\rm A}^*_{m}}}. 
$$

\bl Let $S$ be a decorated surface without boundary. Then there are canonical  isomorphisms of lattices
\be \la{LATT}
\Z[\{\mbox{\rm Faces of the graph $\Gamma_{{\rm A}_{m}}$}\}] = 
\Z[\{\mbox{\rm Non-special faces of the graph $\Gamma_{{\rm A}^*_{m}}$}\} ]= 
\ee
$$
{\rm Ker}\Bigl(H_1({\Gamma_{{\rm A}^*_{m}}};\Z) \stackrel{}{\lra}  H_1(S^\times;\Z) \Bigr) =
{\rm Ker}\Bigl(H_1(\Sigma_{\Gamma_{{\rm A}^*_{m}}};\Z) \stackrel{\pi_*}{\lra}  H_1(S^\times;\Z) \Bigr). 
$$

\el

\begin{proof} The first isomorphism follows from (\ref{11.7.14.1}). 
The second is obvious since $S$ has no boundary, and the 
external faces of ${\Gamma_{{\rm A}^*_{m}}}$ are the ones which have special points inside. 
The third follows since the embedding $\Gamma_{{\rm A}^*_{m}} \hra \Sigma_{\Gamma_{{\rm A}^*_{m}}}$ 
is a homotopy equivalence, and followed by the spectral map, induces the embedding 
$\Gamma_{{\rm A}^*_{m}} \hra S^\times$.  

\end{proof}

\bl 
There is an isomorphism of lattices, where  $\pi_{*}$ is  induced 
by the spectral map: 
\be \la{10.19.14.1}
{\rm Ker}\Bigl(H_1(\Sigma_{\Gamma};\Z) \stackrel{\pi_{*}}{\lra}  H_1(S;\Z) \Bigr) = 
{\rm Ker}\Bigl(H_1(\Gamma;\Z) \stackrel{}{\lra}  H_1(S;\Z) \Bigr). 
\ee
\el

\begin{proof} 
The conjugate surface for the bipartite surface graph $\Gamma$ 
is identified with its spectral surface $\Sigma_{\Gamma}$. 
The conjugate surface of any bipartite ribbon graph $\Gamma$  is homotopy equivalent to 
 $\Gamma$. So
$
H_1(\Gamma;\Z) = H_1(\Sigma_{\Gamma};\Z). 
$ 
Using this identification, the map $H_1({\Gamma};\Z)
\lra H_1({S}; \Z)$ is identified with the map $ H_1(\Sigma_{\Gamma};\Z) 
\stackrel{\pi_{*}}{\lra} H_1({S}; \Z)$. 
\end{proof}

\section{Ideal webs on polygons and configurations of (decorated) flags} \la{SSec4}
\subsection{Ideal webs on a  polygon} \la{sec3.1}

Consider an example when the decorated surface $S$ is a 
polygon $P$, whose vertices are the special points. 
Let us spell the definition of ideal webs in this case. 
Notice that in this case any ideal web is strict ideal - there are no nontrivial covers of the polygon.

\begin{figure}[ht]
\centerline{\epsfbox{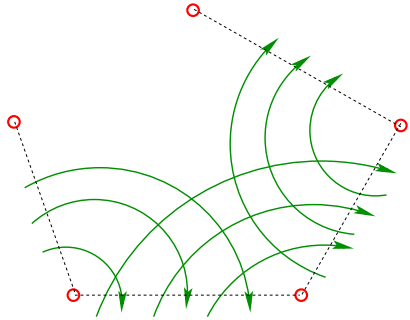}}
\caption{Strands of a rank 4 ideal web around three consecutive (red) vertices.}
\label{gra14}
\end{figure} 

\bd \la{2.4.12.111}
An ideal web in a polygon is a collection of oriented strands in the polygon 
with simple crossings, considered modulo isotopy, such that
\begin{enumerate}

\item As we move along a strand, 
the orientations of the crossing strands alternate.  

\item The web is connected, the complement to the web is a disjoint union of discs.

\item  
Each strand surrounds exactly one vertex of the polygon.

\item Strands have no selfcrossings; there are no parallel bigons.

\item Strands have no parallel halfbigons.

\end{enumerate}
\ed

Recall that two webs are {equivalent} if they are related by a sequence of two by two moves. 

\bt \la{8.13.11.1}
Any rank $m+1$ ideal web  on a polygon 
is either an ${\rm A}_{m}$-web, or a ${\rm A}_{m}^*$-web.  

Any two ${\rm A}_{m}$-webs are equivalent, and any two ${\rm A}_{m}^*$-webs are equivalent
\et

\begin{figure}[ht]
\centerline{\epsfbox{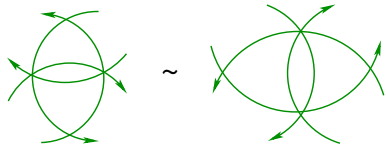}}
\caption{A two by two move of triple point diagrams.}
\label{gc3}
\end{figure}

{\bf Proof}.  Recall that resolving the  intersection points of a triple point diagram in one of the two 
standard ways, see Figure \ref{gc4}, we get 
a web, and this procedure can be reversed. The two by two moves of triple crossing diagrams, see Figure \ref{gc3},  match 
the ones of the webs. 

Therefore an ideal web in a polygon gives rise to  
a minimal triple point diagram. Each strand has 
the ``in'' endpoint, and the  ``out'' endpoint. The operations above do not change the matching 
of the ``in'' and ``out'' points on the boundary provided by the oriented strands, so that the strand is 
oriented from ``in'' to ``out''.

\begin{figure}[ht]
\centerline{\epsfbox{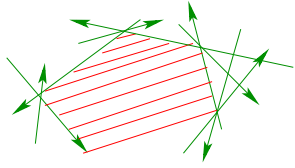}}
\caption{For any triple point diagram, the orientations of the strands induce 
consistent orientations on the complementary regions.}
\label{gc0}
\end{figure}

Let $n$ be the number of strands of a triple point diagram on the disc. 
The 
endpoints alternate as we move around the boundary. 
Indeed, the orientations of the strands induce 
consistent orientations on the complementary regions, see Figure \ref{gc0}, which implies the claim.

A triple point diagram on the disc is {\it minimal} if the number of its intersection points is not bigger 
than for any other diagram 
inducing  the same matching on the boundary. 

Given a collection of $2n$ alternatively oriented points on the boundary, 
a triple crossing diagram inducing such a collection provides a matching between 
the ``in'' and ``out'' points.  

\bt [\cite{T}] \la{DT}
1) In a disc with $2n$ endpoints on the 
boundary, all $n!$ matchings of ``in'' endpoints with ``out'' endpoints are achieved by 
minimal triple point diagrams.

2) 
Any two minimal triple point diagrams on the disc with the same matching on the 
endpoints are equivalent, that is can be related by a sequence of $2 \leftrightarrow 2$ moves, shown on Figure \ref{gc3}.
\et

Let us return to the proof of Theorem \ref{8.13.11.1}. 
Each strand 
goes clockwise around a vertex $v$ of the polygon.  Its ``in'' and ``out'' endpoints are at the 
sides of the polygon shared by $v$. 

\bl \la{12.2.12.1} Consider an ideal rank $m$ 
web on a polygon. Then each side of the polygon supports $2(m-1)$ alternating 
``out'' and  ``in'' endpoints. 
It  starts 
 from an ``out'' endpoint if we count counterclockwise, see Figure \ref{gra14}. 
\el 

\begin{proof}
Due to the minimality there are no selfintersections or 
parallel bigons formed by the strands surrounding a given vertex of the polygon. 
This plus condition 5)  implies (Corollary \ref{2.4.12.11}) that 
the strands surrounding a vertex $v$ of the polygon form ``concentric half-circles''. 
They start at the side preceding $v$, and ends on the next side, following counterclockwise orientation.

Recall that for any triple point diagram on the disc,  as we move along the boundary of the disc, 
the directions of the ends at the 
boundary alternate. 
Therefore for any web on the disc, 
the directions of the ends at the boundary alternate. 

 Finally, take the closest to $v$ strand surrounding a vertex $v$. It surround at least two endpoints. 
Indeed, if it surrounds no points, it can be contracted, which contradicts to the minimality. 
The only other option compatible with the above properties is that it surrounds 
two endpoints. \end{proof}

\bc
All rank $m$ ideal 
webs on a polygon have the same ``in'' and ``out'' matching of the ends, induced by the strands of the web. 
\ec

Therefore thanks to Theorem \ref{DT} the triple point crossing diagrams 
assigned to any two rank $m$ ideal webs on the polygon are equivalent. 
Thus there are exactly two equivalences classes of rank $m$ ideal webs on the polygon. 

The ideal ${\rm A}_{m}$-web on the disc assigned to a triangulation of the polygon
 has the configuration of ``in'' and ``out'' 
endpoints and their matching as described in Lemma \ref{12.2.12.1}. 
A two by two move of an ${\rm A}_{m}$-web is an ${\rm A}_{m}$-web. 
The same is true for the ideal ${\rm A}_{m}^*$-web on the disc assigned to a triangulation of the polygon. 
Theorem \ref{8.13.11.1} is proved.

\subsection{${\rm A}_{m}$ / ${\rm A}_{m}^*$-webs on polygons and cluster coordinates}

\subsubsection{ ${\rm A}_{m}$-webs and cluster coordinates on configurations 
of decorated flags }
Let $(V_m, \Omega_m)$ be an $m$-dimensional vector space with a volume form $\Omega_m$. 
\bd
i) A decorated flag 
in $(V_m, \Omega_m)$ 
is a data $(F_{\bullet}, \alpha_\bullet)$:
\be \la{1}
F_{\bullet}: ~ 0=F_0 \subset F_1  \subset \ldots \subset F_{m-1} \subset F_{m} = V_m, 
\quad {\rm dim}F_{p}=p, ~~ \alpha_p \in {\rm det}F_{p}-0,
\ee
consisting of a flag $F_{\bullet}$ and non-zero $p$-vectors   
$\alpha_p \in {\rm det}F_{p}$  for $p=1,..., m-1$. 

ii) A $\ast$-decorated flag 
in $(V_m, \Omega_m)$ is a similar data $(F_{\bullet}, \alpha_\bullet)$ where $\alpha_p \in {\rm det}F_{p}$  for $p=1,..., m$.
\ed

Instead of the $p$-vectors  $\alpha_p$ one can choose 
non-zero vectors $v_p \in F_{p}/F_{p-1}$,  $p=1,..., m-1$. 
Indeed, given the vectors, the volume forms are $\alpha_p = v_1 \wedge ... \wedge v_p$.

Denote by ${\cal A}_{SL_m}$ 
the moduli space of decorated flags in $(V_m, \Omega_m)$, and 
by ${\rm Conf}_n({\cal A}_{SL_m})$ the moduli space of  
configurations of $n$ decorated flags:
$$
 {\rm Conf}_n({\cal A}_{SL_m}):=  ({\cal A}_{SL_m})^n/{\rm Aut}(V_m, \Omega_m).
$$

Take an oriented $n$-gon $P_n$ whose vertices are numbered 
by $(1, ..., n)$ following the clockwise orientation. 
Consider configurations of decorated flags $(A_1, ..., A_n)$ in  $(V_m, \Omega_m)$
 labeled by the  vertices of $P_n$. 
Let ${\cal W}$ be an ${\rm A}_{m-1}$-web on the polygon $P_n$. 
Given a face $F$,  let $d_i$ be the codistance from 
$F$ to the vertex $v_i$. 
Corollary \ref{1.21.12.1d} tells that  
\be
d_1+\ldots +d_n=m.
\ee 
Given an decorated flag $A= (F_\bullet, \alpha_\bullet)$, see 
(\ref{1}), we use the convention 
\be \la{conv}
\alpha_0(A) =1\in \C = {\rm det}F_0.
\ee
So $\alpha_0(A)$
 is the unit in the exterior algebra of $V_m$. 
\bd \la{6.8.11.2} Given an ${\rm A}_{m-1}$-web ${\cal W}$ on 
$P_n$, we assign to its face $F$ a regular function 
$A^{\cal W}_F$ on the 
space ${\rm Conf}_n({\cal A}_{SL_m})$ given by 
\be \la{6.8.11.3a}
A^{\cal W}_F:= 
\langle \Omega_m, \alpha_{d_1}(A_{1})\wedge \ldots \wedge \alpha_{d_n}(A_{n})\rangle.
\ee
\ed

\bl
The functions $A^{\cal W}_F$ are invariant under the twisted cyclic shift 
$$
t: {\rm Conf}_n({\cal A}_{SL_m}) \lra {\rm Conf}_n({\cal A}_{SL_m}), \quad 
(A_1, ..., A_n) \lms ((-1)^{m-1}A_n, A_1, ..., A_{n-1}).
$$
\el

\begin{proof} Straightforward. \end{proof}

\subsubsection{${\rm A}_{m}^*$-webs  and cluster coordinates on configurations of $\ast$-decorated flags}
Recall the configuration space of $n$ $\ast$-decorated flags for the group $SL_m = {\rm Aut}(V_m, \Omega_m)$:
\be \la{2.22.12.1}
{\rm Conf}_n({\cal A}^*_{SL_m}):= \Bigl({\cal A}_{GL_m} \times \ldots \times 
{\cal A}_{GL_m}\Bigr)/{\rm Aut}(V_m, \Omega_m), ~~~~
{\cal A}_{GL_m} := GL_m/U.
\ee

Given an ${\rm A}_{m-1}^*$-web ${\cal W}_{{\rm A}_{m-1}^*}$ on an $n$-gon, 
just the same construction as above provides
 regular  functions $A^{\cal W}_{F}$ on the space (\ref{2.22.12.1}). 
We assign the coordinates to \underline{all faces} of the ${\rm A}_{m-1}^*$-web ${\cal W}_{{\rm A}_{m-1}^*}$. 
The  number of the faces equals the dimension of the space 
(\ref{2.22.12.1}). 

The coordinates at the
 \underline{internal} faces of ${\cal W}_{{\rm A}_{m-1}^*}$ are nothing else but the coordinates assigned to the associated 
${\rm A}_{m-1}$-web ${\cal W}_{{\rm A}_{m-1}}$: this web is 
obtained by shrinking all 2-valent $\bullet$-vertices, and creating as a 
result a $\circ$-vertex at each vertex of the polygon.

Finally, one has 
\be \la{2.22.12.1}
{\rm Conf}_n({\cal A}_{GL_m}):= \Bigl({\cal A}_{GL_m} \times \ldots \times 
{\cal A}_{GL_m}\Bigr)/{\rm Aut}(V_m) = {\rm Conf}_n({\cal A}^*_{SL_m})/GL(1).
\ee
Here the group $GL(1)$ acts diagonally. The action of an element $t\in GL(1)$ 
amounts to multiplying by $t$ each of the coordinates $\{A^{\cal W}_{F}\}$ 
at the faces of an ${\rm A}_{m-1}^*$-web ${\cal W}_{{\rm A}_{m-1}^*}$.

\subsubsection{${\rm A}_{m}$-webs and cluster Poisson coordinates on configurations 
of  flags}

Let ${\cal B}_m$ be 
the space of flags in $V_m$. 
Consider   the moduli space of  
configurations of $n$ flags:
$$
{\rm Conf}_n({\cal B}_m):= ({\cal B}_m)^n/{\rm Aut}(V_m).
$$

Let $F$, $G$ be two faces of an ${\rm A}_{m-1}$-web ${\cal W}$ on a polygon $P_n$.
 Recall the pairing 
$\langle F, G\rangle_{\cal W}$.

There is a projection ${\cal A}_{SL_m}\to {\cal B}_m$. 
Given a configuration of flags $(F_1, ..., F_n)$, 
choose a configuration $(A_1, ..., A_n)$ of decorated flags 
projecting onto  it. 
Then for any face $G$ of an ideal web ${\cal W}$ we have a function 
$A_G^{\cal W}$. 

\bd \la{6.8.11.2a} Given an ${\rm A}_{m-1}$-web ${\cal W}$ on 
$P_n$, we assign to an $\underline{\mbox{{internal}}}$ face $F$ of ${\cal W}$ a rational function
 $X^{\cal W}_F$ on the space ${\rm Conf}_n({\cal B}_m)$, given by the following product over the faces $G$: 
\be \la{6.8.11.4a}
X^{\cal W}_F:= \prod_{G}{(A_G^{\cal W}})^{\langle F, G\rangle_{\cal W}}.
\ee
\ed

\bl
The function $X^{\cal W}_F$  does not depend 
on the choice of configuration $(A_1, ..., A_n)$. 
\el

\begin{figure}[ht]
\centerline{\epsfbox{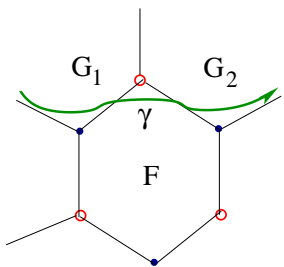}}
\caption{Contribution to $X_F$ of the decorated flag associated with a strand $\gamma$.}
\label{gc10}
\end{figure}

\begin{proof} Take a strand $\gamma$. It is associated with a special point $s$. 
Let $d$ be the codistance from $\gamma$ to $s$. Let us calculate how the coordinate 
$X_F$ changes when we multiply $\alpha_d$  by a non-zero number $\lambda$. 
Clearly $X_F$ might change only if $\gamma$ intersects some edges of the face $F$. 
In this case, it must intersect a pair of consecutive edges $E_1, E_2$ sharing a vertex $v$, or may be several such pairs. 
Let us investigate what happens near $v$. There is a unique  face $G_i$ containing the edge $E_i$ and different then $F$, 
see Figure \ref{gc10}. The faces $G_1, G_2$ share the vertex $v$. They are on the same side of the zig-zag strand $\gamma$. 
Thus their contribution to $X_F$, using notation from Lemma \ref{Combl},  is
$$
(A_{G_1})^{\langle F, G_1\rangle^{(E_1)}}(A_{G_2})^{\langle F, G_1\rangle^{(E_2)}} = (A_{G_1}/A_{G_2})^{\pm 1}. 
$$
Therefore it does no change if we multiply $\alpha_d$ by $\lambda$. \end{proof}

\subsection{Cluster nature of the coordinates assigned to   ${\rm A}_{m}$ / ${\rm A}^*_{m}$-webs }

Recall that what we call  a cluster / cluster Poisson coordinate system was called in \cite{FG2} 
cluster ${\cal A}$ /  cluster ${\cal X}$- coordinate system.  
A collection of cluster / cluster Poisson coordinate systems related by cluster / cluster Poisson 
transformations form a cluster / cluster Poisson atlas. 

\bt \la{MTHW}

a)   ${\rm A}_{m-1}$-webs on an $n$-gon provide a cluster atlas  on 
 ${\rm Conf}_n({\cal A}_{SL_m})$. Precisely:
\begin{enumerate} 
\item
The functions $\{A^{\cal W}_F\}$ at the faces $F$ of a given ${\rm A}_{m-1}$-web ${\cal W}$ 
form a regular coordinate system on  the space ${\rm Conf}_n({\cal A}_{SL_m})$.

\item  Any two ${\rm A}_{m-1}$-webs are related by two by two and shrink / expand moves. 

Any two by two move of webs  $\mu_F: {\cal W} \to {\cal W}'$ centered at a face $F$  of ${\cal W}$  
amounts to a  cluster 
mutation at the face $F$. The shrink / expand moves do not affect the coordinates. 
\end{enumerate}

b)  Internal faces of ${\rm A}_{m-1}$-webs on an $n$-gon provide a cluster Poisson atlas  on 
 ${\rm Conf}_n({\cal B}_m)$:

\begin{enumerate} 

\item 
The functions $\{X^{\cal W}_F\}$ assigned to  the internal faces $F$ of a given ${\rm A}_{m-1}$-web ${\cal W}$  
form a rational coordinate system on  the space ${\rm Conf}_n({\cal B}_m)$.

\item Any two by two move $\mu_F: {\cal W} \to {\cal W}'$    
amounts to a cluster 
Poisson mutation at the face $F$. The shrink / expand moves do not affect the coordinates. 
\end{enumerate}

c)  The ${\rm A}^*_{m-1}$-webs on an $n$-gon provide 
 a cluster atlas  on 
 ${\rm Conf}_n({\cal A}^*_{SL_m})$:

\vskip 2mm
The functions $\{A^{\cal W}_F\}$ at the faces $F$  of a given ${\rm A}^*_{m-1}$-web ${\cal W}$ 
form a regular coordinate system on  ${\rm Conf}_n({\cal A}^*_{SL_m})$. 
 A two by two move  $\mu_F: {\cal W} \to {\cal W}'$ 
amounts to a  cluster 
mutation at  $F$. 

\vskip 3mm

d) Given an   ${\rm A}^*_{m-1}$-web   ${\cal W}$ on an $n$-gon, the 
 functions $\{A_F^{\cal W}\}$ provide  
 a birational isomorphism
$$
\Lambda_{\cal W}\otimes_\Z {\Bbb G}_m \lra {\rm Conf}_n({\cal A}_{GL_m}). 
$$
\et

\begin{proof} Our construction recovers the cluster coordinates 
introduced in \cite[Section 8]{FG1}. 

\bl \la{6.8.11.5a} Let ${\cal W}_T$ be the ${\rm A}_{m}$-web assigned to an 
ideal triangulation $T$ of the polygon. Then

i) For the web ${\cal W}_T$, the functions from Definitions 
\ref{6.8.11.2} and \ref{6.8.11.2a} coincide with the coordinates assigned to $T$
in {\it loc. cit.}

ii) The form $\langle F, G\rangle$ for the web ${\cal W}_T$ coincides with the form 
$\varepsilon_{ij}$ on the set $I$ parametrising the coordinates 
considered in {\it loc. cit.}. 

iii) A two by two move of the web ${\cal W}_T$ corresponds to a mutation in {\it loc. cit.}. 
\el

\begin{proof}  i) For the  cluster coordinates this is clear from the very definition. Then 
formula (\ref{6.8.11.4}) is the standard formula relating 
the  cluster  and cluster Poisson coordinates. 

ii) It is clear from the very definitions. 

iii) Follows from Section \ref{ssec2.1}.  \end{proof} 

a) We proved in \cite{FG1} that the coordinates assigned 
to an ideal triangulation $T$ of the polygon
 are indeed coordinates on the moduli space ${\rm Conf}_n({\cal A}_{SL_m})$. 
The part i) of Lemma \ref{6.8.11.5a} tells that those coordinates coincide with the ones arising from the web assigned to $T$. 
Finally and two ${\rm A}_{m-1}$-webs on the polygon are related by two by two moves by Theorem \ref{8.13.11.1}. 
 Thanks to the Pl\"ucker identity between the minors of a $2 \times 4$ matrix, and similar identities between 
higher order minors induced by it,  
a two by two move amounts to a cluster mutation of the ${\cal A}$-coordinates. The part a) is proved.

b) Follows from a) and Lemma \ref{6.8.11.5a} by using the standard relationship between the 
  and cluster Poisson coordinates, see (\ref{6.8.11.4a}).  

c) Deduced easily from  a). 

d) It is deduced from c). Indeed, 
${\rm Conf}_n({\cal A}_{GL_m}) = {\rm Conf}_n({\cal A}^*_{SL_m})/{\Bbb G}_m$, and on the other hand 
the elements $e_F$ generate the lattice $\Lambda_{\cal W}$ and satisfy a single relation 
$\sum e_F=0$. \end{proof}

\section{Ideal webs and coordinates on moduli spaces of local systems} \la{SSec5}

\subsection{Cluster coordinates for moduli spaces of local systems on surfaces }
Denote by  $T'S$ the bundle of non-zero tangent vectors to the surface $S$. 
A {\it twisted $SL_m$-local system} on $S$ is a local system on $T'S$ which has monodromy $(-1)^{m-1}$ around 
a loop rotating a tangent vector at any point by $360^\circ$. 

Let $G$ be either $SL_m$ or $GL_m$.  Let 
$U$ be a maximal unipotent subgroup of $G$. Let ${\cal L}$  be a $G$-local system on a space. Set 
$$
{\cal L}_{\cal A}:= {\cal L}/U. 
$$
When ${\cal L}$  is an $SL_m$-local system on a space, we also set 
$$
{\cal L}_{\cal A^*}:= {\cal L} \times_{SL_m} (GL_m/U). 
$$
Given a 
$PGL_m$-local system ${\cal L}$ on $S$, set ${\cal L}_{\cal B}:= {\cal L}/B$, where 
$B$ is a Borel subgroup of $PGL_m$. 

We remove a little disc around each special point  
of a decorated surface $S$, getting a surface $S^\circ$. 
Let $\partial_{\rm th} S^\circ$ be the {\it thickened boundary of $S^\circ$} - 
a small neighborhood of the boundary of  $S^\circ$.

\begin{definition} [\cite{FG1}] Let $S$ be a decorated surface. 
\begin{itemize}

\item Let $G$ be either $SL_m$ or $GL_m$. A decoration on a twisted $G$-local system ${\cal L}$ on $S^\circ$  is 
a flat section of restriction of ${\cal L}_{\cal A}$ to $T'(\partial_{\rm th} S^\circ)$. 
The moduli space ${\cal A}_{G,S}$ parametrises decorated twisted $G$-local system on $S^\circ$.

\item A $\ast$-decoration on a twisted $SL_m$-local system ${\cal L}$ on $S^\circ$  is 
a flat section of the restriction of  ${\cal L}_{\cal A^*}$ to  $T'(\partial_{\rm th} S^\circ)$. 
The  moduli space ${\cal A}^*_{SL_m,S}$ parametrises $\ast$-decorated twisted $SL_m$-local systems on $S^\circ$. 

\item  
A framing on a $PGL_m$-local system ${\cal L}$ on $S^\circ$ is 
a  flat section of  the restriction of 
 ${\cal L}_{\cal B}$ to  $\partial_{\rm th} S^\circ$. 
The moduli space ${\cal X}_{PGL_m,S}$ parametrises framed $G$-local systems on $S^\circ$. 
\end{itemize}
\end{definition}

Let ${\cal W}$ be an ${\rm A}_{m-1}$-web on a decorated surface $S$, and 
$\widetilde {\cal W}$ its lift to the universal cover $\widetilde S$. 
We are going to assign the cooridnates to the 
faces / internal faces 
of $\widetilde {\cal W}$. By the construction, the coordinates assigned to the faces which differ 
by the action of the deck transformation group $\pi_1(S)$ will be the same. 
This way we assign the coordinates to the faces of ${\cal W}$. 

If ${\cal W}$ is a strict ${\rm A}_{m-1}$-web, there is no need to go to the universal cover, 
and the coordinates are assigned just to the faces / internal faces of the web ${\cal W}$.

Let $\widetilde {\cal L}$ be the lift of a decorated $SL_m$-local system ${\cal L}$ 
to $\widetilde S$. 
Let $s$ be a special point on $\widetilde S$. 
Denote by $A(s)$ the 
flat section of $\widetilde {\cal L}_{\cal A}$ providing the decoration  $\widetilde {\cal L}$.
Given a zig-zag strand $\gamma$ of $\widetilde {\cal W}$ associated to $s$, we can extend 
$A(s)$ to a flat section of $\widetilde {\cal L}_{\cal A}$ 
on the punctured disc $\widetilde S_{\gamma}-s$.

Let $F$ be a face of $\widetilde {\cal W}$. 
Choose a non-zero tangent vector $v$ at a point inside of  $F$. 
 Connect it by a path in $T'\widetilde S$ with a non-zero tangent vector $v_i$ at
 a point near the special point $s_i$. Transport the section $A(s_i)$ along this path 
from $v_i$ to $v$. 
The resulting decorated flag $A_v(s_i)$ over $v$  is well defined 
since $A(s_i)$ is a flat section of $\widetilde {\cal L}_{\cal A}$  over $\widetilde S_\gamma-s_i$. 
We get decorated flags  $A_v(s_1), \ldots , A_v(s_n)$ 
over $v$. 
Given a face $F$,  let $d_{s_i}$ be the codistance from 
$F$ to  $s_i$. 
Corollary \ref{1.21.12.1d} implies that  
\be
d_{s_1}+\ldots +d_{s_n}=m.
\ee

\bd \la{6.8.11.2bz} Given an ${\rm A}_{m-1}$-web ${\cal W}$ on 
$S$ and a face $F$ of $\widetilde {\cal W}$, we define 
a regular function $A^{\cal W}_F$ on the moduli space ${\cal A}_{SL_m, S}$ by setting
\be \la{6.8.11.3b}
A^{\cal W}_F:= 
\langle \Omega_m, \alpha_{d_1}(A_v(s_1))\wedge \ldots \wedge \alpha_{d_n}(A_v(s_n))\rangle.
\ee
\ed

The functions $A^{\cal W}_F$ are evidently $\pi_1(S)$-invariant, and so 
we may assume that $F$ stands here for the $\pi_1(S)$-orbit in the space of faces of 
$\widetilde {\cal W}$. 

\vskip 3mm
Just the same way, given an ${\rm A}^*_{m-1}$-web ${\cal W}$, we define regular  functions $A^{\cal W}_F$ on the moduli space ${\cal A}^*_{SL_m, S}$ by formula (\ref{6.8.11.3b}). The only difference is that 
now we have extra coordinates at the external faces $F$ containing special points of $S$. 
\vskip 3mm

\bd \la{6.8.11.a2bz}Given an ${\rm A}_{m-1}$-web ${\cal W}$ on 
$S$ and an internal face $F$ of $\widetilde {\cal W}$, we define 
a rational function  $X^{\cal W}_F$ on the space ${\cal X}_{PGL_m, S}$ by the following product over the faces $G$: 
\be \la{6.8.11.4}
X^{\cal W}_F:= \prod_{G}{(A_G^{\cal W}})^{\langle F, G\rangle}, 
\ee
\ed

The analog of Lemma \ref{6.8.11.5a} is valid. In particular:
\bl For the ${\rm A}_{m-1}$-web  ${\cal W}_T$ assigned to an 
ideal triangulation $T$ of $S$, the coordinates from Definitions  
\ref{6.8.11.2bz} and \ref{6.8.11.a2bz} coincide with the coordinates assigned to $T$
in \cite{FG1}, Section 8. 
\el

\bt \la{8.13.11.2} Given a strict ideal ${\rm A}_{m-1}$-web ${\cal W}$ on $S$:

i) The regular functions $\{A_F^{\cal W}\}$ at the faces $F$ of ${\cal W}$ 
  are cluster coordinates on 
 ${\cal A}_{SL_m, S}$ 

ii) Rational functions $\{X_F^{\cal W}\}$ at the faces $F$ of ${\cal W}$ 
are cluster Poisson coordinates on ${\cal X}_{PGL_m, S}$. 

\vskip 2mm
Given a strict ideal  ${\rm A}^*_{m-1}$-web ${\cal W}$ on $S$: 

iii) The regular functions $\{A_F^{{\cal W}}\}$ at the faces $F$ of ${\cal W}$ 
  are cluster coordinates on 
${\cal A}^*_{SL_m, S}$. 
\et

\begin{proof} Given a strict ideal web ${\cal W}$, 
and picking an ideal triangulation $T$ of $S$, the restriction ${\cal W}_t$ of ${\cal W}$ 
to each triangle $t$ of $T$ is 
a strict ideal web on $t$. Therefore by Theorem \ref{MTHW} it can be transformed by elementary transformations to the 
standard web of type  either ${\rm A}_{m-1}$ or ${\rm A}^*_{m-1}$ on $t$. So 
the original web ${\cal W}$ is equivalent to one of the two standard webs ${\cal W}_{{\rm A}_{m-1}, T}$ 
or ${\cal W}_{{\rm A}^*_{m-1}, T}$ related to the triangulation $T$, depending on the type of the web ${\cal W}$.

It was proved in \cite{FG1} that the ${\cal A}$-functions assigned to 
an ideal triangulation $T$ form a regular cluster coordinate system on the  space 
${\cal A}_{SL_m, S}$, and the ${\cal X}$-functions form a rational 
cluster Poisson coordinate system on the  space 
${\cal X}_{PGL_m, S}$. 

A two by two move amounts to a cluster transformation of the ${\cal A}$-coordinates. 
This implies  that a two by two move amounts to a cluster Poisson transformation of the ${\cal X}$-coordinates. 
So we get the parts i) and ii). 
The part iii) reduces easily to i). 
\end{proof}

\paragraph{Remark.} Any spectral web ${\cal W}$ on a decorated surface $S$ satisfying the condition that 
\be \la{3.4.12.1}
\mbox{each 
domain $S_\gamma$ is a disc}
\ee
produces a collection of functions $\{A_F\}$ assigned to the faces $F$ of ${\cal W}$ 
as well as a collection of functions  $\{X_F\}$ assigned to the internal faces $F$ of ${\cal W}$.  
Indeed, the decoration / framing provides a flat section on the disc $S_\gamma$ / punctured disc 
$S_\gamma$. The counting argument related to the degree of the spectral map,  guarantees that we can 
use this flat sections to define the functions. 
 
However the number of faces / internal faces of such a web is typically 
bigger then the dimension of the spaces ${\cal A}_{SL_m, S}$ and, respectively,  
${\cal X}_{PGL_m, S}$. So non-ideal spectral webs satisfying condition (\ref{3.4.12.1})
usually do not produce coordinate systems on these spaces.

\paragraph{${\rm A}_{m-1}^*$-bipartite graphs and ${\cal A}_{GL_m, S}$.}

\bt \la{2.22.12.4} 
Let $\Gamma$ be a strict ${\rm A}_{m-1}^*$-graph on a decorated surface $S$ without boundary. 
Then there is a birational isomorphism
$$
{\cal A}_{GL_m, S} \stackrel{\sim}{\lra} ~ \mbox{\rm the moduli space of 
$GL(1)$-local system on  $\Gamma$}. 
$$
\et

\begin{proof} Theorem \ref{2.22.12.4} is equivalent to Theorem  \ref{8.13.11.2} iii). 
Indeed, one has 
$$
0 \lra \frac{\Z[\{\mbox{\rm faces of $\Gamma$}\}]}{H_2(S, \Z)} \lra H_1(\Gamma, \Z) \lra H_1(S, \Z) \lra 0. 
$$
Thus a line bundle with connection on the graph $\Gamma$ is 
 uniquely determined by the following data:

i) Monodromies around  the faces of the graph $\Gamma$.

ii) The global monodromies 
on $S$. 

Since $GL_m=SL_m \times {\Bbb G}_m$, and thanks to (\ref{2.22.12.1}), 
the space ${\cal A}_{GL_m, S}$ is a product:
\be \la{12.12.12}
{\cal A}_{GL_m, S} = {\cal A}^*_{SL_m, S}/GL(1) ~\times~ \mbox{\rm the moduli space of 
$GL(1)$-local system on  $S$}.
\ee
 
The data ii) describes a $GL(1)$-local system on  $S$. It is responsible for the second factor. 

The data i) accounts for the description  of the space ${\cal A}^*_{SL_m, S}$  
in Theorem  \ref{8.13.11.2} iii) by the coordinates on the  faces of $\Gamma$.  
The quotient by 
the  $GL(1)$ in the first factor in (\ref{12.12.12}) amounts to the quotient 
of $\Z[\{\mbox{\rm faces of $\Gamma$}\}]$ by $H_2(S, \Z) = \Z$. 
\end{proof} 

\subsection{Ideal loops and monodromies of  framed local systems around punctures} 
Recall that, given a decorated surface $S$, 
 a {\it puncture}  is an internal special point on $S$. Recall 
$$
S^\times:= S - \{\mbox{\rm punctures of $S$}\}.
$$ 
The monodromy of a framed $PGL_m$-local system on $S^\times$ around a puncture     
lies in the Cartan group $H$ of $PGL_m$, isomorphic to ${\Bbb G}_m^{m-1}$. 
Monodromies around the punctures provide a map
\be \la{6.8.11.1}
 {\cal X}_{PGL_m, S} \lra H^{\{\mbox{punctures    of $S^\times$}\}}.
\ee
Its fibers are the generic symplectic leaves on ${\cal X}_{PGL_m, S}$. 
Monodromies around the punctures generate the center of the Poisson algebra 
of functions on ${\cal X}_{PGL_m, S}$ \cite{FG1}. 

To each 
zig-zag loop $\gamma$ we assign a function $C_\gamma$ 
given by the  product of the face coordinates assigned to 
the faces $F$ sitting inside of the disc $S_\gamma$:
$$
C_\gamma:= \prod_{F \subset S_\gamma}X_F
$$
Given a strict ${\rm A}_{m-1}$-web ${\cal W}$ on $S$, for every puncture $s$
there are $m-1$ ideal loops around $s$. 

\bt Let ${\cal W}$ be a strict ${\rm A}_{m-1}$-web on $S$. Then 

i) The functions $C_\gamma$ are invariant under the two by two moves. 

ii) The functions $C_\gamma$ are Casimirs. 
Their product is  equal to $1$, and this is the only relation between them. 
They generate the center of the Poisson algebra of functions 
on ${\cal X}_{PGL_m, S}$. 

iii) The map (\ref{6.8.11.1}) is described by the collection of functions $C_\gamma$. 
\et

\begin{proof}  i) Indeed, under a two by two moves the change of the face coordinates 
located inside of the zig-zag strand shown on Fig \ref{gc9} is given by 
$$
A'= A(1+X^{-1})^{-1}, ~~~ B'=B(1+X). 
$$
Thus $ABX= A'B'$, which proves the  claim.

\begin{figure}[ht]
\centerline{\epsfbox{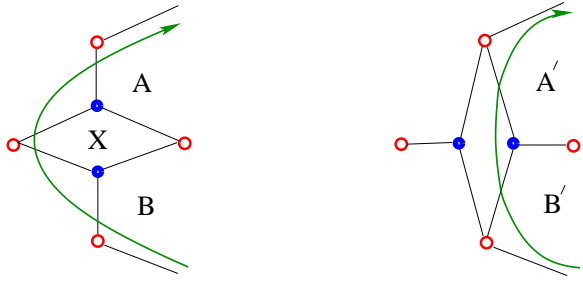}}
\caption{The action of a two by two move on the face coordinates.}
\label{gc9}
\end{figure} 

ii) The zig-zag loops match the
 boundary components on the spectral surface $\Sigma$. This plus the exact sequence (\ref{1.21.12.101}) 
implies all the claims since the Poisson bracket 
is given by the formula
$$
\{X_{F}, X_G\}= \langle F, G\rangle X_FX_G.
$$

iii) For the ${\rm A}_{m-1}$-web related to a triangulation of $S$ this 
is clear 
from the explicit description of the framed local system assigned to the coordinates 
which was given in \cite{FG1}.  So in general it follows from i). \end{proof}

%



\section{Webs, quivers with potentials, and 3d Calabi-Yau categories} \la{sssec2}

\subsection{Quivers  with potentials and  ${\rm CY}_3$ categories with cluster collections} \la{sec1.4}

In Section \ref{sec1.4} we collect for convenience of the reader 
some foundational material, including the relationship between quivers with potentials and 
3d Calabi-Yau categories \cite[Section 8.1]{KS1}. 

\paragraph{1. Quivers and potentials.} 
A {quiver}  ${\bf q}$ without loops and 2-cycles is described by a set of vertices $I$, a subset 
$I^f\subset I$ of frozen vertices, 
and a skewsymmetric function $\varepsilon_{ij}: I\times I \to \Z$, encoding the number of oriented 
arrows $i\to j$ between the vertices. 

Any 
element $k \in I-I^f$ provides  a {mutated in the 
direction $k$} quiver $ {\mathbf q}'$ described by the same set $I$ and 
new function $ \varepsilon'_{ij} = \mu_k(\varepsilon_{ij})$ defined by the Fomin-Zelevinsky formula \cite{FZI}:
\begin{equation} \label{5.11.03.6}
 \varepsilon'_{ij} := \left\{ \begin{array}{lll} 
- \varepsilon_{ij} & \mbox{ if $k \in \{i,j\}$} \\ 
\varepsilon_{ij} & \mbox{ if $\varepsilon_{ik}
\varepsilon_{kj} \leq 0, \quad k \not \in \{i,j\}$} \\
\varepsilon_{ij} + |\varepsilon_{ik}| \cdot \varepsilon_{ kj}& 
\mbox{ if $\varepsilon_{ik}
\varepsilon_{kj} > 0, \quad k \not \in \{i,j\}.$}\end{array}\right.
\end{equation} 
This procedure is involutive: 
the mutation of  $\varepsilon'_{ij}$ at the vertex $k$ is 
the original function $\varepsilon_{ij}$.

There is  a geometric description of a quiver   as  a data 
$\Bigl(\Lambda, \Lambda^f, \{e_i\}, (\ast, \ast)\Bigr)$,  given by  a lattice 
$\Lambda$, a 
basis $\{e_i\}_{i\in I}$ of $\Lambda$, 
a skew-symmetric $\Z$-valued 
bilinear form $(\ast, \ast)$ on $\Lambda$, and 
 a sublattice $\Lambda^f$ generated by the {\it frozen basis vectors}. 
The non-frozen basis vectors generate a subquiver, called the non-frozen part of the original quiver. 

One relates two definitions by setting $\Lambda = \Z[I]$ and $(e_i, e_j) = \varepsilon_{ij}$.

A { mutation of a quiver} ${\bf q}$ 
in the direction of a basis vector $e_k$ is a new quiver $  {\bf q}' = \mu_{e_k}({\bf q})$.
 The lattice $\Lambda$ and the form $(\ast, \ast)$  for $  {\bf q}'$ 
are the same as for ${\bf q}$. 
The basis $\{  e'_i\}$ for $  {\bf q}'$ is given by 
\begin{equation} \label{12.12.04.2a}
 e'_i := 
\left\{ \begin{array}{lll} e_i + (e_{i}, e_k)_+e_k
& \mbox{ if } &  i\not = k\\
-e_k& \mbox{ if } &  i = k.\end{array}\right. \qquad a_+:={\rm max}(a, 0).
\end{equation}
The composition of mutations $\mu_{  e'_k} \circ \mu_{e_k} $ no longer acts as the 
identity on the original basis $\{e_k\}$.  It is rather the reflection of the 
collection $\{e_i\}$ at $e_k$. 

In the ``simply-laced case'' 
  a cluster algebra  \cite{FZI},  and more generally 
a cluster variety \cite{FG2},  are described by a collection of quivers  
related by quiver mutations (\ref{5.11.03.6}). \vskip 3mm

Categorification requires a more elaborate 
version of a quiver, which we use throughout the paper. 
A quiver $Q$ is  given by a 
set of vertices $\{i\}_{i\in I}$ and a collection of finite dimensional vector spaces  
$A(i,j)$, the {\it arrow spaces},  assigned to each pair of vertices $(i,j)$. 
It is customary to choose a basis in each arrow space 
${\rm A}(i,j)$, which we thought of  as a collection of arrows  
 $i \to j$.  

To define the potential we consider the path algebra ${\cal P}(Q)$ of the quiver $Q$:
$$
{\cal P}(Q):= \oplus_{n\geq 2}\bigoplus_{i_1, ..., i_n\in I}A(i_1, i_2) \otimes \ldots \otimes A(i_{n-1}, i_n).
$$
The product is given by concatenation. Let us consider its cyclic envelope 
\be \la{cycen}
{\rm HH}_0({\cal P}(Q)):=  \oplus_{n\geq 2}\bigoplus_{i_1, ..., i_n\in I}\Bigl(A(i_1, i_2) \otimes \ldots \otimes A(i_{n-1}, i_n) 
\otimes A(i_{n}, i_1) \Bigr)_{\Z/n\Z}.
\ee
Here the subscript $\Z/n\Z$ denotes the coinvariants of the cyclic shift.
\bd 
A potential is a linear functional  $P: {\rm HH}_0({\cal P}(Q)) \lra k$.\footnote{To 
motivate this definition observe that any element $P \in 
{\rm HH}_0({\cal P}(Q))$ gives rise to a function $F_P$ 
on the space of representations of the quiver $Q$: we assign to  a cycle $C$ in  (\ref{cycen}) 
a function $F_C$ given by the trace of the product of operators assigned to the arrows of $C$, and 
extend by linearity: 
$F_{C_1+C_2}:= F_{C_1} + F_{C_2}$.}  
\ed

Informally, the space (\ref{cycen})  is 
a vector space with a natural basis given by {\it cyclic words} of arrows of the quiver - the paths ending at their 
starting vertex. So a potential can be thought of as 
 a formal linear sum, possibly infinite, 
of cyclic words. 

\paragraph{2. An example: quivers with potentials from bipartite graphs.} 
A bipartite graph $\Gamma$ on $S$ gives rise to 
a {\it quiver with a canonical potential $(Q_\Gamma, P_\Gamma)$}. 
The quiver $Q_\Gamma$ is given by the dual graph: its vertices are the faces of the graph $\Gamma$; its edges 
are dual to the edges $E$ of $\Gamma$, oriented so that the $\bullet$-vertex of $E$ is on the right.
We call the oriented edges of the quiver {\it arrows}.

Each $\bullet$-vertex $b$ of a bipartite graph $\Gamma$ gives rise to a simple cyclic path $C_b$ on $Q_\Gamma$ 
given by the arrows going clockwise around the vertex $b$. Each $\circ$-vertex $r$ 
of  $\Gamma$ gives rise to a similar simple loop $C_r$ on $Q_\Gamma$ 
given by the arrows going counterclockwise around the vertex $r$. 
Set 
\be \la{8.30.14.1}
P_\Gamma:= \sum_{b: ~\mbox{$\bullet$-vertices of $\Gamma$}} C_b - \sum_{r: ~\mbox{$\circ$-vertices of $\Gamma$}}  
C_r \in {\rm HH}_0({\cal P}(Q_\Gamma)).
\ee
Quivers with potentials assigned to bipartite graphs were invented by physicists, see \cite{FHKVW} and references there. 
 
Below we briefly recall Kontsevich-Soibelman's construction \cite{KS1}, Section 8. 

\paragraph{3. 3d Calabi-Yau categories and cluster collections.} 
A 3d Calabi-Yau category over a characteristic zero field $k$ is a 
(weakly unital) $k$-linear triangulated $A_\infty$-category ${\cal C}$ with the following properties: 

\begin{enumerate} 
\item For any two objects $E$, $F$, the ${\rm Hom}^\bullet(E, F)$ is a $\Z$-graded finite dimensional vector space:
$$
{\rm Hom}^\bullet(E, F)= \oplus_n{\rm Hom}^n(E, F). 
$$

\item  There is a non-degenerate symmetric pairing 
\be\la{4.18.12.1}
(\bullet, \bullet): {\rm Hom}^\bullet(E, F) \otimes {\rm Hom}^\bullet(F, E) \lra k[-3].
\ee

\item  For any objects $E_1, \ldots , E_{n+1}$, there are higher composition maps 
\be\la{4.18.12.2}
m_n: \bigotimes_{1 \leq i \leq n}\Bigl({\rm Hom}^\bullet(E_i, E_{i+1})[1]\Bigr)  \lra {\rm Hom}^\bullet(E_1, E_{n+1})[2]
\ee
satisfying the axioms of the $A_\infty$-category, which we recall later on.

\item For any objects $E_1, ..., E_n$,  $n\geq 2$, we use pairing (\ref{4.18.12.1}) to produce 
 a map ($E_{n+1}=E_1$): 
\be\la{4.18.12.3}
W_n: \bigotimes_{1 \leq i \leq n}\Bigl({\rm Hom}^\bullet(E_i, E_{i+1}[1])\Bigr) \lra k.
$$
$$
W_n(a_1, ..., a_n):= (m_{n-1}(a_1, \ldots , a_{n-1}), a_n).
\ee
The last condition is that this map must be cyclically invariant. 
\end{enumerate}

In particular, the map $m_1:  {\rm Hom}^\bullet(E_1, E_{2}) \lra {\rm Hom}^\bullet(E_1, E_{2})[-1]$ is a differential: $m_1^2=0$. 
The ${\rm Ext}$-groups are its cohomology:
$$
{\rm Ext}^\bullet(E, F):= H^\bullet({\rm Hom}^\bullet(E, F), m_1).
$$

\vskip 3mm
For any object $E$, the {\it total potential of the object $E$},
 is   a formal power series 
$$
W_E(\alpha):= \sum_n\frac{W_n(\alpha, ..., \alpha)}{n}, ~~~~\alpha \in {\rm Hom}^\bullet(E,E)[1].
$$
Consider the cyclic tensor envelope  of a graded  vector space $V$:
\be \la{10.12.14.1}
{\cal C}{\rm T}(V):= ~\oplus_{n=0}^\infty (V^{\otimes n})_{\Z/n\Z}.
\ee
The total potential $W_E$ lies in the dual to the cyclic tensor envelope of ${\rm Hom}^\bullet(E, E)[1]$: 
\be \la{10.12.14.2}
{\rm Hom}\Bigl({\cal C}{\rm T}({\rm Hom}^\bullet(E, E)[1]), k\Bigr). 
\ee
 The inverse of the odd non-degenerate pairing (\ref{4.18.12.1}) is 
an odd bivector, providing an odd non-commutative Poisson bracket $\{\ast, \ast\}$ 
on the graded vector space  (\ref{10.12.14.2}), see e.g. \cite{K92}. 

\vskip 3mm
 
Let ${\cal C}$ be a non-triangulated 3d CY ${\rm A}_\infty$-category 
with finitely many objects $S_1, ..., S_n$. 
Set ${\Bbb S}:= \oplus_{i=1}^nS_i$. 
It is the same thing as a 3d CY ${\rm A}_\infty$-algebra structure 
on the graded space ${\rm Hom}^\bullet({\Bbb S}, {\Bbb S})$ equipped with idempotents 
$p_i: {\rm Hom}^\bullet({\Bbb S}, {\Bbb S}) \to {\rm Hom}^\bullet({\Bbb S}, {S}_i)$. 
Equivalently, it can be described by 
a potential $W_{\Bbb S}$ satisfying the classical master equation. Precisely, 
the total potential $W_{\Bbb S}$ of ${\Bbb S}$ is an element 
\be \la{10.12.14.2a}
W_{\Bbb S} \in {\rm Hom}\Bigl({\cal C}{\rm T}({\rm Hom}^\bullet({\Bbb S}, {\Bbb S})[1]), k\Bigr). 
\ee
It satisfies the {\it classical master equation}
\be \la{10.12.14.2asas}
 \{W_{\Bbb S}, W_{\Bbb S}\}=0.
\ee 

Vice versa,  assuming conditions 1) and 2) above, any 
element $W$ in (\ref{10.12.14.2}) 
satisfying  $\{W, W\}=0$ determines a non-triangulated 3d CY $A_\infty$-category 
${\cal C}_W$ with the objects $S_1, ..., S_n$ such that $W$ is the total potential of 
 the object ${\Bbb S}$ of ${\cal C}_W$. The maps $m_n$ are the components 
of the noncommutative Hamiltonian vector field $\{W, \ast\}$ provided by the Hamiltonian $W$ \cite{K92}.

The non-triangulated 3d CY category ${\cal C}$ gives rise to its triangulated envelope, the 
category ${\rm Tw}({\cal C})$ of twisted complexes 
of ${\cal C}$ \cite{BK}.  
Its objects are given by direct sums $X_1[n_1] \oplus \ldots \oplus X_k[n_k]$
of shifted generators $S_i$ of ${\cal C}$, 
equipped with a matrix $f = (f_{ij})_{i<j}$ of morphisms $f_{ij}: X_i \lra X_j$,  
${\rm deg}(f_{ij})=n_i-n_j+1$, so that $f:= \sum_{i<j}f_{ij}$ satisfies the 
Maurer-Cartan equation
$$
\textstyle\sum_{l=1}^k m_l(f, ..., f) =0. 
$$

\paragraph{4. A quiver with potential from a cluster collection \cite{KS1}.} A cluster collection ${\cal S}$ 
gives rise to a quiver $Q_{{\cal C}, {\cal S}}$  whose vertices $\{i\}$ are the objects $S_i$, and the arrow spaces are 
\be \la{asp}
{\rm A}(i, j):= {\rm Ext}^1(S_i, S_j).
\ee

Alternatively,  we get a quiver by considering a lattice 
 $\Lambda:= K_0({\cal C})$, with a basis given by the classes $[S_i]$, and 
a skew symmetric form given by the negative of the Euler form:
\be \la{9.30.14.1}
\begin{split}
&(\ast, \ast): K_0({\cal C}) \wedge K_0({\cal C}) \lra \Z. \\
&\textstyle (A,B):= -\sum_a(-1)^a{\rm dim}{\rm Ext}^a(A, B).\\
\end{split}
\ee
The total potential $W_{\Bbb S}$ of ${\Bbb S} = \oplus_{i=1}^nS_i$, restricted 
to tensor products of ${\rm Ext}^1(S_i, S_j)[1]$'s, 
 provides  a potential $P_{{\cal C}, {\cal S}}$  of the 
quiver $Q_{{\cal C}, {\cal S}}$:
$$
P_{{\cal C}, {\cal S}}: {\rm HH}_0({\cal P}(Q_{{\cal C}, {\cal S}})) \lra k. 
$$
This way a 3d CY $A_\infty$-category with a cluster collection $({\cal C}, {\cal S})$
gives rise to a quiver with potential $(Q_{{\cal C}, {\cal S}}, P_{{\cal C}, {\cal S}})$. 
The quiver does not have loops and cycles of length two, and the potential is {\it minimal}, 
i.e. it starts from the cubic terms. 
 
Given a quiver with potential $(Q, P)$, the group of automorphisms of the path algebra ${\cal P}(Q)$ 
of the quiver $Q$ 
preserving the idempotents corresponding to the vertices of the quiver 
acts on the potentials. Two potentials are {\it gauge equivalent} 
if belong to the same orbit of this group.

\paragraph{5. The Kontsevich-Soibelman correspondence \cite[Section 8]{KS1}.} 
It asserts that 
assigning to a 3d CY category with a cluster collection  $({\cal C}, {\cal S})$
its quiver with potential $(Q_{{\cal C}, {\cal S}}, P_{{\cal C}, {\cal S}})$ we get   
a 1-1 correspondence  
\be \la{KS}
\{\mbox{3d CY categories + cluster collections /
equivalences preserving CY str. \& cluster collections}\}
\ee
$$
\leftrightarrow \{\mbox{quivers without  length one and two cycles, 
with minimal potentials / gauge equivalence}\}. 
$$

The inverse construction proceeds in two steps.

\paragraph{i)} We assign   
 to a quiver with potential a 3d CY non-triangulated  $A_\infty$-category $\widetilde {\cal C}$ 
with a set of objects $\{S_i\}_{i\in I}$ and $m_1=0$. 
The category 
is determined by the full potential $W_{\Bbb S}$ of the object ${\Bbb S}=\oplus_{i=1}^nS_i$, defined 
as a sum of two terms 
$$
W_{\Bbb S} := W'_{\Bbb S} + W_{\rm can}\in {\rm Hom}\Bigl({\cal C}{\rm T}({\rm Ext}^\bullet({\Bbb S}, {\Bbb S})[1]), k\Bigr).
$$
The 
$W'_{\Bbb S}$ is non-zero only on the cyclic products of ${\rm Ext}^1({\Bbb S}, {\Bbb S})[1]$. 
It is   given by the  quiver potential. 

The 
$W_{\rm can}$  is defined by postulating that its only non-zero components are the maps 
\be \la{10.12.14.3}
{\rm Ext}^a(S_i, S_j)[1] \otimes {\rm Ext}^{3-a}(S_j, S_i)[1] \otimes 
{\rm Ext}^0(S_i, S_i)[1] \lra k, ~~~~a=0, 1, 2. 
\ee
These maps are induced by the pairing  
${\rm Ext}^a(S_i, S_j) \otimes {\rm Ext}^{3-a}(S_j, S_i) \lra k$ for $a=0, 1$ and by its negative for $a=2$. 
The term $W_{\rm can}$ can be defined by a formula. Namely, denote by $x_{ij}^s$ a basis of coordinates on ${\rm Ext}^1(S_i, S_j)[1]$, 
by $\xi_{ji}^s$ the dual basis for  ${\rm Ext}^2(S_j, S_i)[1]$, and by $\alpha_i$ and $a_i$ the natural coordinates on 
${\rm Ext}^0(S_i, S_i)[1]=k[1]$ and ${\rm Ext}^3(S_i, S_i)[1] = k[-2]$. Then 
$$
\textstyle W_{\rm can}(\alpha, x, \xi, a):= 
\sum_{i=1}^n\alpha_i^2a_i + \sum_{i, j=1}^n(\alpha_ix_{ij}^s\xi_{ji}^s - \alpha_i\xi_{ij}^sx_{ji}^s).  
$$

The potential $W_{\rm can}$ provides the graded vector space 
${\rm Ext}^\bullet(\S, \S)$ with a structure of a graded associative 
3d CY algebra with a collection of idempotents $\{p_i\}_{i\in I}$ and a unit $e$:
$$
\textstyle p_i:= 1\in {\rm Ext}^0(S_i, S_i), ~~~~e:= \sum_{i\in I} p_i.
$$
 The multiplication vanishes on the graded component 
$$
{\rm Ext}^1(\S,\S) \otimes {\rm Ext}^1(\S,\S) \to {\rm Ext}^2(\S,\S).
$$ It is given by a non-degenerate 
bilinear form on $$
{\rm Ext}^1(\S,\S) \otimes {\rm Ext}^2(\S,\S) \to {\rm Ext}^3(\S,\S).
$$
This implies that $\{W_{\rm can}, W_{\rm can}\} =0$, which  one can easily check directly. 
Clearly, $\{W'_{\Bbb S}, W'_{\Bbb S}\} =0$.  
One proves that 
$\{W_{\rm can}, W'_{\Bbb S}\} =0$. 
 Thus $\{W_{\Bbb S}, W_{\Bbb S}\} =0.$ 
Therefore we get an $A_\infty$ CY category 
$\widetilde {\cal C}_W$ with the objects $S_1, ..., S_n$, which by construction form a cluster collection.

\paragraph{ii)}  The desired 3d CY category ${\cal C}$ is the 
category ${\rm Tw}(\widetilde {\cal C})$ of twisted complexes 
of $\widetilde {\cal C}$. 

\vskip 3mm
The category ${\rm Tw}(\widetilde {\cal C})$ can be described 
as the DG category of finite dimensional DG modules over the Ginzburg algebra 
related to the quiver with potential $(Q_{{\cal C}, {\cal S}}, P_{{\cal C}, {\cal S}})$.  
Although we are not using this, for convenience of the reader 
we recall the definitions. 

\paragraph{6. The Ginzburg algebra of a quiver with potential \cite{G}.} 
The {\it Ginzburg algebra} \cite{G} associated to a quiver with potential $(Q, P)$  is a DG-algebra $A^\bullet_{Q, P}$ 
concentrated in non-positive degrees, 
with the following generators in the degrees $0, -1, -2$: 
\begin{itemize}
\item Projectors $\{e_i\}$,  $e_i^2=e_i$, $\sum_ie_i=1$, corresponding to vertices $\{i\}$ of the quiver $Q$. 

\item Degree $0$ generators $\{a_E\}$, and degree $-1$ generators $\{a^*_E\}$, assigned to arrows $\{E\}$ of $Q$. 

\item Degree $-2$ generators $\{t_i\}$ assigned to vertices $\{i\}$ of the quiver $Q$. 

\item The differential $d$ acts on the generators by 
\be \la{18.4.12.100}
da_E=0, \quad da^*_E = \frac{\partial P}{\partial a_E}, \quad dt_i = \sum_{E~|~ t(E) =i}a^*_Ea_E - \sum_{E~|~ h(E) =i}a_Ea_E^*.
\ee
\end{itemize}
Here $\frac{\partial P}{\partial a_E}$ is the noncommutative derivative. 
These definitions can be written using only arrow spaces. Indeed, choose a basis $\{a_E\}$ in each arrow space;  
let $\{a_E^*\}$ be the dual basis in the dual space. Then to write the differential we employ the identity elements 
$$
\textstyle \sum_{E}a_E\otimes a^*_E \in {\rm A}(i, j)\otimes {\rm A}(i, j)^*, \qquad 
\sum_{E}a^*_E\otimes a_E \in {\rm A}(i, j)^*\otimes {\rm A}(i, j).
$$
The {\it Jacobian algebra $H^0(A^\bullet_{Q, P})$} 
is generated by the generators $e_v, a_E$ satisfying the relations $\frac{\partial P}{\partial a_E}=0$ for each arrow $E$.

Denote by ${\cal F}_{Q, P}$ the DG-category of 
DG-modules over the Ginzburg algebra $A^\bullet_{Q, P}$,  whose homology are finite dimensional,  
nilpotent modules over $H^0(A^\bullet_{Q, P})$. 
The subcategory of degree zero $A^\bullet_{Q, P}$-modules 
is an abelian heart of this category. It is the category $H^0(A^\bullet_{Q, P})$-modules. 

\bl
The DG category ${\cal F}_{Q, P}$ is canonically equivalent to 
the category of twisted complexes over the CY $A_\infty$-category $\widetilde {\cal C}_{Q, P}$ defined above. 
The cluster collection $\{S_i\}$ in $\widetilde {\cal C}_{Q, P}$ is identified with the 
one dimensional $A^\bullet_{Q, P}$-modules corresponding to the vertices $\{i\}$ of $Q$. 
\el



There is another 3d Calabi-Yau category ${\cal P}_{Q, P}$, the category of perfect complexes over the Ginzburg algebra 
 $A^\bullet_{Q, P}$. It is Koszule dual to the category ${\cal F}_{Q, P}$.  
Namely, the category ${\cal F}_{Q, P}$ is equivalent to the category of DG-functors from the category ${\cal P}_{Q, P}$ 
to the category ${\rm Vect}^\bullet$ of bounded complexes of vector spaces with finite dimensional cohomology, and vice versa:
$$
{\cal F}_{Q, P} \sim {\rm Funct}_{dg}({\cal P}_{Q, P}, {\rm Vect}^\bullet), \qquad 
{\cal P}_{Q, P} \sim {\rm Funct}_{dg}({\cal F}_{Q, P}, {\rm Vect}^\bullet). 
$$

\paragraph{7. Mutations of cluster collections \cite{KS1}.} A spherical object $S$ of a CY category 
 gives rise to the 
Seidel-Thomas  \cite{ST} reflection  functor ${\rm R}_{S}$, 
acting by an autoequivelence 
of the category:
$$
R_S(X):= {\rm Cone}\Bigl({\rm Ext}^\bullet(S, X) \otimes S \lra X\Bigr).
$$  

Given elements $0, i\in I$,  we  write $i<0$  if ${\rm Ext}^1(S_i, S_0)$ is non-zero,
 and $i>0$ otherwise.  
\bd
A mutation of a cluster collection ${\cal S} = \{S_i\}_{i\in I}$ in a 3d CY  category ${\cal C}$ at an object $S_0$ 
is a new spherical collection ${\cal S}' =\{S'_i\}$ in the same category ${\cal C}$ given by 
$$
S'_i=S_i, ~~i<0, ~~~~S'_0=S_0[-1], ~~~~S'_i=R_{S_0}(S_i), ~~i>0,
$$
\ed
Since the shift and reflection functors transform spherical objects to spherical ones, 
$\{S'_i\}_{i\in I}$ is a collection of spherical objects. Since the objects $\{S_i\}$ 
generate the triangulated category ${\cal C}$, 
the objects $\{S'_i\}$  are generators. 
However $\{S'_i\}_{i\in I}$ is not necessarily a cluster collection. 

 Mutations of cluster collections, being projected to 
$K_0({\cal C}) $, recover mutations of bases (\ref{12.12.04.2a}).

\vskip 3mm

Mutations of cluster collections  
translate into {\it mutations of quivers with generic potentials}  introduced earlier by
 Derksen, Weyman and Zelevinsky \cite{DWZ}. 
  Keller and Yang \cite{KY}
 promoted mutations of quivers with generic potentials $(Q,P) \to (Q',P')$ 
to {\it mutation functors} -- explicit equivalences ${\cal F}_{Q, P} \lra {\cal F}_{Q', P'}$ 
between the corresponding 3d CY categories. 

In the Kontsevich-Soibelman picture 
a mutation does not affect the 3d CY category, altering only the cluster collection. 
It is build on the cluster collection mutation  ansatz, which a posteriori implies the Derksen-Weyman-Zelevinsky 
mutation rule  for quivers with potentials. The ansatz itself is a categorification 
of the geometric form of quiver mutation (\ref{12.12.04.2a}), and 
related to tiltings of $t$-structures \cite{Br1}. 

\vskip 3mm

In symplectic topology 
we are given a 3d CY category - the Fukaya category, while construction of a cluster collection,  
realised by special Lagrangian spheres, requires some choices. 

To produce a 3d CY category combinatorially, we need as an input 
a quiver with potential. However there is no preferred quiver.
Furthermore, there is  an issue: 

\begin{itemize}

\item {\it A mutation of a cluster collection may deliver a non-cluster  collection.}

\item 

 {\it Mutations of quivers with potentials are defined only for 
quivers with \underline{generic} potentials}.
\end{itemize}

Given a
quiver with potential, it is hard to determine whether after a finite number of mutations 
we will get a quiver with potential which does not allow mutations in some directions.

\paragraph{8. Cluster varieties and quivers with \underline{canonical} potentials.} 
In general, there is no canonical potential assigned to an arbitrary  quiver. 
The combinatorial part of our proposal in Section \ref{anover}  boils down to the following: 

\begin{itemize}
\item Every cluster variety ${\cal V}$ which appears 
in representation theory, geometry  and physics admits a ``tame'' collection of quivers, 
equipped with \underline{canonical} potentials. 

The corresponding 3d CY category is the 
combinatorial category ${\cal C}_{\cal V}$. 

\item 
Any two quivers with potentials from this collection can be connected by mutations 
inside of the collection; elements of the group $\Gamma_{\cal V}$ are realised by compositions of mutations. 

\item The equivalence (\ref{equiva}) identifies the cluster collections in ${\cal C}_{\cal V}$ 
 with the cluster collections provided by 
special Lagrangian spheres in the CY threefold ${\cal Y}_b$. 
\end{itemize}

Recall  that collections of bipartite surface graphs related by two by two moves 
 give rise to a cluster Poisson variety \cite{GK}. 
The quiver of a bipartite surface graph 
comes with a canonical potential (\ref{8.30.14.1}). 
We suggest that  this is a universal source of cluster varieties in the examples related to simple groups of type $A$.

\subsection{The symmetry group of a ${\rm CY}_3$ category with a cluster collection} \la{SSec2}

Let ${\cal S}=\{S_i\}$ be a cluster collection in ${\rm CY}_3$ category ${\cal C}$. 
Consider the quotient group
\be  \la{autcc}
{\rm Auteq}({\cal C}; {\cal S}):= \frac{\mbox{\rm Autoequivalences of the category ${\cal C}$}}
{\mbox{\rm Autoequivalences preserving objects of the cluster collection ${\cal S}$}}.
\ee
Our first goal is to define a group of symmetries of a pair $({\cal C}; {\cal S})$:
$$
\Gamma_{{\cal C}, {\cal S}}\subset {\rm Auteq}({\cal C}; {\cal S}).
$$

\paragraph{1. The braid group of a cluster collection.} 
\bd 
Given a cluster collection ${\cal S} = \{S_i\}_{i\in I}$ in ${\cal C}$, the 
reflection functors ${\rm R}_{S_i}$ and their inverses 
generate a subgroup, called the {\it braid group of the cluster collection ${\cal S}$}:    
$$
{\rm Br}_{{\cal C}, {\cal S}} \subset {\rm Auteq}({\cal C}; {\cal S}).
$$ 
\ed 
By Propositions 2.12 -  2.13  in \cite{ST}, the 
functors ${\rm R}_{S_i}$ satisfy 
the ``braid relations'':
$$
{\rm R}_{S_a}{\rm R}_{S_b} = {\rm R}_{S_b} {\rm R}_{S_a} ~~\mbox{if ${\rm Ext}^1(S_a, S_b)=0$}.
$$
$$
{\rm R}_{S_a}{\rm R}_{S_b}{\rm R}_{S_a} = {\rm R}_{S_b} {\rm R}_{S_a}{\rm R}_{S_b}
 ~~\mbox{if ${\rm dim}{\rm Ext}^1(S_a, S_b)=1$}.
$$
\paragraph{2. Categorified cluster modular groupoid.} 
\bd \la{10.1.14.1} The  categorified cluster modular groupoid ${\cal C}{\rm Mod}$ is the following  groupoid:

\begin{itemize}

\item Its objects 
are pairs $({\cal C}, {\cal S})$, where ${\cal S}$ is a cluster collection
in a 3d CY category ${\cal C}$. 

\item Morphisms $\psi: ({\cal C}, {\cal S}) \lra ({\cal C}', {\cal S}')$ are given by certain pairs 
$(\alpha, \beta)$, where $\alpha: {\cal C} \lra {\cal C}'$ is an equivalence of categories, and 
$\beta: {\cal S} \to {\cal S}'$ is a bijection. Precisely:

The morphisms are generated by the following elementary ones: 

i)  mutations $\mu_{S_k}: ({\cal C}, {\cal S}) \lra ({\cal C}, \mu_{S_k}({\cal S}))$, 
where ${\cal C} \to {\cal C}$ is the  identity functor, and 
$\mu_{S_k}({\cal S})$ is the mutation of the cluster collection ${\cal S}$ 
in the direction $S_k$;

ii) equivalences  $\varphi: ({\cal C}_1, {\cal S}_1) \stackrel{\sim}{\lra} ({\cal C}_2, \varphi({\cal S}_1))$, 
where $\varphi: {\cal C}_1 \lra {\cal C}_2$ is an any equivalence.  

Compositions of morphisms are defined in the obvious way. 
\item 
Relations: two compositions of generating morphisms 
$\psi_1, \psi_2: ({\cal C}_1, {\cal S}_1) \lra ({\cal C}_2, {\cal S}_2)$ are equal 
if the underlying equivalences coincide, and 
 $\beta_1(S_i) = \beta_2(S_i)$, $\forall i\in I$. 

\end{itemize}

\ed

We stress that mutations 
do not affect the category, changing the cluster collections only.

\bl \la{sqmrf} In the groupoid ${\cal C}{\rm Mod}$ 
the square of a mutation is the reflection functor:  
\be \la{9.30.14.10}
\mu_{S[-1]}\circ \mu_{S} = {\rm R}_S~~~~\forall S \in {\cal S}.
\ee
\el

\begin{proof}
The composition $\mu_{S[-1]}\circ \mu_{S}$ is given by a pair $(\alpha, \beta)$ where $\alpha$ 
 is the identity functor on ${\cal C}$, and $\beta$ 
transforms the cluster collection ${\cal S}$ 
to the one   ${\rm R}_S({\cal S})$. 
The  reflection functor ${\rm R}_S$  
is an autoequivalence of the category ${\cal C}$, which has the same effect 
on the cluster collection. This just means that 
we have (\ref{9.30.14.10}). 
\end{proof}

Let ${\rm Aut}_{{\cal C}{\rm Mod}}({\cal C}, {\cal S})$ be the automorphism group 
of the object  $({\cal C}, {\cal S})$ in the groupoid ${\cal C}{\rm Mod}$. 

\bl There is a canonical injective homomorphism 
\be \la{10.2.14.1}
{\rm Aut}_{{\cal C}{\rm Mod}}({\cal C}, {\cal S}) \lra {\rm Auteq}({\cal C}; {\cal S}).
\ee
\el

\begin{proof} Any element of 
${\rm Aut}_{{\cal C}{\rm Mod}}({\cal C}, {\cal S})$ is given by a sequence of 
generating morphisms
$$
({\cal C}, {\cal S}) \stackrel{(\alpha_0, \beta_0)}{\lra} ({\cal C}_1, {\cal S}_1)
\stackrel{(\alpha_1, \beta_1)}{\lra}  ({\cal C}_2, {\cal S}_2) \stackrel{(\alpha_2, \beta_2)}{\lra}  \cdots \stackrel{(\alpha_{n-1}, \beta_{n-1})}{\lra}  ({\cal C}_n, {\cal S}_n) 
\stackrel{(\alpha_n, \beta_n)}{\lra} 
({\cal C}, {\cal S}).
$$
The projection of the composition $\alpha_n\circ \ldots \circ \alpha_1\circ \alpha_0 \in {\rm Auteq}({\cal C})$ 
 to the group 
${\rm Auteq}({\cal C}, {\cal S})$ provides a group homomorphism 
(\ref{10.2.14.1}). It  is injective since, by definition, two 
sequences resulting in the same element of the group ${\rm Auteq}({\cal C}; {\cal S})$ determine the same 
morphism of the groupoid ${\cal C}{\rm Mod}$. 
\end{proof}

\paragraph{3. Comparing groupoids ${\cal C}{\rm Mod}$ and ${\rm Mod}$.}
The categorified cluster modular groupoid ${\cal C}{\rm Mod}$ is similar to the cluster modular groupoid 
${\rm Mod}$ 
defined in \cite{FG2}. 

The objects of the groupoid ${\rm Mod}$ are  quivers 
$(\Lambda, \{e_i\}_{i\in I}, (\ast, \ast))$. Its morphisms are generated by the elementary ones of two types: 
quiver mutations (\ref{12.12.04.2a}) and isomorphisms of quivers. 
Similarly to Definition \ref{10.1.14.1}, 
quiver mutations do not affect the lattice with the form $(\Lambda, (\ast, \ast))$, changing bases only. 
Two compositions of elementary morphisms are equal if they induce the same 
cluster transformation of the corresponding quantum cluster variety \cite{FG2}. 

\bt There  is a canonical functor, see (\ref{9.30.14.1}): 
\be \la{10.1.14.5}
{\cal C}{\rm Mod} \lra {\rm Mod}, ~~~~({\cal C}, {\cal S}) \lra (K_0({\cal C}), \{[S_i]\}, [\ast, \ast]).
\ee
\et

\begin{proof}
The cluster collection mutation $\mu_{S_k}$ induces in $K_0({\cal C})$ 
the halfreflection (\ref{12.12.04.2a}). Therefore the assignment (\ref{10.1.14.5}) 
is defined on the generating morphisms. The fact that it sends the 
relations between the generating morphisms ${\cal C}{\rm Mod}$ 
to the relations in the cluster modular groupoid ${\rm Mod}$ is  implied by   \cite[Theorem 5.2]{K11}.
\end{proof}

The difference between the two groupoids is that in the cluster modular groupoid the square of a mutation is the identity: 
$\mu_{-e_k}\mu_{e_k}={\rm Id}$, while in  ${\cal C}{\rm Mod}$ it is the reflection functor: 

\bl The group ${\rm Aut}_{{\cal C}{\rm Mod}}({\cal C}, {\cal S})$ is an extension of the cluster modular group 
by the cluster braid group:
\be \la{9.30.14.2}
1 \lra {\rm Br}_{{\cal C}, {\cal S}} \lra {\rm Aut}_{{\cal C}{\rm Mod}}({\cal C}, {\cal S}) \lra 
{\rm Aut}_{{\rm Mod}}(K_0({\cal C}), [{\cal S}], (\ast, \ast))\lra 1.
\ee
\el

\begin{proof} Thanks to relation (\ref{9.30.14.10}), 
the cluster braid group ${\rm Br}_{{\cal C}, {\cal S}}$ is a subgroup 
of ${\rm Aut}_{{\cal C}{\rm Mod}}({\cal C}, {\cal S})$. 
\end{proof}

\subsection{${\rm A}_{m}$-webs and  extended mapping class groups} \la{SSec2.3}

\paragraph{1. Categorified two by two moves.} 
A web ${\cal W}$ on $S$ gives rise to a quiver with potential $(Q_{\cal W}, P_{\cal W})$, and therefore to a 
 3d CY category ${\cal C}_{\cal W}:= {\cal C}_{Q_{\cal W}, P_{\cal W}}$ with a cluster collection 
${\cal S}_{\cal W}$, whose objects are  parametrised by the faces of the web ${\cal W}$. 

A two by two move 
$\mu_F: {\cal W}\lra {\cal W}'$ of 
 webs,  
centered at a face $F$ of  ${\cal W}$, 
gives rise to a quiver with potential $(Q_{{\cal W}'}, P_{{\cal W}'})$, and  hence to a 
3d CY category with a cluster collection $({\cal C}_{\cal W'}, {\cal S}_{\cal W'})$. 

\bp \la{2b2c} The quiver with potential $(Q_{{\cal W}'}, P_{{\cal W}'})$ 
describes the mutated cluster collection  
$\mu_F({\cal S}_{\cal W})$ in the category ${\cal C}_{\cal W}$. 
In particular 
the categories  ${\cal C}_{\cal W'}$  and ${\cal C}_{\cal W}$ are equivalent.
\ep

Proposition \ref{2b2c} is proved in Section \ref{sec5.4}

\vskip 3mm
Therefore we arrived at the categorified two by two move:
$$
\mu_F: ({\cal C}_{\cal W}, {\cal S}_{\cal W}) \lra ({\cal C}_{\cal W}, \mu_F({\cal S}_{\cal W})) \sim 
({\cal C}_{\cal W'}, {\cal S}_{\cal W'}) .
$$

Let us explain now how we construct symmetries of the category with cluster collection related to 
the group $PGL_m$ and a decorated surface $S$.

\paragraph{2. Groupoid of regular ideal triangulations.} 
An ideal triangulation of a decorated surface $S$ has triangles of two kinds: 
{\it regular triangles}, with three distinct sides, and {\it folded triangles}, where 
two sides are glued together. Some or all vertices of a regular 
ideal triangle can be glued together. An example is given by an ideal triangulation of a punctured torus.

Below we consider only ideal triangulations 
which do not have folded triangles - we call them {\it regular ideal triangulations}. 
We define the {\it groupoid ${\rm Tr}_S$ of  regular ideal triangulations of $S$}: 

\bd Let $S$ be a decorated surface. 

\begin{itemize}

\item The objects of the groupoid ${\rm Tr}_S$ 
are regular ideal triangulations of $S$. 

\item 
The morphisms are generated by 
isomorphisms of decorated surfaces and flips:

i) Any isomorphism 
$i:S \to S'$ gives rise to morphisms $i: (S,T) \to (S', i(T))$. 

ii) Any internal edge $E$ of a regular ideal triangulation $T$ provides a morphism $\mu_E: T \to T_E$ 
where $T_E$ is the unique modulo isotopy ideal triangulation obtained by flipping the edge $E$. 

\item The relations between the morphisms are generated by the following: 

a) Isomorphisms commute with the flips: $i \circ \mu_E = \mu_{i(E)}\circ i$. 

b) The square  relations: flips at disjoint edges commute. 

c) The pentagon relations. 
\end{itemize}
\ed

The  mapping class group $\Gamma_S$ is isomorphic to the automorphism  group  ${\rm Aut}_{{\rm Tr}_S}(T)$ 
of the object  of 
 groupoid ${\cal T}_S$ given a regular ideal triangulation $T$. 
Indeed, for any $g \in \Gamma_{S}$, there is a sequence of flips of regular ideal triangulations 
$T\to T_1 \to \ldots \to g(T)$ 
connecting $T$ with $g(T)$. Combining with 
the isomorphism $g^{-1}: g(T) \to T$, we get an automorphism $a_{g,T} \in {\rm Aut}_{{\cal T}_S}(T)$, 
 given by the composition
$$
T\to T_1 \to \ldots \to g(T) \stackrel{g^{-1}}{\lra} T.
$$
Different sequences of flips connecting $T$ with $g(T)$ can be homotoped to each other 
by using the square and pentagon relations. So the element $a_{g,T}$ is well defined. 
The map $g \lms a_{g,T}$ is a group homomorphism. It is known to be an isomorphism.

\paragraph{3. The extended mapping class group $\Gamma_{G, S}$.} 
We assigned  to a regular ideal triangulation $T$ of $S$ an ideal  
${\rm A}_{m-1}$-bipartite graph, see Figure \ref{gra10i}, and hence an ideal ${\rm A}_{m-1}$-web ${\cal W}_T$. 
So we get a quiver with potential, and thus  a 3d CY category 
with cluster collection $({\cal C}_{G, S}^T, {\cal S}_{G, S}^T)$, $G=PGL_m$. 

\bd Let $S$  be a decorated surface  and $G=PGL_m$. 
Then   
$$
\widetilde \Gamma_{G, S}:= {\rm Aut}_{{\cal C}{\rm Mod}}({\cal C}_{G, S}^T, {\cal S}_{G, S}^T)
$$ 
is  the automorphism group of the object $({\cal C}_{G, S}^T, {\cal S}_{G, S}^T)$ of the  groupoid ${\cal C}{\rm Mod}$.
\ed

\bl
The isomorphism class 
of the group $\widetilde \Gamma_{G, S}$ does not depend on the choice of $T$. 
 \el

\begin{proof} Any two regular ideal triangulations of $S$ are related by a sequence of flips, 
A flip $T\to T'$ gives rise 
to a sequence of  $(m-1)^2$ two by two moves, and 
shrink / expand moves, transforming 
the ${\rm A}_{m-1}$-web ${\cal W}_{T}$  to  ${\cal W}_{T'}$, see  Section \ref{2by2moves}.   
Proposition \ref{2b2c} 
categorifies the two by two moves  to 
mutations of cluster collections. 
Their composition is a morphism 
 in the groupoid ${\cal C}{\rm Mod}$:
\be \la{9.6.14.1}
{f}^m_T: ({\cal C}_{G, S}^T, {\cal S}_{G, S}^T)\lra ({\cal C}_{G, S}^{T'},{\cal S}_{G, S}^{T'}).
\ee 
So all objects $({\cal C}_{G, S}^T, {\cal S}_{G, S}^T)$ belong to the same connected 
component of the groupoid ${\cal C}{\rm Mod}$. 
\end{proof} 

The 
cluster braid subgroup ${\rm Br}_{G, S} \subset {\rm Aut}_{{\cal C}{\rm Mod}}({\cal C}_{G, S}^T, {\cal S}_{G, S}^T)$ provides an extension 
\be \la{10.1.14.10}
1 \lra {\rm Br}_{G, S} \lra \widetilde \Gamma_{G, S} \lra \widetilde \Gamma_{G, S}/{\rm Br}_{G, S}\lra  1.
\ee

Denote by $\overline {{\cal C}{\rm Mod}}$ the groupoid obtained from ${{\cal C}{\rm Mod}}$ 
by imposing conditions $\mu_{S_k[-1]}\circ\mu_{S_k}={\rm Id}$.  
Lemma \ref{sqmrf} implies that ${\rm Aut}_{\overline {{\cal C}{\rm Mod}}}({\cal C}, {\cal S}) = 
\widetilde \Gamma_{G, S}/{\rm Br}_{G, S}$.

 \bt \la{10.1.14.6} 
There is a morphism of groupoids 
$\tau: {\rm Tr}_S \lra \overline {{\cal C}{\rm Mod}}$ providing a group map 
\be \la{10.1.14.11}
\tau: \Gamma_S \lra \widetilde \Gamma_{G, S}/{\rm Br}_{G, S}. 
\ee
\et

\begin{proof} 
We define a functor $\tau$ on objects by setting 
$
\tau: T \lms ({\cal C}_{G, S}^T, {\cal S}_{G, S}^T).
$ 
Its action on the generating morphisms given by by the 
flips is given by $\tau(T\to T'):= f^m_T$. It remains to show that $\tau$ transforms 
the compositions of flips related 
to the square and pentagon relations to the identity. 
This is obvious  for the square relations, but it is  a rather non-trivial task for the pentagon relations. 
We prove the pentagon relation in Section \ref{slmpentagon}. First, we prove the basic pentagon relation 
in the cluster set-up. 
Then we use a geometric result of \cite{DGG}, translated 
from the language of octahedral moves into the language of web mutations to deduce 
the general pentagon relation to the basic pentagon relation.
\end{proof}

\bd \la{Def5.15}
The group $\Gamma_{G, S}$ is the pull back via map (\ref{10.1.14.11}) of  extension (\ref{10.1.14.10}):
$$
\begin{array}{ccccccccc}
1 &\lra &{\rm Br}_{G, S} &\lra &\widetilde \Gamma_{G, S} &\lra &\widetilde \Gamma_{G, S}/{\rm Br}_{G, S}&\lra & 1\\
&&\uparrow = &&\uparrow &&\uparrow \pi &&\\
1 &\lra &{\rm Br}_{G, S}&\lra &\Gamma_{G, S} &\lra &\Gamma_S &\lra & 1
\end{array}
$$
\ed

We conclude that the two by two moves of ideal ${\rm A}_{m-1}$-webs on a decorated surface $S$ 
give rise to an action of 
 the group $\Gamma_{G, S}$ by symmetries of the category 
${\cal C}_{G, S}^T$, i.e. to a homomorphism
\be \la{9.30.14.3}
\Gamma_{G, S} \lra {\rm Aut}_{{\cal C}{\rm Mod}}({\cal C}_{G, S}^T, {\cal S}_{G, S}^T).
\ee

\bt
Let $G=SL_m$. Then the homomorphism (\ref{9.30.14.3}) is injective. 
\et

This follows from a general result of Keller. 

\begin{figure}[ht]
\centerline{\epsfbox{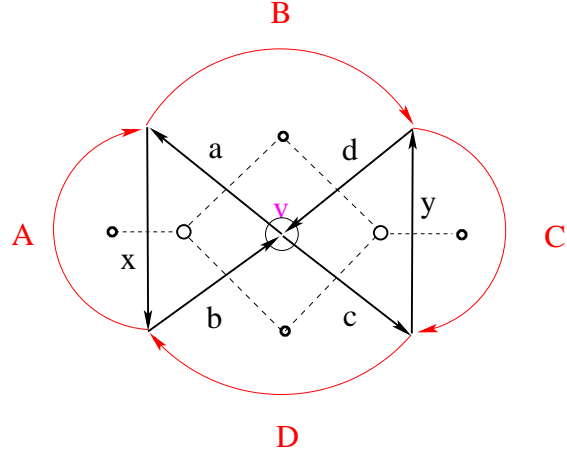}}
\caption{The potential $W= xab + Ax+ycd+yC+Bad+Dcb + W'$.}
\label{gcp1}
\end{figure}

\subsection{The potentials and  two by two moves} \la{slmpentagon} \la{sec5.4}

Denote by $h(a)$ the head of an arrow $a$, and by $t(a)$ the tale. We write $ab$ for the composition of two arrows $b$ and $a$, where $b$ is followed by $a$, 
so that $h(b) = t(a)$. Recall that given a potential $W$ on a quiver $Q$, a mutation of the pair $(Q, W)$ at a vertex $v$ is performed as follows \cite{DWZ}. 

\begin{itemize}
\item For each arrow $p$ incident to $v$: reverse its order to the opposite one, getting an arrow  $p^*$.

\item For each loop $\alpha$ passing through the vertex $v$: for each consecutive pair of arrows $p,q$ such  that $h(p)=t(q)=v$, so that $\alpha = Aqp$,
  alter the loop  $\alpha$ to the loop $A[qp]$. 

 \item For each pair of arrows $p,q$ 
such that $h(p)=t(q)=v$: add new arrow $[qp]$ to the graph, and add a new loop $[qp]p^*q^*$ to the potential. 
 
\end{itemize}

We can alter a quiver with potential, getting an {\it equivalent quiver with potential}, in two ways: 

1. By applying to it an automorphism of the path algebra of the quiver.

2. Presenting a quiver with potential as a direct sum of the two quivers with potentials 
$$
(Q, W) = (Q_1, W_1) \oplus (Q_2, W_2),
$$ where $W_1$ is a sum of two-edge loops: $W_1 = \sum_i a_ib_i$. The "direct sum" means  
 that the $a_i, b_i$ do not enter to the $W_2$. Then 
 the quiver with potential $(Q, W)$ is said to be equivalent to 
$(Q_2, W_2)$.

Consider a two by two move on a bipartite graph $\Gamma \to \Gamma'$ centered at a vertex $v$.  

The  potential of  the original quiver $Q_\Gamma$ can be written  as follows, see Figure \ref{gcp1}:
$$
W= xab + Ax+ycd+yC+Bad+Dcb + W'.
$$
Here the $W'$ is the part of the potential  free from $a,b,c,d,x,y$. 
The $A$ is a linear combination of the paths each of  which start at $h(x)$ and end at $t(x)$, 
so the composition $Ax$ makes sense. Similarly $B, C, D$ are linear combinations of paths such that 
the compositions $Bad$, $Cy$, $Dcb$ make sense.

\begin{figure}[ht]
\centerline{\epsfbox{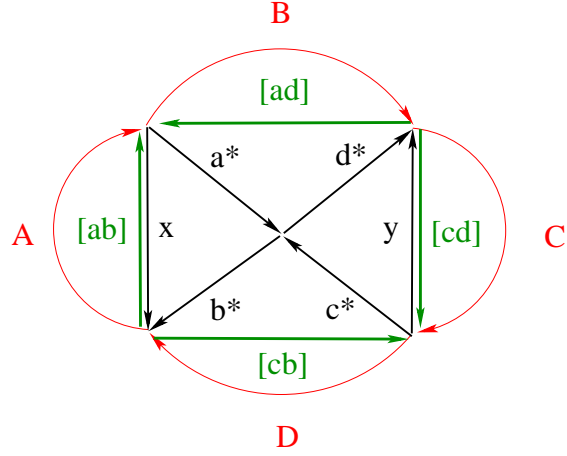}}
\caption{The  mutated potential 
$\widetilde W$.}
\label{gcp2}
\end{figure}
Let us mutate the quiver with potential $W$ at the central vertex $v$ of the quiver $Q_\Gamma$, shown by a circle on Figure \ref{gcp1}. 
We get a new  quiver with a potential $\widetilde W$:
$$
\widetilde W= x[ab] +[cd]y+ Ax+ B[ad]+yC+ D[cb] +[cd]d^*c^* + [ad]d^*a^*+[cb]b^*c^*+ [ab]b^*a^*+ W'.
$$

Consider an automorphism of the path algebra which acts on the generators identically except on the following ones: 
$$
[ab] \lms [ab] - A, ~~~~x \lms x-b^*a^*, ~~~~[cd] \lms [cd] - C, ~~~~y \lms y - d^*c^*.
$$
It sends  the  potential $\widetilde  W$ to the following one: 
\be
\begin{split}
&(x-b^*a^*)([ab]-A) +([cd]-C)(y-d^*c^*)+ A(x-b^*a^*)+ B[ad]+(y-d^*c^*)C+ D[cb] +\\
&([cd]-C)d^*c^* + [ad]d^*a^*+[cb]b^*c^*+ ([ab]-A)b^*a^*+ W'   = \\
&x[ab] +[cd]y+ B[ad]+ D[cb] -Cd^*c^* + [ad]d^*a^*+[cb]b^*c^*-Ab^*a^*+ W'. \\
\end{split}
\ee
This potential is a direct sum of the potential $x[ab] +[cd]y$ and the potential 
$$
\widetilde W^*= B[ad]+ D[cb] -Cd^*c^* -Ab^*a^*+ [ad]d^*a^*+[cb]b^*c^*+ W'.
$$
\begin{figure}[ht]
\centerline{\epsfbox{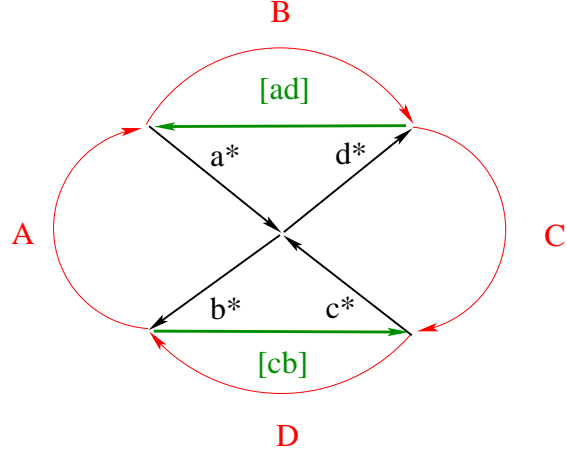}}
\caption{The  equivalent quiver with potential 
$W^*$.}
\label{gcp3}
\end{figure}
Let us consider the automorphism of the path algebra which acts on the generators identically except the following ones:
$$
b^* \lms -b^*, ~~~~d^* \lms -d^*, 
$$
It transforms the potential $\widetilde W^*$ to the following one:
$$
W^*= B[ad]+ D[cb] +Cd^*c^* +Ab^*a^*- [ad]d^*a^*-[cb]b^*c^*+ W'.
$$
\begin{figure}[ht]
\centerline{\epsfbox{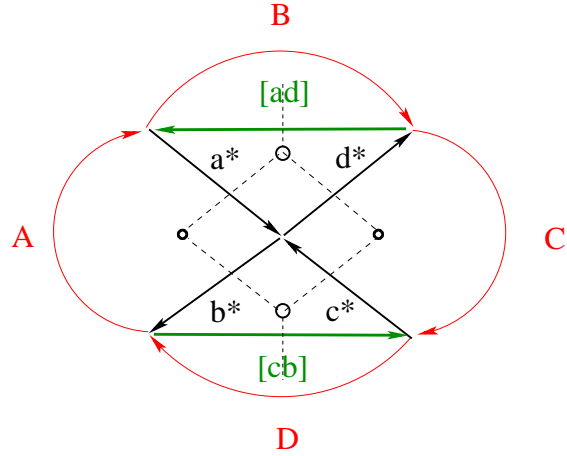}}
\caption{The  mutated potential 
$W^*$ is the potential for the quiver $Q_{\Gamma'}$.}
\label{gcp4}
\end{figure}
But this is nothing else but the potential for the quiver $Q_{\Gamma'}$ for the bipartite graph  $\Gamma'$ obtained from $\Gamma$ by the two my two move 
at the vertex $v$, see Figure \ref{gcp4}.

\section{3d Calabi-Yau categories from ideal webs and open Calabi-Yau threefolds} \la{SSec2.4}
\subsection{Geometry of the Hitchin fibration} \la{SSEEC5.3}

 \paragraph{1. The universal Hitchin base.}  Let $C$ be  a  smooth projective genus $g$ curve with a set of punctures $P:= \{p_1, ..., p_n\}$.  The moduli space of pairs $(C, P)$ is denoted by ${\cal M}_{g,n}$.   
Let $G$ be a simply laced split reductive group with connected center, e.g. $PGL_m$ or $GL_m$. 
Denote by ${\rm Lie }({\rm H})$ the Lie algebra of the Cartan group of $G$, and by $W$ the Weyl group of $G$. Let $\Omega_C(P)$ be the sheaf of meromorphic differentials on $C$ with logarithmic poles at the punctures. The 
 Hitchin base ${\rm B}_{  G, C, P}$ for a group $G$ and a punctured curve $(C, P)$ is given by
$$
{\rm B}_{G, C, P}:= H^0(C, (\Omega_{C}(P)\otimes {\rm Lie}({\rm H}))/W). 
$$
When the pair $(C, P)$ varies over the moduli space ${\cal M}_{g, n}$, 
we get the {\it universal Hitchin base ${\rm B}_{  G, g,n}$}. 
It is fibered over the moduli space the ${\cal M}_{g, n}$ with the fiber ${\rm B}_{  G, C, P}$ over 
a point $(C, P)\in {\cal M}_{g, n}$:
$$
\begin{array}{ccc}
{\rm B}_{  G, C, P}& \hra &{\rm B}_{  G, g,n}\\
\downarrow &&\downarrow \\
(C, P) & \hra &{\cal M}_{g, n}
\end{array}
$$
\paragraph{2. The Hitchin base and the spectral cover for $PGL_m$.} 
A point $t$ of the Hitchin's base of  $(C, P)$ is given by the following data:
$$
t = (C, P; t_2, t_3, \ldots t_m), \quad t_k \in \Omega_C(P)^{\otimes k}. 
$$
The associated spectral curve is a curve in $T^*C$ given by solutions of the polynomial equation:
\be \la{eqsc}
\Sigma_t:= 
\{\lambda \in T^*C~|~ \lambda^m + t_2\lambda^{m-2} + \ldots + t_{m-1}\lambda +t_m=0\} \subset T^*C. 
\ee
The projection $T^*C\to C$ induces the spectral cover 
$\pi_t:\Sigma_t \to C$. It is an $m:1$ ramified cover of $C$. 
For generic $t$ it has only simple ramification points. Their projections are 
precisely the zeros of the discriminant $\Delta_m= (\prod_{1 \leq i<j\leq m}(\lambda_i-\lambda_j))^2$ 
of   equation (\ref{eqsc}), where $\{\lambda_i\}_{i =1, ..., m}$ are the roots. 
For example, $\Delta_3 = 4t_2^3-27t_3^2$.  So $\Delta_m \in \Omega_C(P)^{m(m-1)}$. 
So the total number of ramification points is $m(m-1) (2g-2+n)$.

\paragraph{3. A family of open CY threefolds over the universal Hitchin base.} 
It is a family ${\rm Y}_{G, C, P}$ of open Calabi-Yau complex algebraic threefolds over  
 the Hitchin base ${\rm B}_{G, C, P}$:  
\be \la{2.19.12.1}
{\rm Y}_{G, C, P} \lra {\rm B}_{G, C, P}.
\ee

When the set of punctures is empty, 
it was studied  in \cite{DDDHP}, \cite{DDP}, \cite{G}. Then the
 intermediate  Jacobians of the fibers of fibration (\ref{2.19.12.1}) 
provide Hitchin's integrable system related to $G$ \cite{DDP}. 
A general construction for an arbitrary decorated surface, possibly 
with boundary components carrying special points, is given  \cite[Section 8]{KS3}.

When the pair $(C, P)$ varies over  ${\cal M}_{g,n}$, we get the  universal family of 
open  threefolds: 
\be \la{2.19.12.1ab}
{\rm Y}_{G, g,n} \lra {\rm B}_{G, g,n}.
\ee

\paragraph{Example.} Let $G=PGL_2$. 
Then  Hitchin's base parametrises pairs $(C, P, t)$, where  
$t$ is a quadratic differential on $C$ with the order two poles at the punctures.
Consider a threefold 
\be \la{4.4.12.1}
{\rm Y}^t_{C, P}:= \{x\in C-P; \alpha_1, \alpha_2, \alpha_3 \in \Omega^1_x(C) ~|~ \alpha_1^{\otimes 2}+
 \alpha_2^{\otimes 2}+ \alpha_3^{\otimes 2} =t\}. 
\ee
The four dimensional variety of the quadruples $(x; \alpha_1, \alpha_2, \alpha_3)$ 
has a canonical volume form. Indeed, the determinant of its cotangent space  is 
identified with $T_x(C)^{\otimes 3}\otimes T^*_x(C) = T_x(C)^{\otimes 2}$. 
The quadratic differential $t$ trivialises it. Therefore the hypersurface given by equation (\ref{4.4.12.1}) 
inherits a volume form. This way we get an open Calabi-Yau threefold.

\paragraph{4. The Hitchin integrable system for $G=PGL_m$.} 
A Higgs bundle on $(C, P)$ is a holomorphic $G$-bundle ${\cal E}$ on  $C$  
with an operator $\Phi: {\cal E} \lra {\cal E}\otimes\Omega_C(P)$. 
A Higgs bundle $({\cal E}, \Phi)$ on $C$ is the same thing as a coherent sheaf ${\cal L}_{\cal E}$ on the  surface 
$\Omega_C(P)$. Informally ${\cal L}_{\cal E}$ is given by the eigenvalues of the Higgs field operator 
$\Phi: {\cal E} \lra {\cal E}\otimes\Omega_C(P)$. Let $\pi: \Omega_C(P)\lra C$ be the canonical projection. 
Then $\pi_{*}{\cal L}_{\cal E} = {\cal E}$. 
Denote by 
${\cal H}_{G, g,n}$ the moduli space of Higgs $G$-bundles over the universal curve 
$(C, P)$ over ${\cal M}_{g, n}$. Then 
there is the Hitchin integrable system 
$$
p: {\cal H}_{G, g,n} \lra 
{\rm B}_{G, g,n}~~~~({\cal E}, \Phi)\lms ({\rm Tr}\Phi^2, {\rm Tr}\Phi^3, ..., {\rm Tr}\Phi^m).
$$

Let ${\rm B}^\times_{G, g,n}$ be the complement to the discriminant, that is the open part of the universal 
Hitchin base formed by  the points over which the map $p$ is non-singular. Then the fiber  
is an Abelian variety $A_t$.  
It is described as follows. The spectral cover  
$\pi_t: \Sigma_t \lra C$ gives rise to a map of  Jacobians $\pi_{t*}: J_{\Sigma_t} \lra J_C$. One has 
$$
A_t = {\rm Ker}(\pi_{t*}: J_{\Sigma_t} \lra J_C).
$$
There is a local system of lattices 
 over  ${\rm B}^\times_{  G, g,n}$ with the fibers 
$$
\Lambda_t:= H_1(A_t, \Z) = {\rm Ker}(H_1(\Sigma_t, \Z)  \lra H_1(C, \Z)). 
$$
The space ${\rm B}^\times_{  G, g,n}$ is fibered over the moduli space ${\cal M}_{g,n}$ 
with a fiber ${\rm B}^\times_{  G, C, P}$, the complement to the discriminant in the Hitchin base.  
The fundamental group $\pi_1({\rm B}^\times_{  G, g,n}, t)$ 
is an extension 
$$
1\lra \pi_1({\rm B}^\times_{  G, C, P}) \lra \pi_1({\rm B}^\times_{  G, g,n}, t) \lra \pi_1({\cal M}_{g,n}, (C,P))\lra 1.
$$
The fundamental group $\pi_1({\rm B}^\times_{  G, g,n}, t)$ 
acts on the lattice $\Lambda_t$. 

\paragraph{5. Framed Hitchin base.} We also need a $W:1$ cover  $\widetilde {\cal B}_{G, S}$ of the universal Hitching base:
 $$
 \widetilde {\cal B}_{G, S}  \lra  \widetilde {\rm B}_{G, S}. 
 $$
 When $G = PGL_m$, a point of $\widetilde {\cal B}_{G, S}$ is given by a point $t \in \widetilde {\rm B}_{G, S}$ plus for each puncture $p$ on $S$ a choice  
 of an ordering of the points of the spectral curve over the point $p$. Equivalently, it is given by an ordering of the roots  
 $\lambda_1, ..., \lambda_m$ of the polynomial (\ref{eqsc}) near the point $p$. The $\lambda_i$'s are 1-forms. So we get  an ordered set of their residues 
 $c_i:= {\rm Res}_p\lambda_i$ at $p$.  
 
  We call the $\widetilde {\cal B}_{G, S}$ the framed universal Hitchin base. Dealing with the moduli space ${\cal X}_{G, S}$ of framed 
  $G$-local systems on $S$ we always deal with the framed universal Hitchin base   $\widetilde {\cal B}_{G, S}$.

  \subsection{${\rm CY}_3$ categories from ideal webs and open Calabi-Yau threefolds} \la{SSec2.4}

\paragraph{1. Geometry of the topological spectral cover.} 
Denote by ${S}$  
  a genus $g$ oriented topological surface with $n$ punctures. 
Let ${\cal W}$ be an ideal ${\rm A}_{m-1}$-web on ${S}$. 
Denote by $\Gamma_{\cal W}$ the ${\rm A}_{m-1}$-bipartite graph on  $S$ corresponding 
to  ${\cal W}$. 
The spectral 
surface $\Sigma_{\cal W}$ related to the 
bipartite graph $\Gamma_{\cal W}$ is a surface with $mn$ punctures. 
The spectral map is a ramified $m:1$ map
\be \la{PSMAP}
\pi_{{\cal W}}: \Sigma_{{\cal W}} \lra S.
\ee
It is unramified over punctures: there are $m$ punctures on  $\Sigma_{{\cal W}}$ over each puncture on $S$.

\bl
For a generic $t$, the spectral cover $\pi_t:\Sigma_t \to C$ has the same number of ramification points and the same genus 
as the topological 
spectral cover $\pi_{\cal W}: {\Sigma}_{\cal W} \to { S}$. 
\el

\begin{proof} 
The ramification points of the spectral map $\pi_{\cal W}$ are the 
 $\bullet$-vertices of the web 
${\cal W}$ on ${S}$ of valency $\geq 3$: a $\bullet$-vertex of valency $v$ 
has ramification index $v-2$. To count them, take the ideal ${\rm A}_{m-1}$-web 
corresponding to an ideal triangulation of $S$. 
Then the number of $3$-valent $\bullet$-vertices in a single triangle is ${m \choose 2}$, while the number of ideal triangles is 
$2(2g-2+n)$. So the number $3$-valent $\bullet$-vertices is  ${m \choose 2} 2(2g-2+n) = m(m-1)(2g-2+n)$. 
This matches the number of ramification points of the spectral map $\pi_t$, see Section \ref{SSEEC5.3}. 
The degree of each of the covers is $m$. So the  
surface $\Sigma_{\cal W}$ has the same genus and the same number of punctures, $mn$, 
as the  $\Sigma_{t}$. 
\end{proof}

\paragraph{2. A conjectural combinatorial description of Fukaya categories via ideal webs.} The fiber over a point $t\in {\rm B}^\times_{  G, g,n}$  of the complement to  
discriminant of fibration (\ref{2.19.12.1ab}) is a smooth
open Calabi-Yau threefold ${\cal Y}_t$. 
Its Fukaya category ${\cal F}({\cal Y}_t)$ 
is a Calabi-Yau category. We do not 
discuss its definition here. 
We arrive at a family of Calabi-Yau categories 
over ${\rm B}^\times_{  G, g,n}$. 
 The fundamental group  
$\pi_1({\rm B}^\times_{  G, g,n}, t)$ 
acts by symmetries of the Fukaya category ${\cal F}({\cal Y}_t)$. 
There is a map
\be \la{9.6.14.2}
 \pi_1({\rm B}^\times_{  G, g,n}, t) \lra \pi_1({\cal M}_{g,n}; (C, P)) = \Gamma_S.
\ee
So there are two kinds of 3d CY categories, coming with certain groups of 
symmetries:

\begin{itemize}

\item The combinatorial 
category ${\cal C}_{\cal W}$ assigned to an ideal ${\rm A}_{m-1}$-web ${\cal W}$ on $S$.  

The group $\Gamma_{G, S}$ from Definition \ref{Def5.15} acts by  its symmetries.

\item The Fukaya category ${\cal F}({\cal Y}_t)$ of a complex smooth open CY threefold ${\cal Y}_t$, for $G=PGL_m$, 
  over a generic point $t\in {\rm B}_{G, g,n}^\times$ of the universal Hitchin base.

The fundamental group $\pi_1({\rm B}_{G, g,n}^\times) $ acts  by its symmetries. 
\end{itemize}

We conjecture that the 
easily defined category ${\cal C}_{\cal W}$ provides a transparent combinatorial model 
of a full subcategory of  the Fukaya category.

\bcon
Let ${\cal W}$ be an ${\rm A}_{m-1}$-web  on $S$. Then there is a fully faithful functor 
$$
\varphi: {\cal C}_{\cal W} \lra {\cal F}({\cal Y}_t).
$$ 
It transforms  cluster collections in ${\cal C}_{\cal W}$ to cluster 
collections in the Fukaya category provided by special Lagrangian spheres.  
The functor $\varphi$ intertwines the action of the symmetry group 
$\Gamma_{G, S}$ with the action of $\pi_1({\rm B}_{G, g,n}^\times)$ on 
the Fukaya category, providing a commutative diagram
\be \la{9.6.14.10}
\begin{array}{ccc}
   \Gamma_{G, S} &\lra & \Gamma_{S}\\
\downarrow  &&\downarrow =\\
\pi_1({\rm B}_{G, g,n}^\times) & \lra &\pi_1({\cal M}_{g,n})\\
\end{array}
\ee
\econ

One should have a similar story when $S$ has special points on the boundary.

 
Let us summarise how the combinatorial data provided by the web 
should match the one on the Fukaya category side. 

\begin{center}\begin{tabular}{| c | c | c |}
\bf{Complex geometry:}&{\bf Topology:}\\
\hline 
\bf{Hitchin system for $PGL_m$}&{\bf Moduli space}\\
& {\bf of framed $PGL_m$-local systems}\\
\hline\hline  
A Riemann surface $C$ & A punctured oriented surface $S$\\ 
\hline
A domain ${\cal B}_{\cal W}$ in the complement to the &An ${\rm A}_{m-1}$-web ${\cal W}$ on $S$ = \\
discriminant of the universal Hitchin base & its bipartite graph $\Gamma_{\cal W}$\\
\hline
Spectral curve $\Sigma_t$, where $t\in {\cal B}_{\cal W}$, and & Spectral curve $\Sigma_{\cal W}$, and the\\
the spectral cover $\pi_t: \Sigma_t \lra C$ &  spectral cover $\pi_{\cal W}: \Sigma_{\cal W} \lra S$\\
\hline
 Ramification points of the spectral cover & $\bullet$-points of the graph $\Gamma_{\cal W}$\\
\hline
 Lattice 
$\Lambda_t:= {\rm Ker} (H_1(\Sigma_t;\Z ) \to H_1(C;\Z))$ & Lattice 
$\Lambda_{\cal W}:= {\rm Ker} (H_1(\Sigma_{\cal W};\Z ) \to H_1(S;\Z))$ \\
\hline
Gauss-Manin connection on lattices $\Lambda_t$ & Two by two moves of webs\\
\hline
Fukaya category of the open Calabi-Yau& 3d CY ${\rm A}_\infty$-category ${\cal C}_{\cal W}$ of the  \\
  threefold ${\cal Y}_t$ related to the Hitchin system&  quiver with potential asigned to  ${\cal W}$ \\\hline
A collection of special Lagrangian & A cluster collection of spherical objects \\
spheres generating Fukaya category & in ${\cal C}_{\cal W}$, 
parametrised by internal faces of $\Gamma_{\cal W}$\\\hline
\end{tabular}\end{center}

\section{A complex manifold   of cluster stability conditions }  

Below we construct a manifold whose fundamental group should realize the generalized braid group. 

\paragraph{1. The positive and negative tropical domains.} Let ${\cal X}$ be a cluster Poisson variety. A cluster   coordinate system 
$\{X^{\bf q}_i\}$  assigned to a quiver ${\bf q}$ gives rise to a tropicalized 
 cluster Poisson coordinate system $\{x^{\bf q}_i\}$ on the space of real tropical points ${\cal X}(\R^t)$. It  determines 
a  cone in  ${\cal X}(\R^t)$ consisting of the points 
with \underline{non negative} coordinates in the  coordinate system $\{x^{\bf q}_i\}$:
$$
{\cal X}^{+}_{\bf q}(\R^t):= \{l\in {\cal X}(\R^t)~|~ x^{\bf q}_i(l) \geq 0\} \subset {\cal X}(\R^t).
$$

  It was conjectured in \cite{FG4} and proved in \cite{GHKK} that when ${\bf q}$ runs through all quivers obtained from a given one by mutations, the  
 internal parts of the domains ${\cal X}^{+}_{\bf q}(\R^t)$ are disjoint. 
The domains  ${\cal X}^{+}_{\bf q}(\R^t)$  meet along the walls 
where one of the coordinates become zero. The wall $x^{\bf q}_k =0$ separates 
the domain ${\cal X}^+_{\bf q}(\R^t)$ and ${\cal X}^+_{\bf q'}(\R^t)$ where the quivers ${\bf q}$ and ${\bf q}'$ are 
related by a mutation in the direction $k$. 
The domains ${\cal X}^+_{\bf q}(\R^t)$ usually do not cover the space ${\cal X}(\R^t)$. 

For example, the space ${\cal X}_{PGL_2, T}(\R^t)$ for a punctured torus $T$ is illustrated in \cite[Fig 1]{FG4}. 
For any punctured surface $S$, the complement ${\cal X}_{PGL_2, S}(\R^t) - {\cal X}^+_{PGL_2, S}(\R^t)$ is of measure zero.

\bd The union of the domains ${\cal X}^+_{\bf q}(\R^t)$  is  called the {\it positive tropical domain}:  
$$
{\cal X}^+ (\R^t) \subset {\cal X}(\R^t).
$$
\ed
Evidently, it is a connected domain. 
Similarly, there is a {\it negative tropical domain}
$$
{\cal X}^- (\R^t) \subset {\cal X}(\R^t).
$$
It is the union of the domains
$ 
{\cal X}^{-}_{\bf q}(\R^t):= \{l\in {\cal X}(\R^t)~|~ x^{\bf q}_i(l) \leq  0\} \subset {\cal X}(\R^t).
$ 

The cluster modular group $\Gamma$ 
acts by automorphisms of  the domains ${\cal X}^\pm (\R^t)$.

 Given any mutation-connected set ${\cal S}$ of quivers, 
there are subdomains
\be \la{GLU}
{\cal X}^\pm_{\cal S}(\R^t):= \coprod_{{\bf q}\in {\cal S}}  {\cal X}^\pm_{\bf q}(\R^t) \subset {\cal X}(\R^t). 
\ee

\paragraph{2. The complex domain ${\cal U}_{\cal X}$ of cluster stability conditions.} Our goal is to define 
 a $\Gamma$-equivariant complex domain ${\cal U}_{\cal X}$ equipped with  a  $\Gamma$-equivariant surjective projection
$$
p: {\cal U}_{\cal X} \lra {\cal X}^+ (\R^t).
$$
It  
sits  in the space of stability conditions ${\rm Stab}_{\cal X}/{\rm Br}_{\cal X}$ for the 3d CY category related to the cluster variety ${\cal X}$ with a generic potential, i.e. there is a 
$\Gamma$-equivariant embedding:
$$
{\cal U}_{\cal X} \subset {\rm Stab}_{\cal X}/{\rm Br}_{\cal X}.
$$
So each point $u \in {\cal U}_{\cal X}$ has a neighborhood isomorphic to an open subset in ${\rm Hom}(\Lambda, \C)$. 

Take an upper half plane
$$
{\cal H} := \{r{\rm exp}(i\varphi) ~|~ r>0, 0 < \varphi \leq \pi\}.
$$

A quiver ${\bf q} = (\Lambda, \{e_i\}_{i\in {\rm I}}, (\ast, \ast))$ 
 gives rise to a complex domain ${\cal U}_{\bf q}$   parametrizing homomorphisms of abelian groups  
$f: \Lambda \to \C$  (the central charges)  with $f(e_i) \in {\cal H}$:
\be \la{DSIL}
{\cal U}_{\bf q} :=\{f: \Lambda \lra \C ~|~ f(a+b) = f(a) + f(b), ~f(e_i) \in {\cal H}, ~\forall i \in {\rm I}\}.
\ee
So the domain ${\cal U}_{\bf q}$ is the product of the domains ${\cal H}$ over the set $I$ of vertices of  ${\bf q}$:
\be \la{domD1}
{\cal U}_{\bf q} := \prod_{j\in {\rm I}} {\cal H}_{j},
\qquad {\cal H}_j := \{z_j= r_j{\rm exp}(i\varphi_j) ~|~ r_j>0, 0 < \varphi_j \leq \pi\}.
\ee

Recall that we assign to a quiver ${\bf q}$ a set of {tropical cluster Poisson coordinates} $\{x_i\}_{i\in {\rm I}}$. 
The imaginary part map 
${\rm Im}: {\cal H}   \lra \R_{\geq 0}$,  $z
  \lms  {\rm Im}(z)$ provides  a  projection 
\be \la{4.26.16.1}
\begin{split}
&{\rm Im}: {\cal U}_{\bf q} \lra {\cal X}^+_{\bf q}(\R^t), ~~~~ \{z_i\} \lms \{x_i:= {\rm Im}(z_i)\}.\\
\end{split}
\ee
 So we assign to a quiver ${\bf q}$ a triple 
 \be \la{triple}
 \Bigl({\cal U}_{\bf q}, ~{\cal X}^+_{\bf q}(\R^t), ~{\rm Im}: {\cal U}_{\bf q} \lra {\cal X}^+_{\bf q}(\R^t)\Bigr).
 \ee 
 An isomorphism of  quivers ${\bf q}\lra {\bf q}'$ gives rise to a unique isomorphism of the triples. 

Let us glue the triples (\ref{triple}) assigned to the quivers obtained by mutations of a quiver ${\bf q}$. 

{\it Mutations}. Denote by $\overline {\cal U}_{{\bf q}}$ the partial closure of the domain ${\cal U}_{{\bf q}}$ 
defined by using the inequalities $ 0 \leq \varphi_j \leq \pi$ in (\ref{domD1}) while keeping $r>0$. The domain $\overline {\cal U}_{{\bf q}}$ 
has boundary walls of real codimension one given by setting either $\varphi_j = \pi$ or $\varphi_j = 0$ for some $j \in {\rm I}$. 

Given  a mutation of quivers ${\bf q} \stackrel{k}{\longrightarrow} {\bf q}'$ 
  in the direction $k$,   
let us glue the  domain $\overline{\cal U}_{{\bf q}}$ to the domain  $\overline{\cal U}_{{\bf q}'}$ 
along  the two  real codimension one walls $\varphi_k=\pi$ and $\varphi_k=0$:
 \be
\begin{split}
& G_{{\bf q}\to {\bf q'}}: \{\mbox{the wall $\varphi_{k}=\pi$ in $\overline{\cal U}_{{\bf q}}$}\} ~\stackrel{\rm glued}{\longrightarrow} ~\{\mbox{the wall 
$\varphi'_{k}=0$ in $\overline{\cal U}_{{\bf q'}}$}\},\\
 &G_{{\bf q}\to {\bf q'}}: \{\mbox{the wall $\varphi_{k}=0$ in $\overline{\cal U}_{{\bf q}}$}\} ~ \stackrel{\rm glued}{\longrightarrow} ~\{\mbox{the wall $\varphi'_{k}=\pi$ in $\overline{\cal U}_{{\bf q'}}$}\}.\\
\end{split}
\ee
Namely, let $\{z_i\}_{i\in {\rm I}}$ be the coordinates in the domain ${\cal U}_{{\bf q}}$, and $\{z'_i\}_{i\in {\rm I}}$  the ones in ${\cal U}_{{\bf q'}}$. We define  
the gluing maps as follows:  
\begin{equation} \label{5.11.03.6z}
\begin{split}
& z_i':= (G_{{\bf q}\to {\bf q'}})^*(z_{i}) := \left\{ \begin{array}{lll} 
- z_k & \mbox{ if $i = k$}, \\ 
z_i+ [\varepsilon_{ik}]_+ z_k& \mbox{ if $i \not = k$}. \\
\end{array}\right.\\
\end{split}
\end{equation} 
Observe that if ${\rm Im}(z_k)=0$, and ${ \rm Im}(z_i) >0$ for $i \not = k$, then ${\rm Im}(z'_k)=0$ and ${\rm  Im}(z'_i) >0$ for $i \not = k$. 
Furthermore, if $\varphi_k=0$ then $\varphi_k'=\pi$, and vice versa,    if $\varphi_k=\pi$ then $\varphi_k'=0$.

Here is an equivalent way to look at the gluing of the domains $\overline{\cal U}_{{\bf q}}$ and  $\overline{\cal U}_{{\bf q}'}$.

 Recall that the basis $\{e_i'\}$ of the  quiver ${\bf q}'$ is related to the basis $\{e_i\}$ of the quiver ${\bf q}$  by 
\begin{equation} \label{5.11.03.6z2}
\begin{split}
& e_i':= \left\{ \begin{array}{lll} 
- e_k & \mbox{ if $i = k$}, \\ 
e_i+ [\varepsilon_{ik}]_+ e_k& \mbox{ if $i \not = k$}. \\
\end{array}\right.\\
\end{split}
\end{equation} 
Both domains ${\cal U}_{{\bf q}}$ and ${\cal U}_{{\bf q'}}$ sit in the space ${\rm Hom}(\Lambda, \C)$ via (\ref{DSIL}). It is clear from (\ref{5.11.03.6z2}) that they are glued 
by the map (\ref{5.11.03.6z}) precisely the way they intersect there; the result of  gluing is identified  with  the union of  domains 
$$
{\cal U}_{{\bf q}} \cup {\cal U}_{{\bf q'}} \subset {\rm Hom}(\Lambda, \C).
$$
In particular the  ${\cal U}_{{\bf q}} \cup {\cal U}_{{\bf q'}}$ is a manifold.

Denote  by $\{x'_i\}$ the tropical cluster Poisson coordinates assigned  to ${\bf q'}$.    They are related to the ones  $\{x_i\}_{i\in {\rm I}}$ for the quiver ${\bf q}$ by the formula
\begin{equation} \label{5.11.03.1xtr}
\mu_k^*(x'_{i}) :=\left\{\begin{array}{ll} -x_k& \mbox{if }  i=k \\
 x_i-\varepsilon_{ik}{\rm min}\{0, -{\rm sgn} (\varepsilon_{ik})x_k\} & \mbox{\rm otherwise}.
\end{array} \right.
\end{equation}

  \bl \la{4.30.16.1}
  The projection (\ref{4.26.16.1}) is compatible with mutations. So we get a projection
 \be \la{CANPRO}
  {\cal U}_{{\bf q}} \cup {\cal U}_{{\bf q'}} \lra   {\cal X}^+_{{\bf q}}(\R^t) \cup {\cal X}^+_{{\bf q'}}(\R^t).    
  \ee
  \el

\begin{proof} If $x_k\geq 0$, the  mutation  formula (\ref{5.11.03.1xtr}) for  tropical cluster Poisson coordinates reduces to:
\begin{equation} \label{5.11.03.1xtrasd}
\mu_k^*(x'_{i}) :=\left\{\begin{array}{ll} -x_k& \mbox{if }  i=k \\
 x_i + [\varepsilon_{ik}]_+ x_k & \mbox{\rm otherwise}.
\end{array} \right.
\end{equation}
So taking the imaginary part of the gluing map  (\ref{5.11.03.6z}) we recover (\ref{5.11.03.1xtrasd}). 
\end{proof}
        The result of gluing of the two domains corresponding to a pair of quivers ${\bf q} \stackrel{k}{\longleftrightarrow} {\bf q}'$ is a domain which is 
homotopy equivalent to a circle. Indeed, it is  identified with 
$$
 {\cal U}_{{\bf q}} \cup  {\cal U}_{{\bf q'}} = \C^* \times \prod_{j \not = k}{\cal H}_j.
$$         
         
Consider now the composition of two mutations in the same direction $k$:
$$
{\bf q} \stackrel{k}{\lra} {\bf q}' \stackrel{k}{\lra} {\bf q}''.
$$
It is described by  mutations of the  basis $\{e_i\}$ defining the quiver ${\bf q}$:
$$
\{e_i\} \stackrel{k}{\lra} \{e'_i\}  \stackrel{k}{\lra}  \{e''_i\} .
$$
The basis $\{e''_i\}$ is usually different then the one $\{e_i\}$. However quivers ${\bf q''}$ and ${\bf q}$ are canonically isomorphic: 
there is a unique isomorphism of lattices $\varphi: \Lambda \lra \Lambda$  
preserving the form $(\ast, \ast)$ such that $\varphi(e_i)= e''_i$. It defines an isomorphism of quivers, and hence domains:
$$
\varphi: {\cal U}_{\bf q} \lra {\cal U}_{\bf q''}.
$$
It is well known that the cones ${\cal X}_{\bf q}^+(\R^t)$ and ${\cal X}_{\bf q''}^+(\R^t)$  coincide. 
So we get a commutative diagram
  \begin{displaymath}
    \xymatrix{
        {\cal U}_{\bf q}\ar[r]^{\varphi}  \ar[d]^{{\rm Im}_{\bf q}}  & {\cal U}_{\bf q''} 
          \ar[d]^{{\rm Im}_{\bf q''}} \\
         {\cal X}_{\bf q}^+(\R^t)   \ar[r]^{=} &{\cal X}_{\bf q''}^+(\R^t) }
         \end{displaymath}

         {\it Cluster transformations.} Consider a tree ${\rm T}_n$ whose  edges incident 
to a given vertex are 
parametrized by the set 
${\rm I}$.  
We assign to each vertex ${a}$ of  ${\rm T}_n$ a quiver ${\bf q}_a$   
{considered up to an isomorphism}, so that  the  quivers ${\bf q}_a$ and ${\bf q}_b$ assigned to 
the vertices of an edge   labeled by an element $k \in {\rm I}$ 
 are related by the mutation $\mu_k$: 
\be \la{muasse}
{\bf q}_b = \mu_{k}({\bf q}_a), ~~~~ {\bf q}_a = \mu_{k}({\bf q}_b). 
\ee
 We assign to each vertex  $a$ of the tree ${\rm T}_n$   
 the domain ${\cal U}_{{\bf q}_a}$.  
Gluing the pairs of domains ${\cal U}_{{\bf q}_a}$ and ${\cal U}_{{\bf q}_b}$ assigned to each edge of the tree ${\rm T}_n$  as   explained above,  
we get a manifold denoted $\widetilde{\cal U}_{\cal X}$.  Lemma \ref{4.30.16.1} implies that  projections (\ref{4.26.16.1}) are glued into a projection
\be \la{CANPRO1}
{\rm Im}: \widetilde{\cal U}_{\cal X} \lra {\cal X}(\R^t).
\ee
For any pair of vertices $a, b$ of the tree ${\rm T}_n$, there is a unique  path ${\bf i}$ on the tree   connecting them:
\be \la{CPATH}
{\bf i}:  ~~a= a_0 \stackrel{k_1}{\lra} a_1 \stackrel{k_2}{\lra}  \ldots \stackrel{k_{m-1}}{\lra} a_{m-1} \stackrel{k_{m}}{\lra}  a_m=b.
\ee
Moreover, 
 there is a unique automorphism of trees $i_{a\to b}: {\rm T}_n \lra {\rm T}_n$ which sends $a$ to $b$ and preserves the decoration of the edges by the set ${\rm I}$. 
 
 Now take a pair of vertices $a$ and $b$ of the tree ${\rm T}_n$ such that

\begin{itemize}

\item The  quivers ${\bf q}_a$ and $ {\bf q}_b$    are isomorphic.
\end{itemize}
Then  for any vertex $c$, the automorphism  $i_{a\to b}: {\rm T}_n \lra {\rm T}_n$ induces an isomorphism of the quivers assigned to the vertices 
$c$ and $i_{a\to b}(c)$. Therefore it induces an automorphism 
$$
\varphi_{a\to b}: \widetilde{\cal U}_{\cal X} \lra \widetilde{\cal U}_{\cal X}.
$$
On the other hand, it induces an automorphism of the cluster Poisson variety: 
\be \la{AUTPHI}
\Phi_{a\to b}: {\cal X} \lra {\cal X}.
\ee
Its tropicalisation   induces an automorphism of the space of real tropical points:
$$
\Phi^t_{a\to b}: {\cal X}(\R^t) \lra {\cal X}(\R^t).
$$
We have a commutative diagram:
 \begin{displaymath}
    \xymatrix{
        \widetilde {\cal U}_{\cal X}  \ar[r]^{{\varphi}_{{a}\to {b}}}  \ar[d]^{{\rm Im} }  & \widetilde {\cal U}_{\cal X} 
          \ar[d]^{{\rm Im}} \\
         {\cal X}^+(\R^t) \ar[r]^{\Phi^t_{a\to b}} &  {\cal X}^+(\R^t) }
         \end{displaymath}

 Now take a pair of vertices $a$ and $b$ of the tree ${\rm T}_n$ such that 
 \begin{itemize} \item The   quivers ${\bf q}_a$ and $ {\bf q}_b$ are isomorphic.   The map     $\Phi^t_{a\to b}$ is the identity.  
\end{itemize}
The corresponding automorphisms ${\varphi}_{{a}\to {b}}$ of the manifold $ \widetilde{\cal U}_{\cal X}$ 
form a group denoted by $\Gamma_{\rm triv}$.   

\bd
The   orbifold ${\cal U}_{\cal X}$ is the quotient of the manifold $ \widetilde{\cal U}_{\cal X}$ by the action of the group $\Gamma_{\rm triv}$:
$$
{\cal U}_{\cal X}:= \widetilde{\cal U}_{\cal X}/\Gamma_{\rm triv}.
$$
\ed
Then, by the  definition,   the projection (\ref{CANPRO1}) induces a canonical projection
\be \la{IM}
{\rm Im}: 
{\cal U}_{\cal X} \lra {\cal X}^+(\R^t).
\ee

By the   definition, the automorphisms (\ref{AUTPHI}) generate the cluster modular group $\Gamma$. 
Evidently, the projection (\ref{IM}) is $\Gamma$-equivariant. 
The complex manifold   ${\cal U}_{\cal X}$ has a $\Gamma$-equivariant  Poisson structure $\{\ast, \ast\}$ which in any cluster coordinate system $\{z_i\}$ is given by 
$$
\{z_i, z_j\} = (e_i, e_j). 
$$
The real tropical space also has a Poisson structure given by the same formula $\{x_i, x_j\} = (e_i, e_j)$ in any tropical cluster coordinate system $\{x_i\}$. 
The projection (\ref{IM}) is evidently Poisson.

 We summarize the main result of this Section.
 
 \bt \la{MTHU}
 Let ${\cal X}$ be a cluster Poisson variety. Then the domains ${\cal U}_{\bf q}$  are glued into a complex Poisson orbifold ${\cal U}_{\cal X}$. 
 The cluster modular group $\Gamma$ acts on ${\cal U}_{\cal X}$. 
  The projections (\ref{4.26.16.1}) give rise to  a   $\Gamma$-equivariant Poisson projection on the space of   positive real tropical points of ${\cal X}$:
 $$
 p: {\cal U}_{\cal X} \lra {\cal X}^+(\R^t).
 $$
 There is a $\Gamma$-equivariant  embedding of the orbifold ${\cal U}_{\cal X}$   into the space of   stability conditions ${\rm Stab}_{\cal X}/{\rm Br}_{\cal X}$ on 
 the 3d Calaby-Yau category related to a generic potential 
 of a quiver ${\bf q}$:
\be \la{DSCon}
 {\cal U}_{\cal X} \subset   {\rm Stab}_{\cal X}/{\rm Br}_{\cal X}.
 \ee
  \et
  
  \begin{proof}  The only claim which was not discussed yet is the $\Gamma$-equivariance of the projection $p$. 
  Let $\widetilde \Gamma$ be the group of automorphisms $i_{a\to b}$ of the tree. Denote by $\Gamma'_{\rm triv}$ the subgroup 
  of the automorphisms $i_{a\to b}$ which induce the trivial automorphisms $\Phi_{a\to b}$ of the cluster variety.    
  The cluster modular group is given by $\Gamma :=  \widetilde \Gamma/  \Gamma'_{\rm triv}$.  However the group which acts naturally 
  on $\widetilde {\cal U}_{\cal X}$ is $ \widetilde \Gamma/  \Gamma_{\rm triv}$. Evidently $\Gamma'_{\rm triv} \subset \Gamma_{\rm triv}$: 
  if a cluster Poisson transformation   is trivial, its tropicalisation is also trivial. The crucial claim that in fact $\Gamma'_{\rm triv} = \Gamma_{\rm triv}$    
  is given by Theorem \ref{basiclamDT}. 
  \bt \la{basiclamDT} 
Let $\sigma: {\bf q} \to {\bf q'}$ be a cluster transformation which induces the identity transformation on the set of the positive tropical points. 
Then the quivers ${\bf q}$ and ${\bf q'}$ are isomorphic, and the 
 corresponding quantum cluster transformation $\Phi(\sigma)$  is the identity. 
\et

Theorem   \ref{basiclamDT}  is derived  in \cite[Section 2]{GS16} from a  theorem  of Keller  \cite{K11}, which was stated  
  using 
 {\it $c$-vectors} and {\it $C$-matrices}  
 rather than the tropical points of 
cluster Poisson varieties.  
  \end{proof}

   There is a similar    complex Poisson manifold ${\cal U}^-_{\cal X}$ equipped  with a $\Gamma$-equivariant projection  
 $$
  p: {\cal U}^-_{\cal X} \lra {\cal X}^-(\R^t). 
  $$  
  
  \paragraph{An example.} Consider the cluster Poisson variety ${\cal X}_{\rm A_2}$ of type ${\rm A}_2$. Its tropicalisation ${\cal X}_{\rm A_2}(\R^t)$ 
  is identified with the real plane $\R^2$ decomposed into the union of five cones, see Figure \ref{A2}, denoted by $C_1, ..., C_5$, numbered clockwise. 
  We assign to each cone $C_i$ a complex domain ${\cal U}_i$, isomorphic to ${\cal H}\times {\cal H} $. Then we glue these domains 
  in the cyclic order:
  $$
  {\cal U}_1 \longleftrightarrow {\cal U}_2 \longleftrightarrow{\cal U}_3 \longleftrightarrow {\cal U}_4 \longleftrightarrow {\cal U}_5 \longleftrightarrow {\cal U}_1.
  $$  
  The mapping class group is isomorphic to $\Z/5\Z$. Its generator acts by shifting $ {\cal U}_i$ to ${\cal U}_{i+1}$.

\end{document}